\documentclass[12pt,a4paper]{article} 

\def\endofps{EndOfTheIncludedPostscriptMagicCookie}
\chardef\other=12
\newwrite\psdumphandle 
\outer\def\psdump#1{\par\medbreak
  \immediate\openout\psdumphandle=#1
  \copytoblankline}
\def\copytoblankline{\begingroup\setupcopy\copypsline}
\def\setupcopy{\def\do##1{\catcode`##1=\other}\dospecials
  \catcode`\\=\other \obeylines}
{\obeylines \gdef\copypsline#1
  {\def\next{#1}%
  \ifx\next\endofps\let\next=\endgroup %
  \else\immediate\write\psdumphandle{\next} \let\next=\copypsline\fi\next}}
\outer\def\closepsdump{
  \immediate\closeout\psdumphandle}

\psdump{constr.eps}
statusdict begin/waittimeout 0 def end{/MGXPS_2.8_EPS begin}stopped{/MGXPS_2.8_EPS
200 dict def}if MGXPS_2.8_EPS begin/bd{bind def}bind def/xd{exch def}bd/ld{load
def}bd/M/moveto ld/m/rmoveto ld/L/lineto ld/l/rlineto ld/w/setlinewidth
ld/n/newpath ld/P/closepath ld/tr/translate ld/gs/gsave ld/gr/grestore
ld/lc/setlinecap ld/lj/setlinejoin ld/ml/setmiterlimit ld/TM/concat
ld/F/fill ld/bz/curveto ld/_m matrix def/_cm matrix def/_d 0 def/_o
0 def/_pb 0 def/_hb -1 def currentscreen/_P xd dup/_A xd/__A xd dup/_F
xd/__F xd/srgbc{_neg 0 ne{3{neg 1 add 3 1 roll}repeat}if setrgbcolor}bd/ssc{gs
0 0 M 0 _yp neg _xp 0 0 _yp 3 D closepath fill gr 0 0 0 srgbc}bd/Cne{/_bs
xd/_neg where{/_neg get/_neg xd pop}{pop 0/_neg xd}ifelse dup 0 ne{setflat}{pop}ifelse
dup -1 eq{pop}{pop pop}ifelse/_rz xd/_yp xd/_xp xd/_a 0 def 0 ne{/_a
90 def -90 rotate _xp _rz div 72 mul neg 0 tr}if 0 _yp _rz div 72 mul
tr 72 _rz div dup neg scale/#copies xd 2 copy/_ym xd/_xm xd tr n 1
lj 1 lc 1 w 0 8{dup}repeat sc fc bc}bd/clpg{n M dup neg 3 1 roll 0
l 0 exch l 0 l P clip n}bd/SB{_o 0 ne _d 0 ne and{gs []0 setdash _br
_bg _bb srgbc stroke gr}if}bd/S{_r _g _b srgbc SB stroke}bd/SX{_r _g
_b srgbc _cm currentmatrix pop _m setmatrix SB stroke _cm setmatrix}bd/d{/_d
xd setdash}bd/c{3{255 div 3 1 roll}repeat}bd/sc{c/_b xd/_g xd/_r xd}bd/bc{c/_bb
xd/_bg xd/_br xd}bd/fc{c/_B xd/_G xd/_R xd}bd/tc{c srgbc}bd/f{0 eq{eofill}{fill}ifelse}bd/D{{l}repeat}bd/PT{dup
statusdict exch known{cvx exec}{pop}ifelse}bd/TT{dup statusdict exch
known{cvx exec}{pop pop}ifelse}bd/SM{_m currentmatrix pop}bd/RM{_m
setmatrix}bd/sp{dup 0 eq{/_pb xd}{255 div/_gr xd/_pb xd 100 div/_f
xd}ifelse}bd/B{/_t2 xd/_t3 xd _pb _t2 _t3 8 idiv add get 1 7 _t3 8
mod sub bitshift and 0 ne}bd/SF{2{1 add 4 mul cvi exch}repeat B{0}{1}ifelse}bd/pb{_gr
setgray _f _a/SF load setscreen}bd
/s{100 div/_t xd[_t 0 0 _t 1 _t sub _xp 2 div mul 1 _t sub _yp 2 div
mul]concat}bd
/_am matrix def/a{_am currentmatrix pop/_t1 xd/_t2 xd/_t3 exch 10
div def/_t4 exch 10 div def tr _t2 _t1 neg scale/_t xd _t 1 eq{0 0
M}if 0 0 1 _t4 _t3 arc _t 0 ne{P}if _am setmatrix}bd/an{_am currentmatrix
pop/_t1 xd/_t2 xd/_t3 exch 10 div def/_t4 exch 10 div def tr _t2 _t1
neg scale/_t xd _t 1 eq{0 0 M}if 0 0 1 _t4 _t3 arcn _t 0 ne{P}if _am
setmatrix}bd
/_em matrix def/E{_em currentmatrix pop/_t1 xd/_t2 xd/_t3 xd/_t4 xd
_t4 _t3 _t1 sub M _t4 _t3 tr _t2 _t1 neg scale 0 0 1 90 450 arc _em
setmatrix}bd
/_rm matrix def/cn{_rm currentmatrix pop tr x3 y3 scale arc _rm setmatrix}bd/rr{2
div/y3 xd 2 div/x3 xd/y2 xd/x2 xd/y1 xd/x1 xd x1 y1 y3 add M 0 0 1
180 270 x1 x3 add y1 y3 add cn x2 x3 sub y1 L 0 0 1 270 0 x2 x3 sub
y1 y3 add cn x2 y2 y3 sub L 0 0 1 0 90 x2 x3 sub y2 y3 sub cn x1 x3
add y2 L 0 0 1 90 180 x1 x3 add y2 y3 sub cn}bd
/cp{/_t xd{{gs _t 0 eq{eoclip}{clip}ifelse}stopped{gr currentflat
1 add setflat}{exit}ifelse}loop n}bd/p{P dup 0 eq{S pop}{_R _G _B srgbc
_pb 0 eq{2 eq{gs f gr S}{f}ifelse}{_pb 1 eq{_o 1 eq{gs _br _bg _bb
srgbc 1 index f gr}if _R _G _B srgbc 2 eq{gs hf gr S}{hf}ifelse n}{2
eq{true pbc S}{false pbc n}ifelse}ifelse}ifelse}ifelse}bd/px{P dup
0 eq{SX pop}{_R _G _B srgbc _pb 0 eq{2 eq{gs f gr SX}{f}ifelse}{_pb
1 eq{_o 1 eq{gs _br _bg _bb srgbc 1 index f gr}if 2 eq{gs _m setmatrix
hf gr SX}{gs _m setmatrix hf gr}ifelse n}{2 eq{true pbc SX}{false pbc
n}ifelse}ifelse}ifelse}ifelse}bd
/spc{dup 0 eq{/_pb xd}{1 eq{255 div/_gr xd/_pb xd/_pg 0 def/_pr 0
def}{/_pb xd/_pg xd/_pr xd}ifelse 1 add/_pq xd}ifelse}bd/_im{[_t 0
0 _t 0 0]}bd/bbx{_pq 75 mul _rz div/_t xd pathbbox 1 add cvi 31 or/_ury
xd 1 add cvi 31 or/_urx xd 1 sub cvi -32 and/_lly xd 1 sub cvi -32
and/_llx xd/_dH _ury _lly sub def/_dW _urx _llx sub def}bd/xp{/_pbs
xd _row _pbs mul _pbs getinterval/_s xd 0 _pbs 8 _pbs div dup _sW add
1 sub exch div cvi dup/str exch string def _pbs sub{str exch _s putinterval}for
str}bd/pbc{/__t save def{_m setmatrix}if 0 eq{eoclip}{clip}ifelse/_row
0 def bbx _llx _lly tr _t _dW mul ceiling cvi 8 add dup 8 mod sub dup/_sW
xd _t _dH mul ceiling cvi 8 add dup 8 mod sub _pg 0 ne{8 _im{_pr 8
xp}{_pg 8 xp}{_pb 8 xp/_row _row 1 add 8 mod def}true 3 colorimage}{1
_o eq{gs _br _bg _bb srgbc 2 copy true _im{_pb 1 xp/_row _row 1 add
8 mod def}imagemask gr}if _R _G _B srgbc false _im{_pb 1 xp/_row _row
1 add 8 mod def}imagemask}ifelse __t restore}bd
/hp{dup/_hb xd -1 eq{/_pb 0 def}{1 add/_pq xd/_pb 1 def}ifelse}bd/vert{X0
_w X1{dup Y0 M Y1 L stroke}for}bd/horz{Y0 _w Y1{dup X0 exch M X1 exch
L stroke}for}bd/fdiag{X0 _w X1{Y0 M X1 X0 sub dup l stroke}for Y0 _w
Y1{X0 exch M Y1 Y0 sub dup l stroke}for}bd/bdiag{X0 _w X1{Y1 M X1 X0
sub dup neg l stroke}for Y0 _w Y1{X0 exch M Y1 Y0 sub dup neg l stroke}for}bd/AU{1
add cvi 31 or}bd/AD{1 sub cvi -32 and}bd/hf{pathbbox 1 add cvi 31 or/Y1
xd 1 add cvi 31 or/X1 xd 1 sub cvi -32 and/Y0 xd 1 sub cvi -32 and/X0
xd 2 w [] 0 setdash/_w _rz 20 div 8 div 2 mul _pq div round 8 mul def
cp _hb 0 eq{horz}if _hb 1 eq{vert}if _hb 2 eq{fdiag}if _hb 3 eq{bdiag}if
_hb 4 eq{horz vert}if _hb 5 eq{fdiag bdiag}if gr}bd
/sp{dup 0 eq{/_pb xd}{255 div/_gr xd/_pb xd 100 div/_f xd}ifelse}bd/B{/_t2
xd/_t3 xd _pb _t2 _t3 8 idiv add get 1 7 _t3 8 mod sub bitshift and
0 ne}bd/SF{2{1 add 4 mul cvi exch}repeat B{0}{1}ifelse}bd/pb{_gr setgray
_f _a/SF load setscreen}bd
/SCA{10 div dup dup/_a xd/_A xd/__na xd currentscreen/__p xd/__a xd/__f
xd __f __na/__p load setscreen}bd
/SCF{100 div dup dup/_F xd/__nf xd currentscreen/__p xd/__a xd/__f
xd __nf __a/__p load setscreen}bd
/makeOblique 0 def/bC 1 def/makeoutlinedict 7 dict def/mbf{makeoutlinedict
begin/uniqueid exch def/strokewidth exch def/newfontname exch def/basefontname
exch def/basefontdict basefontname findfont def/numentries basefontdict
maxlength 1 add def basefontdict/UniqueID known not{/numentries numentries
1 add def}if/outfontdict numentries dict def basefontdict{exch dup/FID
ne{exch outfontdict 3 1 roll put}{pop pop}ifelse}forall outfontdict/FontName
newfontname put outfontdict/PaintType 2 put outfontdict/StrokeWidth
strokewidth put outfontdict/UniqueID uniqueid put newfontname outfontdict
definefont pop end}bd/usf{findfont exch makeOblique 0 ne{/makeOblique
makeOblique sin neg def [1 0 makeOblique 1 0 0] matrix concatmatrix/makeOblique
0 def}if makefont setfont}bd/sf{/_nf xd FontDirectory _nf known{pop}{findfont
dup maxlength dict/_nfd xd{exch dup/FID ne{exch _nfd 3 1 roll put}{pop
pop}ifelse}forall _nfd dup/FontName _nf put/Encoding WE put _nf _nfd
definefont pop}ifelse _nf usf}bd/us{gs [] 0 setdash 0 setlinecap w
0 exch M _w _be add 0 l currentrgbcolor sC S gr}bd/ob{n/_h xd _cx _cy
M _w 0 l 0 _h l _w neg 0 l P tc fill}bd/st{dup/_be xd 0 ne{/_bc xd}if/_vt
xd dup stringwidth pop/_sw xd dup length dup 1 gt{1 sub}if/_sl xd _be
0 eq{_w _sw sub _sl div 0 _vt{exch}if 3 -1 roll ashow}{_be _bc div
0 _vt{exch}if 32 _w _sw sub _sl div 0 _vt{exch}if 6 -1 roll awidthshow}ifelse}bd/sb{gs
M 10 div rotate/_w xd currentpoint/_cy xd/_cx xd 0 ne{ob}if tc _cx
_cy 3 -1 roll add tr 0 0 M st dup 4 and 4 eq{4 sub/_dx _cx _px sub
def _dx neg 0 tr/_w _w _dx add def}if dup 0 eq{pop}{dup 1 eq{pop us}{2
eq{us}{dup 3 1 roll us us}ifelse}ifelse}ifelse _cx/_px xd gr}bd
/cl{gs n M dup neg 3 1 roll 0 l 0 exch l 0 l P clip n}bd/clx{gs n
_cm currentmatrix pop _m setmatrix n M dup neg 3 1 roll 0 l 0 exch
l 0 l P clip _cm setmatrix n}bd
/sbb{/_t save def _neg 0 eq{0}{1}ifelse setgray tr/_dW xd/_dH xd dup/_sW
xd 8 _bs div dup _sW add 1 sub exch div cvi/str exch string def exch
dup/_sH xd 3 -1 roll{/_im{[_sW _dW div 0 0 _sH neg _dH div 0 _sH]}bd}{/_im{[_sW
_dW div 0 0 _sH _dH div 0 0]}bd}ifelse _bs 1 eq{false}{_bs}ifelse _im{currentfile
str readhexstring pop}_bs 1 eq{imagemask}{image}ifelse _t restore}bd
/scbb{/_t save def tr/_dW xd/_dH xd dup/_sW xd 8 _bs div dup _sW add
1 sub exch div cvi 3 mul/str exch string def exch dup/_sH xd _bs 4
-1 roll{/_im{[_sW _dW div 0 0 _sH neg _dH div 0 _sH]}bd}{/_im{[_sW
_dW div 0 0 _sH _dH div 0 0]}bd}ifelse _im{currentfile str readhexstring
pop}false 3 colorimage _t restore}bd
/mr{2 copy -1 eq{_yp _ym 2 mul sub}{0}ifelse exch -1 eq{_xp _xm 2
mul sub}{0}ifelse exch tr scale}bd
/WE 256 array def StandardEncoding WE copy pop WE dup
0[/grave/acute/circumflex/tilde/macron/breve/dotaccent/dieresis/ring/cedilla
/hungarumlaut/ogonek/caron/dotlessi
]putinterval dup
39/quotesingle put dup  96/grave put dup 124/bar put dup
145[/quoteleft/quoteright/quotedblleft/quotedblright/bullet/endash/emdash]putinterval 
160[/space/exclamdown/cent/sterling/currency/yen/brokenbar/section/dieresis
/copyright/ordfeminine/guillemotleft/logicalnot/endash/registered/macron/ring
/plusminus/twosuperior/threesuperior/acute/mu/paragraph/periodcentered
/cedilla/onesuperior/ordmasculine/guillemotright/onequarter/onehalf
/threequarters/questiondown/Agrave/Aacute/Acircumflex/Atilde/Adieresis/Aring
/AE/Ccedilla/Egrave/Eacute/Ecircumflex/Edieresis/Igrave/Iacute/Icircumflex
/Idieresis/Eth/Ntilde/Ograve/Oacute/Ocircumflex/Otilde/Odieresis/circumflex
/Oslash/Ugrave/Uacute/Ucircumflex/Udieresis/Yacute/Thorn/germandbls/agrave
/aacute/acircumflex/atilde/adieresis/aring/ae/ccedilla/egrave/eacute
/ecircumflex/edieresis/igrave/iacute/icircumflex/idieresis/eth/ntilde/ograve
/oacute/ocircumflex/otilde/odieresis/tilde/oslash/ugrave/uacute/ucircumflex
/udieresis/yacute/thorn/ydieresis]putinterval
/TTEV 256 array def TTEV
0 [/G00/G01/G02/G03/G04/G05/G06/G07/G08/G09/G0a/G0b/G0c/G0d/G0e/G0f
   /G10/G11/G12/G13/G14/G15/G16/G17/G18/G19/G1a/G1b/G1c/G1d/G1e/G1f
   /G20/G21/G22/G23/G24/G25/G26/G27/G28/G29/G2a/G2b/G2c/G2d/G2e/G2f
   /G30/G31/G32/G33/G34/G35/G36/G37/G38/G39/G3a/G3b/G3c/G3d/G3e/G3f
   /G40/G41/G42/G43/G44/G45/G46/G47/G48/G49/G4a/G4b/G4c/G4d/G4e/G4f
   /G50/G51/G52/G53/G54/G55/G56/G57/G58/G59/G5a/G5b/G5c/G5d/G5e/G5f
   /G60/G61/G62/G63/G64/G65/G66/G67/G68/G69/G6a/G6b/G6c/G6d/G6e/G6f
   /G70/G71/G72/G73/G74/G75/G76/G77/G78/G79/G7a/G7b/G7c/G7d/G7e/G7f
   /G80/G81/G82/G83/G84/G85/G86/G87/G88/G89/G8a/G8b/G8c/G8d/G8e/G8f
   /G90/G91/G92/G93/G94/G95/G96/G97/G98/G99/G9a/G9b/G9c/G9d/G9e/G9f
   /Ga0/Ga1/Ga2/Ga3/Ga4/Ga5/Ga6/Ga7/Ga8/Ga9/Gaa/Gab/Gac/Gad/Gae/Gaf
   /Gb0/Gb1/Gb2/Gb3/Gb4/Gb5/Gb6/Gb7/Gb8/Gb9/Gba/Gbb/Gbc/Gbd/Gbe/Gbf
   /Gc0/Gc1/Gc2/Gc3/Gc4/Gc5/Gc6/Gc7/Gc8/Gc9/Gca/Gcb/Gcc/Gcd/Gce/Gcf
   /Gd0/Gd1/Gd2/Gd3/Gd4/Gd5/Gd6/Gd7/Gd8/Gd9/Gda/Gdb/Gdc/Gdd/Gde/Gdf
   /Ge0/Ge1/Ge2/Ge3/Ge4/Ge5/Ge6/Ge7/Ge8/Ge9/Gea/Geb/Gec/Ged/Gee/Gef
   /Gf0/Gf1/Gf2/Gf3/Gf4/Gf5/Gf6/Gf7/Gf8/Gf9/Gfa/Gfb/Gfc/Gfd/Gfe/Gff
] putinterval
end
MGXPS_2.8_EPS begin
0 0 1 0 3307 4677 400 -1 0 0 1 Cne
4551 3134 0 0 clpg
0 0 0 sc
/_o 0 def
255 255 255 bc
1009 139 M 1131 104 1149 104 1271 139 bz
S
1019 144 M 0 0 -11 -4 2 D
S
1012 128 M 0 0 -4 12 2 D
S
1007 232 M 1129 267 1147 267 1269 232 bz
S
1265 243 M 0 0 6 -12 2 D
S
1258 226 M 0 0 12 5 2 D
S
0 230 M 213 125 395 84 554 107 bz
713 129 842 213 948 363 bz
S
165 74 M 228 183 303 232 390 221 bz
476 207 622 -18 684 5 bz
744 31 769 148 757 368 bz
S
6 w
10433 1000 div ml
758 178 M 764 261 760 302 756 365 bz
S
194 150 M 367 94 568 66 758 177 bz
S
1 w
1368 186 M 0 0 927 0 2 D
S
1558 96 M 1613 253 1679 324 1759 314 bz
1837 301 1962 23 2033 32 bz
2103 42 2152 155 2179 374 bz
S
6 w
2143 186 M 2164 242 2175 323 2179 374 bz
S
1561 186 M 0 0 580 0 2 D
S
end
EndOfTheIncludedPostscriptMagicCookie

\closepsdump

\psdump{coord.eps}
statusdict begin/waittimeout 0 def end{/MGXPS_2.8_EPS begin}stopped{/MGXPS_2.8_EPS
200 dict def}if MGXPS_2.8_EPS begin/bd{bind def}bind def/xd{exch def}bd/ld{load
def}bd/M/moveto ld/m/rmoveto ld/L/lineto ld/l/rlineto ld/w/setlinewidth
ld/n/newpath ld/P/closepath ld/tr/translate ld/gs/gsave ld/gr/grestore
ld/lc/setlinecap ld/lj/setlinejoin ld/ml/setmiterlimit ld/TM/concat
ld/F/fill ld/bz/curveto ld/_m matrix def/_cm matrix def/_d 0 def/_o
0 def/_pb 0 def/_hb -1 def currentscreen/_P xd dup/_A xd/__A xd dup/_F
xd/__F xd/srgbc{_neg 0 ne{3{neg 1 add 3 1 roll}repeat}if setrgbcolor}bd/ssc{gs
0 0 M 0 _yp neg _xp 0 0 _yp 3 D closepath fill gr 0 0 0 srgbc}bd/Cne{/_bs
xd/_neg where{/_neg get/_neg xd pop}{pop 0/_neg xd}ifelse dup 0 ne{setflat}{pop}ifelse
dup -1 eq{pop}{pop pop}ifelse/_rz xd/_yp xd/_xp xd/_a 0 def 0 ne{/_a
90 def -90 rotate _xp _rz div 72 mul neg 0 tr}if 0 _yp _rz div 72 mul
tr 72 _rz div dup neg scale/#copies xd 2 copy/_ym xd/_xm xd tr n 1
lj 1 lc 1 w 0 8{dup}repeat sc fc bc}bd/clpg{n M dup neg 3 1 roll 0
l 0 exch l 0 l P clip n}bd/SB{_o 0 ne _d 0 ne and{gs []0 setdash _br
_bg _bb srgbc stroke gr}if}bd/S{_r _g _b srgbc SB stroke}bd/SX{_r _g
_b srgbc _cm currentmatrix pop _m setmatrix SB stroke _cm setmatrix}bd/d{/_d
xd setdash}bd/c{3{255 div 3 1 roll}repeat}bd/sc{c/_b xd/_g xd/_r xd}bd/bc{c/_bb
xd/_bg xd/_br xd}bd/fc{c/_B xd/_G xd/_R xd}bd/tc{c srgbc}bd/f{0 eq{eofill}{fill}ifelse}bd/D{{l}repeat}bd/PT{dup
statusdict exch known{cvx exec}{pop}ifelse}bd/TT{dup statusdict exch
known{cvx exec}{pop pop}ifelse}bd/SM{_m currentmatrix pop}bd/RM{_m
setmatrix}bd/sp{dup 0 eq{/_pb xd}{255 div/_gr xd/_pb xd 100 div/_f
xd}ifelse}bd/B{/_t2 xd/_t3 xd _pb _t2 _t3 8 idiv add get 1 7 _t3 8
mod sub bitshift and 0 ne}bd/SF{2{1 add 4 mul cvi exch}repeat B{0}{1}ifelse}bd/pb{_gr
setgray _f _a/SF load setscreen}bd
/s{100 div/_t xd[_t 0 0 _t 1 _t sub _xp 2 div mul 1 _t sub _yp 2 div
mul]concat}bd
/_am matrix def/a{_am currentmatrix pop/_t1 xd/_t2 xd/_t3 exch 10
div def/_t4 exch 10 div def tr _t2 _t1 neg scale/_t xd _t 1 eq{0 0
M}if 0 0 1 _t4 _t3 arc _t 0 ne{P}if _am setmatrix}bd/an{_am currentmatrix
pop/_t1 xd/_t2 xd/_t3 exch 10 div def/_t4 exch 10 div def tr _t2 _t1
neg scale/_t xd _t 1 eq{0 0 M}if 0 0 1 _t4 _t3 arcn _t 0 ne{P}if _am
setmatrix}bd
/_em matrix def/E{_em currentmatrix pop/_t1 xd/_t2 xd/_t3 xd/_t4 xd
_t4 _t3 _t1 sub M _t4 _t3 tr _t2 _t1 neg scale 0 0 1 90 450 arc _em
setmatrix}bd
/_rm matrix def/cn{_rm currentmatrix pop tr x3 y3 scale arc _rm setmatrix}bd/rr{2
div/y3 xd 2 div/x3 xd/y2 xd/x2 xd/y1 xd/x1 xd x1 y1 y3 add M 0 0 1
180 270 x1 x3 add y1 y3 add cn x2 x3 sub y1 L 0 0 1 270 0 x2 x3 sub
y1 y3 add cn x2 y2 y3 sub L 0 0 1 0 90 x2 x3 sub y2 y3 sub cn x1 x3
add y2 L 0 0 1 90 180 x1 x3 add y2 y3 sub cn}bd
/cp{/_t xd{{gs _t 0 eq{eoclip}{clip}ifelse}stopped{gr currentflat
1 add setflat}{exit}ifelse}loop n}bd/p{P dup 0 eq{S pop}{_R _G _B srgbc
_pb 0 eq{2 eq{gs f gr S}{f}ifelse}{_pb 1 eq{_o 1 eq{gs _br _bg _bb
srgbc 1 index f gr}if _R _G _B srgbc 2 eq{gs hf gr S}{hf}ifelse n}{2
eq{true pbc S}{false pbc n}ifelse}ifelse}ifelse}ifelse}bd/px{P dup
0 eq{SX pop}{_R _G _B srgbc _pb 0 eq{2 eq{gs f gr SX}{f}ifelse}{_pb
1 eq{_o 1 eq{gs _br _bg _bb srgbc 1 index f gr}if 2 eq{gs _m setmatrix
hf gr SX}{gs _m setmatrix hf gr}ifelse n}{2 eq{true pbc SX}{false pbc
n}ifelse}ifelse}ifelse}ifelse}bd
/spc{dup 0 eq{/_pb xd}{1 eq{255 div/_gr xd/_pb xd/_pg 0 def/_pr 0
def}{/_pb xd/_pg xd/_pr xd}ifelse 1 add/_pq xd}ifelse}bd/_im{[_t 0
0 _t 0 0]}bd/bbx{_pq 75 mul _rz div/_t xd pathbbox 1 add cvi 31 or/_ury
xd 1 add cvi 31 or/_urx xd 1 sub cvi -32 and/_lly xd 1 sub cvi -32
and/_llx xd/_dH _ury _lly sub def/_dW _urx _llx sub def}bd/xp{/_pbs
xd _row _pbs mul _pbs getinterval/_s xd 0 _pbs 8 _pbs div dup _sW add
1 sub exch div cvi dup/str exch string def _pbs sub{str exch _s putinterval}for
str}bd/pbc{/__t save def{_m setmatrix}if 0 eq{eoclip}{clip}ifelse/_row
0 def bbx _llx _lly tr _t _dW mul ceiling cvi 8 add dup 8 mod sub dup/_sW
xd _t _dH mul ceiling cvi 8 add dup 8 mod sub _pg 0 ne{8 _im{_pr 8
xp}{_pg 8 xp}{_pb 8 xp/_row _row 1 add 8 mod def}true 3 colorimage}{1
_o eq{gs _br _bg _bb srgbc 2 copy true _im{_pb 1 xp/_row _row 1 add
8 mod def}imagemask gr}if _R _G _B srgbc false _im{_pb 1 xp/_row _row
1 add 8 mod def}imagemask}ifelse __t restore}bd
/hp{dup/_hb xd -1 eq{/_pb 0 def}{1 add/_pq xd/_pb 1 def}ifelse}bd/vert{X0
_w X1{dup Y0 M Y1 L stroke}for}bd/horz{Y0 _w Y1{dup X0 exch M X1 exch
L stroke}for}bd/fdiag{X0 _w X1{Y0 M X1 X0 sub dup l stroke}for Y0 _w
Y1{X0 exch M Y1 Y0 sub dup l stroke}for}bd/bdiag{X0 _w X1{Y1 M X1 X0
sub dup neg l stroke}for Y0 _w Y1{X0 exch M Y1 Y0 sub dup neg l stroke}for}bd/AU{1
add cvi 31 or}bd/AD{1 sub cvi -32 and}bd/hf{pathbbox 1 add cvi 31 or/Y1
xd 1 add cvi 31 or/X1 xd 1 sub cvi -32 and/Y0 xd 1 sub cvi -32 and/X0
xd 2 w [] 0 setdash/_w _rz 20 div 8 div 2 mul _pq div round 8 mul def
cp _hb 0 eq{horz}if _hb 1 eq{vert}if _hb 2 eq{fdiag}if _hb 3 eq{bdiag}if
_hb 4 eq{horz vert}if _hb 5 eq{fdiag bdiag}if gr}bd
/sp{dup 0 eq{/_pb xd}{255 div/_gr xd/_pb xd 100 div/_f xd}ifelse}bd/B{/_t2
xd/_t3 xd _pb _t2 _t3 8 idiv add get 1 7 _t3 8 mod sub bitshift and
0 ne}bd/SF{2{1 add 4 mul cvi exch}repeat B{0}{1}ifelse}bd/pb{_gr setgray
_f _a/SF load setscreen}bd
/SCA{10 div dup dup/_a xd/_A xd/__na xd currentscreen/__p xd/__a xd/__f
xd __f __na/__p load setscreen}bd
/SCF{100 div dup dup/_F xd/__nf xd currentscreen/__p xd/__a xd/__f
xd __nf __a/__p load setscreen}bd
/makeOblique 0 def/bC 1 def/makeoutlinedict 7 dict def/mbf{makeoutlinedict
begin/uniqueid exch def/strokewidth exch def/newfontname exch def/basefontname
exch def/basefontdict basefontname findfont def/numentries basefontdict
maxlength 1 add def basefontdict/UniqueID known not{/numentries numentries
1 add def}if/outfontdict numentries dict def basefontdict{exch dup/FID
ne{exch outfontdict 3 1 roll put}{pop pop}ifelse}forall outfontdict/FontName
newfontname put outfontdict/PaintType 2 put outfontdict/StrokeWidth
strokewidth put outfontdict/UniqueID uniqueid put newfontname outfontdict
definefont pop end}bd/usf{findfont exch makeOblique 0 ne{/makeOblique
makeOblique sin neg def [1 0 makeOblique 1 0 0] matrix concatmatrix/makeOblique
0 def}if makefont setfont}bd/sf{/_nf xd FontDirectory _nf known{pop}{findfont
dup maxlength dict/_nfd xd{exch dup/FID ne{exch _nfd 3 1 roll put}{pop
pop}ifelse}forall _nfd dup/FontName _nf put/Encoding WE put _nf _nfd
definefont pop}ifelse _nf usf}bd/us{gs [] 0 setdash 0 setlinecap w
0 exch M _w _be add 0 l currentrgbcolor sC S gr}bd/ob{n/_h xd _cx _cy
M _w 0 l 0 _h l _w neg 0 l P tc fill}bd/st{dup/_be xd 0 ne{/_bc xd}if/_vt
xd dup stringwidth pop/_sw xd dup length dup 1 gt{1 sub}if/_sl xd _be
0 eq{_w _sw sub _sl div 0 _vt{exch}if 3 -1 roll ashow}{_be _bc div
0 _vt{exch}if 32 _w _sw sub _sl div 0 _vt{exch}if 6 -1 roll awidthshow}ifelse}bd/sb{gs
M 10 div rotate/_w xd currentpoint/_cy xd/_cx xd 0 ne{ob}if tc _cx
_cy 3 -1 roll add tr 0 0 M st dup 4 and 4 eq{4 sub/_dx _cx _px sub
def _dx neg 0 tr/_w _w _dx add def}if dup 0 eq{pop}{dup 1 eq{pop us}{2
eq{us}{dup 3 1 roll us us}ifelse}ifelse}ifelse _cx/_px xd gr}bd
/cl{gs n M dup neg 3 1 roll 0 l 0 exch l 0 l P clip n}bd/clx{gs n
_cm currentmatrix pop _m setmatrix n M dup neg 3 1 roll 0 l 0 exch
l 0 l P clip _cm setmatrix n}bd
/sbb{/_t save def _neg 0 eq{0}{1}ifelse setgray tr/_dW xd/_dH xd dup/_sW
xd 8 _bs div dup _sW add 1 sub exch div cvi/str exch string def exch
dup/_sH xd 3 -1 roll{/_im{[_sW _dW div 0 0 _sH neg _dH div 0 _sH]}bd}{/_im{[_sW
_dW div 0 0 _sH _dH div 0 0]}bd}ifelse _bs 1 eq{false}{_bs}ifelse _im{currentfile
str readhexstring pop}_bs 1 eq{imagemask}{image}ifelse _t restore}bd
/scbb{/_t save def tr/_dW xd/_dH xd dup/_sW xd 8 _bs div dup _sW add
1 sub exch div cvi 3 mul/str exch string def exch dup/_sH xd _bs 4
-1 roll{/_im{[_sW _dW div 0 0 _sH neg _dH div 0 _sH]}bd}{/_im{[_sW
_dW div 0 0 _sH _dH div 0 0]}bd}ifelse _im{currentfile str readhexstring
pop}false 3 colorimage _t restore}bd
/mr{2 copy -1 eq{_yp _ym 2 mul sub}{0}ifelse exch -1 eq{_xp _xm 2
mul sub}{0}ifelse exch tr scale}bd
/WE 256 array def StandardEncoding WE copy pop WE dup
0[/grave/acute/circumflex/tilde/macron/breve/dotaccent/dieresis/ring/cedilla
/hungarumlaut/ogonek/caron/dotlessi
]putinterval dup
39/quotesingle put dup  96/grave put dup 124/bar put dup
145[/quoteleft/quoteright/quotedblleft/quotedblright/bullet/endash/emdash]putinterval 
160[/space/exclamdown/cent/sterling/currency/yen/brokenbar/section/dieresis
/copyright/ordfeminine/guillemotleft/logicalnot/endash/registered/macron/ring
/plusminus/twosuperior/threesuperior/acute/mu/paragraph/periodcentered
/cedilla/onesuperior/ordmasculine/guillemotright/onequarter/onehalf
/threequarters/questiondown/Agrave/Aacute/Acircumflex/Atilde/Adieresis/Aring
/AE/Ccedilla/Egrave/Eacute/Ecircumflex/Edieresis/Igrave/Iacute/Icircumflex
/Idieresis/Eth/Ntilde/Ograve/Oacute/Ocircumflex/Otilde/Odieresis/circumflex
/Oslash/Ugrave/Uacute/Ucircumflex/Udieresis/Yacute/Thorn/germandbls/agrave
/aacute/acircumflex/atilde/adieresis/aring/ae/ccedilla/egrave/eacute
/ecircumflex/edieresis/igrave/iacute/icircumflex/idieresis/eth/ntilde/ograve
/oacute/ocircumflex/otilde/odieresis/tilde/oslash/ugrave/uacute/ucircumflex
/udieresis/yacute/thorn/ydieresis]putinterval
/TTEV 256 array def TTEV
0 [/G00/G01/G02/G03/G04/G05/G06/G07/G08/G09/G0a/G0b/G0c/G0d/G0e/G0f
   /G10/G11/G12/G13/G14/G15/G16/G17/G18/G19/G1a/G1b/G1c/G1d/G1e/G1f
   /G20/G21/G22/G23/G24/G25/G26/G27/G28/G29/G2a/G2b/G2c/G2d/G2e/G2f
   /G30/G31/G32/G33/G34/G35/G36/G37/G38/G39/G3a/G3b/G3c/G3d/G3e/G3f
   /G40/G41/G42/G43/G44/G45/G46/G47/G48/G49/G4a/G4b/G4c/G4d/G4e/G4f
   /G50/G51/G52/G53/G54/G55/G56/G57/G58/G59/G5a/G5b/G5c/G5d/G5e/G5f
   /G60/G61/G62/G63/G64/G65/G66/G67/G68/G69/G6a/G6b/G6c/G6d/G6e/G6f
   /G70/G71/G72/G73/G74/G75/G76/G77/G78/G79/G7a/G7b/G7c/G7d/G7e/G7f
   /G80/G81/G82/G83/G84/G85/G86/G87/G88/G89/G8a/G8b/G8c/G8d/G8e/G8f
   /G90/G91/G92/G93/G94/G95/G96/G97/G98/G99/G9a/G9b/G9c/G9d/G9e/G9f
   /Ga0/Ga1/Ga2/Ga3/Ga4/Ga5/Ga6/Ga7/Ga8/Ga9/Gaa/Gab/Gac/Gad/Gae/Gaf
   /Gb0/Gb1/Gb2/Gb3/Gb4/Gb5/Gb6/Gb7/Gb8/Gb9/Gba/Gbb/Gbc/Gbd/Gbe/Gbf
   /Gc0/Gc1/Gc2/Gc3/Gc4/Gc5/Gc6/Gc7/Gc8/Gc9/Gca/Gcb/Gcc/Gcd/Gce/Gcf
   /Gd0/Gd1/Gd2/Gd3/Gd4/Gd5/Gd6/Gd7/Gd8/Gd9/Gda/Gdb/Gdc/Gdd/Gde/Gdf
   /Ge0/Ge1/Ge2/Ge3/Ge4/Ge5/Ge6/Ge7/Ge8/Ge9/Gea/Geb/Gec/Ged/Gee/Gef
   /Gf0/Gf1/Gf2/Gf3/Gf4/Gf5/Gf6/Gf7/Gf8/Gf9/Gfa/Gfb/Gfc/Gfd/Gfe/Gff
] putinterval
end
MGXPS_2.8_EPS begin
0 0 1 0 3400 4400 400 -1 0 0 1 Cne
4280 3240 0 0 clpg
0 0 0 sc
/_o 0 def
255 255 255 bc
854 0 M 1044 216 L
1002 264 998 270 1002 311 bz
1007 347 1047 403 1077 370 bz
1151 424 1183 441 1219 462 bz
1212 532 1228 607 1235 665 bz
1240 716 1219 769 1199 779 bz
1168 793 788 814 676 799 bz
565 782 567 783 526 775 bz
580 708 748 522 726 503 bz
700 480 521 635 496 769 bz
452 768 322 712 283 679 bz
234 639 221 612 307 529 bz
377 457 600 282 663 216 bz
720 166 854 0 854 0 bz
192 192 192 fc
0 2 p
516 349 M 607 282 690 274 765 279 bz
841 286 911 357 961 399 bz
1007 440 1019 463 1040 481 bz
1062 495 1074 499 1110 472 bz
1110 472 1169 530 1222 573 bz
1176 634 1106 643 1055 646 bz
989 646 926 611 878 576 bz
829 538 775 453 731 447 bz
686 441 647 466 601 489 bz
516 349 L
255 255 255 fc
0 2 p
277 679 M 208 777 98 824 53 832 bz
7 838 -9 797 4 712 bz
S
1216 459 M 1268 376 1309 310 1335 163 bz
S
6 w
10433 1000 div ml
278 679 M 401 492 615 366 738 366 bz
849 369 920 529 1011 558 bz
1102 584 1167 549 1216 460 bz
S
1 w
617 392 7 7 E
0 0 0 fc
0 2 p
729 322 M 0 0 -140 -159 2 D
S
end
EndOfTheIncludedPostscriptMagicCookie

\closepsdump

\psdump{distfn.eps}
statusdict begin/waittimeout 0 def end{/MGXPS_2.8_EPS begin}stopped{/MGXPS_2.8_EPS
200 dict def}if MGXPS_2.8_EPS begin/bd{bind def}bind def/xd{exch def}bd/ld{load
def}bd/M/moveto ld/m/rmoveto ld/L/lineto ld/l/rlineto ld/w/setlinewidth
ld/n/newpath ld/P/closepath ld/tr/translate ld/gs/gsave ld/gr/grestore
ld/lc/setlinecap ld/lj/setlinejoin ld/ml/setmiterlimit ld/TM/concat
ld/F/fill ld/bz/curveto ld/_m matrix def/_cm matrix def/_d 0 def/_o
0 def/_pb 0 def/_hb -1 def currentscreen/_P xd dup/_A xd/__A xd dup/_F
xd/__F xd/srgbc{_neg 0 ne{3{neg 1 add 3 1 roll}repeat}if setrgbcolor}bd/ssc{gs
0 0 M 0 _yp neg _xp 0 0 _yp 3 D closepath fill gr 0 0 0 srgbc}bd/Cne{/_bs
xd/_neg where{/_neg get/_neg xd pop}{pop 0/_neg xd}ifelse dup 0 ne{setflat}{pop}ifelse
dup -1 eq{pop}{pop pop}ifelse/_rz xd/_yp xd/_xp xd/_a 0 def 0 ne{/_a
90 def -90 rotate _xp _rz div 72 mul neg 0 tr}if 0 _yp _rz div 72 mul
tr 72 _rz div dup neg scale/#copies xd 2 copy/_ym xd/_xm xd tr n 1
lj 1 lc 1 w 0 8{dup}repeat sc fc bc}bd/clpg{n M dup neg 3 1 roll 0
l 0 exch l 0 l P clip n}bd/SB{_o 0 ne _d 0 ne and{gs []0 setdash _br
_bg _bb srgbc stroke gr}if}bd/S{_r _g _b srgbc SB stroke}bd/SX{_r _g
_b srgbc _cm currentmatrix pop _m setmatrix SB stroke _cm setmatrix}bd/d{/_d
xd setdash}bd/c{3{255 div 3 1 roll}repeat}bd/sc{c/_b xd/_g xd/_r xd}bd/bc{c/_bb
xd/_bg xd/_br xd}bd/fc{c/_B xd/_G xd/_R xd}bd/tc{c srgbc}bd/f{0 eq{eofill}{fill}ifelse}bd/D{{l}repeat}bd/PT{dup
statusdict exch known{cvx exec}{pop}ifelse}bd/TT{dup statusdict exch
known{cvx exec}{pop pop}ifelse}bd/SM{_m currentmatrix pop}bd/RM{_m
setmatrix}bd/sp{dup 0 eq{/_pb xd}{255 div/_gr xd/_pb xd 100 div/_f
xd}ifelse}bd/B{/_t2 xd/_t3 xd _pb _t2 _t3 8 idiv add get 1 7 _t3 8
mod sub bitshift and 0 ne}bd/SF{2{1 add 4 mul cvi exch}repeat B{0}{1}ifelse}bd/pb{_gr
setgray _f _a/SF load setscreen}bd
/s{100 div/_t xd[_t 0 0 _t 1 _t sub _xp 2 div mul 1 _t sub _yp 2 div
mul]concat}bd
/_am matrix def/a{_am currentmatrix pop/_t1 xd/_t2 xd/_t3 exch 10
div def/_t4 exch 10 div def tr _t2 _t1 neg scale/_t xd _t 1 eq{0 0
M}if 0 0 1 _t4 _t3 arc _t 0 ne{P}if _am setmatrix}bd/an{_am currentmatrix
pop/_t1 xd/_t2 xd/_t3 exch 10 div def/_t4 exch 10 div def tr _t2 _t1
neg scale/_t xd _t 1 eq{0 0 M}if 0 0 1 _t4 _t3 arcn _t 0 ne{P}if _am
setmatrix}bd
/_em matrix def/E{_em currentmatrix pop/_t1 xd/_t2 xd/_t3 xd/_t4 xd
_t4 _t3 _t1 sub M _t4 _t3 tr _t2 _t1 neg scale 0 0 1 90 450 arc _em
setmatrix}bd
/_rm matrix def/cn{_rm currentmatrix pop tr x3 y3 scale arc _rm setmatrix}bd/rr{2
div/y3 xd 2 div/x3 xd/y2 xd/x2 xd/y1 xd/x1 xd x1 y1 y3 add M 0 0 1
180 270 x1 x3 add y1 y3 add cn x2 x3 sub y1 L 0 0 1 270 0 x2 x3 sub
y1 y3 add cn x2 y2 y3 sub L 0 0 1 0 90 x2 x3 sub y2 y3 sub cn x1 x3
add y2 L 0 0 1 90 180 x1 x3 add y2 y3 sub cn}bd
/cp{/_t xd{{gs _t 0 eq{eoclip}{clip}ifelse}stopped{gr currentflat
1 add setflat}{exit}ifelse}loop n}bd/p{P dup 0 eq{S pop}{_R _G _B srgbc
_pb 0 eq{2 eq{gs f gr S}{f}ifelse}{_pb 1 eq{_o 1 eq{gs _br _bg _bb
srgbc 1 index f gr}if _R _G _B srgbc 2 eq{gs hf gr S}{hf}ifelse n}{2
eq{true pbc S}{false pbc n}ifelse}ifelse}ifelse}ifelse}bd/px{P dup
0 eq{SX pop}{_R _G _B srgbc _pb 0 eq{2 eq{gs f gr SX}{f}ifelse}{_pb
1 eq{_o 1 eq{gs _br _bg _bb srgbc 1 index f gr}if 2 eq{gs _m setmatrix
hf gr SX}{gs _m setmatrix hf gr}ifelse n}{2 eq{true pbc SX}{false pbc
n}ifelse}ifelse}ifelse}ifelse}bd
/spc{dup 0 eq{/_pb xd}{1 eq{255 div/_gr xd/_pb xd/_pg 0 def/_pr 0
def}{/_pb xd/_pg xd/_pr xd}ifelse 1 add/_pq xd}ifelse}bd/_im{[_t 0
0 _t 0 0]}bd/bbx{_pq 75 mul _rz div/_t xd pathbbox 1 add cvi 31 or/_ury
xd 1 add cvi 31 or/_urx xd 1 sub cvi -32 and/_lly xd 1 sub cvi -32
and/_llx xd/_dH _ury _lly sub def/_dW _urx _llx sub def}bd/xp{/_pbs
xd _row _pbs mul _pbs getinterval/_s xd 0 _pbs 8 _pbs div dup _sW add
1 sub exch div cvi dup/str exch string def _pbs sub{str exch _s putinterval}for
str}bd/pbc{/__t save def{_m setmatrix}if 0 eq{eoclip}{clip}ifelse/_row
0 def bbx _llx _lly tr _t _dW mul ceiling cvi 8 add dup 8 mod sub dup/_sW
xd _t _dH mul ceiling cvi 8 add dup 8 mod sub _pg 0 ne{8 _im{_pr 8
xp}{_pg 8 xp}{_pb 8 xp/_row _row 1 add 8 mod def}true 3 colorimage}{1
_o eq{gs _br _bg _bb srgbc 2 copy true _im{_pb 1 xp/_row _row 1 add
8 mod def}imagemask gr}if _R _G _B srgbc false _im{_pb 1 xp/_row _row
1 add 8 mod def}imagemask}ifelse __t restore}bd
/hp{dup/_hb xd -1 eq{/_pb 0 def}{1 add/_pq xd/_pb 1 def}ifelse}bd/vert{X0
_w X1{dup Y0 M Y1 L stroke}for}bd/horz{Y0 _w Y1{dup X0 exch M X1 exch
L stroke}for}bd/fdiag{X0 _w X1{Y0 M X1 X0 sub dup l stroke}for Y0 _w
Y1{X0 exch M Y1 Y0 sub dup l stroke}for}bd/bdiag{X0 _w X1{Y1 M X1 X0
sub dup neg l stroke}for Y0 _w Y1{X0 exch M Y1 Y0 sub dup neg l stroke}for}bd/AU{1
add cvi 31 or}bd/AD{1 sub cvi -32 and}bd/hf{pathbbox 1 add cvi 31 or/Y1
xd 1 add cvi 31 or/X1 xd 1 sub cvi -32 and/Y0 xd 1 sub cvi -32 and/X0
xd 2 w [] 0 setdash/_w _rz 20 div 8 div 2 mul _pq div round 8 mul def
cp _hb 0 eq{horz}if _hb 1 eq{vert}if _hb 2 eq{fdiag}if _hb 3 eq{bdiag}if
_hb 4 eq{horz vert}if _hb 5 eq{fdiag bdiag}if gr}bd
/sp{dup 0 eq{/_pb xd}{255 div/_gr xd/_pb xd 100 div/_f xd}ifelse}bd/B{/_t2
xd/_t3 xd _pb _t2 _t3 8 idiv add get 1 7 _t3 8 mod sub bitshift and
0 ne}bd/SF{2{1 add 4 mul cvi exch}repeat B{0}{1}ifelse}bd/pb{_gr setgray
_f _a/SF load setscreen}bd
/SCA{10 div dup dup/_a xd/_A xd/__na xd currentscreen/__p xd/__a xd/__f
xd __f __na/__p load setscreen}bd
/SCF{100 div dup dup/_F xd/__nf xd currentscreen/__p xd/__a xd/__f
xd __nf __a/__p load setscreen}bd
/makeOblique 0 def/bC 1 def/makeoutlinedict 7 dict def/mbf{makeoutlinedict
begin/uniqueid exch def/strokewidth exch def/newfontname exch def/basefontname
exch def/basefontdict basefontname findfont def/numentries basefontdict
maxlength 1 add def basefontdict/UniqueID known not{/numentries numentries
1 add def}if/outfontdict numentries dict def basefontdict{exch dup/FID
ne{exch outfontdict 3 1 roll put}{pop pop}ifelse}forall outfontdict/FontName
newfontname put outfontdict/PaintType 2 put outfontdict/StrokeWidth
strokewidth put outfontdict/UniqueID uniqueid put newfontname outfontdict
definefont pop end}bd/usf{findfont exch makeOblique 0 ne{/makeOblique
makeOblique sin neg def [1 0 makeOblique 1 0 0] matrix concatmatrix/makeOblique
0 def}if makefont setfont}bd/sf{/_nf xd FontDirectory _nf known{pop}{findfont
dup maxlength dict/_nfd xd{exch dup/FID ne{exch _nfd 3 1 roll put}{pop
pop}ifelse}forall _nfd dup/FontName _nf put/Encoding WE put _nf _nfd
definefont pop}ifelse _nf usf}bd/us{gs [] 0 setdash 0 setlinecap w
0 exch M _w _be add 0 l currentrgbcolor sC S gr}bd/ob{n/_h xd _cx _cy
M _w 0 l 0 _h l _w neg 0 l P tc fill}bd/st{dup/_be xd 0 ne{/_bc xd}if/_vt
xd dup stringwidth pop/_sw xd dup length dup 1 gt{1 sub}if/_sl xd _be
0 eq{_w _sw sub _sl div 0 _vt{exch}if 3 -1 roll ashow}{_be _bc div
0 _vt{exch}if 32 _w _sw sub _sl div 0 _vt{exch}if 6 -1 roll awidthshow}ifelse}bd/sb{gs
M 10 div rotate/_w xd currentpoint/_cy xd/_cx xd 0 ne{ob}if tc _cx
_cy 3 -1 roll add tr 0 0 M st dup 4 and 4 eq{4 sub/_dx _cx _px sub
def _dx neg 0 tr/_w _w _dx add def}if dup 0 eq{pop}{dup 1 eq{pop us}{2
eq{us}{dup 3 1 roll us us}ifelse}ifelse}ifelse _cx/_px xd gr}bd
/cl{gs n M dup neg 3 1 roll 0 l 0 exch l 0 l P clip n}bd/clx{gs n
_cm currentmatrix pop _m setmatrix n M dup neg 3 1 roll 0 l 0 exch
l 0 l P clip _cm setmatrix n}bd
/sbb{/_t save def _neg 0 eq{0}{1}ifelse setgray tr/_dW xd/_dH xd dup/_sW
xd 8 _bs div dup _sW add 1 sub exch div cvi/str exch string def exch
dup/_sH xd 3 -1 roll{/_im{[_sW _dW div 0 0 _sH neg _dH div 0 _sH]}bd}{/_im{[_sW
_dW div 0 0 _sH _dH div 0 0]}bd}ifelse _bs 1 eq{false}{_bs}ifelse _im{currentfile
str readhexstring pop}_bs 1 eq{imagemask}{image}ifelse _t restore}bd
/scbb{/_t save def tr/_dW xd/_dH xd dup/_sW xd 8 _bs div dup _sW add
1 sub exch div cvi 3 mul/str exch string def exch dup/_sH xd _bs 4
-1 roll{/_im{[_sW _dW div 0 0 _sH neg _dH div 0 _sH]}bd}{/_im{[_sW
_dW div 0 0 _sH _dH div 0 0]}bd}ifelse _im{currentfile str readhexstring
pop}false 3 colorimage _t restore}bd
/mr{2 copy -1 eq{_yp _ym 2 mul sub}{0}ifelse exch -1 eq{_xp _xm 2
mul sub}{0}ifelse exch tr scale}bd
/WE 256 array def StandardEncoding WE copy pop WE dup
0[/grave/acute/circumflex/tilde/macron/breve/dotaccent/dieresis/ring/cedilla
/hungarumlaut/ogonek/caron/dotlessi
]putinterval dup
39/quotesingle put dup  96/grave put dup 124/bar put dup
145[/quoteleft/quoteright/quotedblleft/quotedblright/bullet/endash/emdash]putinterval 
160[/space/exclamdown/cent/sterling/currency/yen/brokenbar/section/dieresis
/copyright/ordfeminine/guillemotleft/logicalnot/endash/registered/macron/ring
/plusminus/twosuperior/threesuperior/acute/mu/paragraph/periodcentered
/cedilla/onesuperior/ordmasculine/guillemotright/onequarter/onehalf
/threequarters/questiondown/Agrave/Aacute/Acircumflex/Atilde/Adieresis/Aring
/AE/Ccedilla/Egrave/Eacute/Ecircumflex/Edieresis/Igrave/Iacute/Icircumflex
/Idieresis/Eth/Ntilde/Ograve/Oacute/Ocircumflex/Otilde/Odieresis/circumflex
/Oslash/Ugrave/Uacute/Ucircumflex/Udieresis/Yacute/Thorn/germandbls/agrave
/aacute/acircumflex/atilde/adieresis/aring/ae/ccedilla/egrave/eacute
/ecircumflex/edieresis/igrave/iacute/icircumflex/idieresis/eth/ntilde/ograve
/oacute/ocircumflex/otilde/odieresis/tilde/oslash/ugrave/uacute/ucircumflex
/udieresis/yacute/thorn/ydieresis]putinterval
/TTEV 256 array def TTEV
0 [/G00/G01/G02/G03/G04/G05/G06/G07/G08/G09/G0a/G0b/G0c/G0d/G0e/G0f
   /G10/G11/G12/G13/G14/G15/G16/G17/G18/G19/G1a/G1b/G1c/G1d/G1e/G1f
   /G20/G21/G22/G23/G24/G25/G26/G27/G28/G29/G2a/G2b/G2c/G2d/G2e/G2f
   /G30/G31/G32/G33/G34/G35/G36/G37/G38/G39/G3a/G3b/G3c/G3d/G3e/G3f
   /G40/G41/G42/G43/G44/G45/G46/G47/G48/G49/G4a/G4b/G4c/G4d/G4e/G4f
   /G50/G51/G52/G53/G54/G55/G56/G57/G58/G59/G5a/G5b/G5c/G5d/G5e/G5f
   /G60/G61/G62/G63/G64/G65/G66/G67/G68/G69/G6a/G6b/G6c/G6d/G6e/G6f
   /G70/G71/G72/G73/G74/G75/G76/G77/G78/G79/G7a/G7b/G7c/G7d/G7e/G7f
   /G80/G81/G82/G83/G84/G85/G86/G87/G88/G89/G8a/G8b/G8c/G8d/G8e/G8f
   /G90/G91/G92/G93/G94/G95/G96/G97/G98/G99/G9a/G9b/G9c/G9d/G9e/G9f
   /Ga0/Ga1/Ga2/Ga3/Ga4/Ga5/Ga6/Ga7/Ga8/Ga9/Gaa/Gab/Gac/Gad/Gae/Gaf
   /Gb0/Gb1/Gb2/Gb3/Gb4/Gb5/Gb6/Gb7/Gb8/Gb9/Gba/Gbb/Gbc/Gbd/Gbe/Gbf
   /Gc0/Gc1/Gc2/Gc3/Gc4/Gc5/Gc6/Gc7/Gc8/Gc9/Gca/Gcb/Gcc/Gcd/Gce/Gcf
   /Gd0/Gd1/Gd2/Gd3/Gd4/Gd5/Gd6/Gd7/Gd8/Gd9/Gda/Gdb/Gdc/Gdd/Gde/Gdf
   /Ge0/Ge1/Ge2/Ge3/Ge4/Ge5/Ge6/Ge7/Ge8/Ge9/Gea/Geb/Gec/Ged/Gee/Gef
   /Gf0/Gf1/Gf2/Gf3/Gf4/Gf5/Gf6/Gf7/Gf8/Gf9/Gfa/Gfb/Gfc/Gfd/Gfe/Gff
] putinterval
end
MGXPS_2.8_EPS begin
0 0 1 0 3400 4400 400 -1 0 0 1 Cne
4280 3240 0 0 clpg
0 0 0 sc
/_o 0 def
255 255 255 bc
260 347 78 78 E
0 p
260 343 4 3 E
0 0 0 fc
0 2 p
[5 13 ]0 2 d
149 252 M 300 115 407 57 479 68 bz
547 80 534 283 565 329 bz
597 370 650 370 668 338 bz
685 305 665 182 658 145 bz
655 136 610 60 685 46 bz
830 -8 666 142 738 182 bz
813 220 1008 256 1105 352 bz
1198 447 1327 642 1288 748 bz
1245 852 979 998 855 971 bz
729 941 665 639 545 578 bz
422 517 248 628 136 605 bz
10 575 56 335 149 252 bz
0 p
[]0 0 d
153 175 M 328 37 434 7 513 40 bz
582 71 597 210 604 287 bz
605 364 640 306 625 238 bz
615 188 613 149 607 108 bz
602 64 624 15 709 0 bz
899 9 718 138 785 159 bz
1041 228 1055 237 1162 342 bz
1264 447 1380 617 1344 776 bz
1296 890 1031 1046 871 1034 bz
711 1024 557 632 516 625 bz
476 594 151 705 71 640 bz
-37 550 -30 306 153 175 bz
0 p
175 161 M 0 0 33 38 2 D
S
210 201 M 321 168 421 253 418 348 bz
415 434 350 506 260 506 bz
166 506 88 424 100 318 bz
S
47 291 M 0 0 51 27 2 D
S
2 477 M 0 0 56 4 2 D
S
60 482 M 101 561 193 600 277 593 bz
S
278 593 M 0 0 15 54 2 D
S
367 635 M 0 0 -13 -61 2 D
S
355 573 M 564 493 563 181 322 119 bz
S
321 118 M 0 0 -24 -40 2 D
S
416 31 M 0 0 3 5 2 D
S
418 35 M 0 0 8 36 2 D
S
425 72 M 466 94 500 116 531 159 bz
S
531 160 M 0 0 53 -8 2 D
S
609 79 M 0 0 35 4 2 D
S
644 83 M 687 72 703 67 740 79 bz
S
740 79 M 0 0 45 18 2 D
S
768 148 M 0 0 -47 12 2 D
S
720 159 M 697 149 683 162 663 162 bz
S
661 164 M 0 0 -44 4 2 D
S
610 322 M 0 0 -12 34 2 D
S
597 356 M 591 465 591 488 521 569 bz
S
521 570 M 0 0 -15 52 2 D
S
553 667 M 0 0 53 -36 2 D
S
606 631 M 693 508 690 477 698 414 bz
702 363 708 348 671 337 bz
S
668 339 M 0 0 -45 -23 2 D
S
628 255 M 0 0 44 -5 2 D
S
674 249 M 721 249 750 266 781 296 bz
809 326 808 341 804 412 bz
801 499 738 654 665 720 bz
S
665 722 M 0 0 -58 34 2 D
S
662 848 M 0 0 54 -29 2 D
S
716 819 M 857 699 919 509 923 403 bz
926 313 902 270 863 227 bz
S
863 227 M 0 0 17 -42 2 D
S
1046 247 M 0 0 -24 43 2 D
S
1020 291 M 1048 392 1029 532 1002 612 bz
970 695 889 836 782 919 bz
S
915 1033 M 0 0 -6 -59 2 D
S
909 973 M 1055 863 1155 674 1155 407 bz
S
1105 352 M 0 0 137 -113 2 D
S
1195 379 M 0 0 -40 29 2 D
S
1268 594 M 0 0 49 -29 2 D
S
1268 594 M 1262 673 1224 821 1154 878 bz
S
1154 878 M 0 0 36 49 2 D
S
781 919 M 746 956 L
S
end
EndOfTheIncludedPostscriptMagicCookie

\closepsdump

\psdump{nurhopi.eps}
statusdict begin/waittimeout 0 def end{/MGXPS_2.8_EPS begin}stopped{/MGXPS_2.8_EPS
200 dict def}if MGXPS_2.8_EPS begin/bd{bind def}bind def/xd{exch def}bd/ld{load
def}bd/M/moveto ld/m/rmoveto ld/L/lineto ld/l/rlineto ld/w/setlinewidth
ld/n/newpath ld/P/closepath ld/tr/translate ld/gs/gsave ld/gr/grestore
ld/lc/setlinecap ld/lj/setlinejoin ld/ml/setmiterlimit ld/TM/concat
ld/F/fill ld/bz/curveto ld/_m matrix def/_cm matrix def/_d 0 def/_o
0 def/_pb 0 def/_hb -1 def currentscreen/_P xd dup/_A xd/__A xd dup/_F
xd/__F xd/srgbc{_neg 0 ne{3{neg 1 add 3 1 roll}repeat}if setrgbcolor}bd/ssc{gs
0 0 M 0 _yp neg _xp 0 0 _yp 3 D closepath fill gr 0 0 0 srgbc}bd/Cne{/_bs
xd/_neg where{/_neg get/_neg xd pop}{pop 0/_neg xd}ifelse dup 0 ne{setflat}{pop}ifelse
dup -1 eq{pop}{pop pop}ifelse/_rz xd/_yp xd/_xp xd/_a 0 def 0 ne{/_a
90 def -90 rotate _xp _rz div 72 mul neg 0 tr}if 0 _yp _rz div 72 mul
tr 72 _rz div dup neg scale/#copies xd 2 copy/_ym xd/_xm xd tr n 1
lj 1 lc 1 w 0 8{dup}repeat sc fc bc}bd/clpg{n M dup neg 3 1 roll 0
l 0 exch l 0 l P clip n}bd/SB{_o 0 ne _d 0 ne and{gs []0 setdash _br
_bg _bb srgbc stroke gr}if}bd/S{_r _g _b srgbc SB stroke}bd/SX{_r _g
_b srgbc _cm currentmatrix pop _m setmatrix SB stroke _cm setmatrix}bd/d{/_d
xd setdash}bd/c{3{255 div 3 1 roll}repeat}bd/sc{c/_b xd/_g xd/_r xd}bd/bc{c/_bb
xd/_bg xd/_br xd}bd/fc{c/_B xd/_G xd/_R xd}bd/tc{c srgbc}bd/f{0 eq{eofill}{fill}ifelse}bd/D{{l}repeat}bd/PT{dup
statusdict exch known{cvx exec}{pop}ifelse}bd/TT{dup statusdict exch
known{cvx exec}{pop pop}ifelse}bd/SM{_m currentmatrix pop}bd/RM{_m
setmatrix}bd/sp{dup 0 eq{/_pb xd}{255 div/_gr xd/_pb xd 100 div/_f
xd}ifelse}bd/B{/_t2 xd/_t3 xd _pb _t2 _t3 8 idiv add get 1 7 _t3 8
mod sub bitshift and 0 ne}bd/SF{2{1 add 4 mul cvi exch}repeat B{0}{1}ifelse}bd/pb{_gr
setgray _f _a/SF load setscreen}bd
/s{100 div/_t xd[_t 0 0 _t 1 _t sub _xp 2 div mul 1 _t sub _yp 2 div
mul]concat}bd
/_am matrix def/a{_am currentmatrix pop/_t1 xd/_t2 xd/_t3 exch 10
div def/_t4 exch 10 div def tr _t2 _t1 neg scale/_t xd _t 1 eq{0 0
M}if 0 0 1 _t4 _t3 arc _t 0 ne{P}if _am setmatrix}bd/an{_am currentmatrix
pop/_t1 xd/_t2 xd/_t3 exch 10 div def/_t4 exch 10 div def tr _t2 _t1
neg scale/_t xd _t 1 eq{0 0 M}if 0 0 1 _t4 _t3 arcn _t 0 ne{P}if _am
setmatrix}bd
/_em matrix def/E{_em currentmatrix pop/_t1 xd/_t2 xd/_t3 xd/_t4 xd
_t4 _t3 _t1 sub M _t4 _t3 tr _t2 _t1 neg scale 0 0 1 90 450 arc _em
setmatrix}bd
/_rm matrix def/cn{_rm currentmatrix pop tr x3 y3 scale arc _rm setmatrix}bd/rr{2
div/y3 xd 2 div/x3 xd/y2 xd/x2 xd/y1 xd/x1 xd x1 y1 y3 add M 0 0 1
180 270 x1 x3 add y1 y3 add cn x2 x3 sub y1 L 0 0 1 270 0 x2 x3 sub
y1 y3 add cn x2 y2 y3 sub L 0 0 1 0 90 x2 x3 sub y2 y3 sub cn x1 x3
add y2 L 0 0 1 90 180 x1 x3 add y2 y3 sub cn}bd
/cp{/_t xd{{gs _t 0 eq{eoclip}{clip}ifelse}stopped{gr currentflat
1 add setflat}{exit}ifelse}loop n}bd/p{P dup 0 eq{S pop}{_R _G _B srgbc
_pb 0 eq{2 eq{gs f gr S}{f}ifelse}{_pb 1 eq{_o 1 eq{gs _br _bg _bb
srgbc 1 index f gr}if _R _G _B srgbc 2 eq{gs hf gr S}{hf}ifelse n}{2
eq{true pbc S}{false pbc n}ifelse}ifelse}ifelse}ifelse}bd/px{P dup
0 eq{SX pop}{_R _G _B srgbc _pb 0 eq{2 eq{gs f gr SX}{f}ifelse}{_pb
1 eq{_o 1 eq{gs _br _bg _bb srgbc 1 index f gr}if 2 eq{gs _m setmatrix
hf gr SX}{gs _m setmatrix hf gr}ifelse n}{2 eq{true pbc SX}{false pbc
n}ifelse}ifelse}ifelse}ifelse}bd
/spc{dup 0 eq{/_pb xd}{1 eq{255 div/_gr xd/_pb xd/_pg 0 def/_pr 0
def}{/_pb xd/_pg xd/_pr xd}ifelse 1 add/_pq xd}ifelse}bd/_im{[_t 0
0 _t 0 0]}bd/bbx{_pq 75 mul _rz div/_t xd pathbbox 1 add cvi 31 or/_ury
xd 1 add cvi 31 or/_urx xd 1 sub cvi -32 and/_lly xd 1 sub cvi -32
and/_llx xd/_dH _ury _lly sub def/_dW _urx _llx sub def}bd/xp{/_pbs
xd _row _pbs mul _pbs getinterval/_s xd 0 _pbs 8 _pbs div dup _sW add
1 sub exch div cvi dup/str exch string def _pbs sub{str exch _s putinterval}for
str}bd/pbc{/__t save def{_m setmatrix}if 0 eq{eoclip}{clip}ifelse/_row
0 def bbx _llx _lly tr _t _dW mul ceiling cvi 8 add dup 8 mod sub dup/_sW
xd _t _dH mul ceiling cvi 8 add dup 8 mod sub _pg 0 ne{8 _im{_pr 8
xp}{_pg 8 xp}{_pb 8 xp/_row _row 1 add 8 mod def}true 3 colorimage}{1
_o eq{gs _br _bg _bb srgbc 2 copy true _im{_pb 1 xp/_row _row 1 add
8 mod def}imagemask gr}if _R _G _B srgbc false _im{_pb 1 xp/_row _row
1 add 8 mod def}imagemask}ifelse __t restore}bd
/hp{dup/_hb xd -1 eq{/_pb 0 def}{1 add/_pq xd/_pb 1 def}ifelse}bd/vert{X0
_w X1{dup Y0 M Y1 L stroke}for}bd/horz{Y0 _w Y1{dup X0 exch M X1 exch
L stroke}for}bd/fdiag{X0 _w X1{Y0 M X1 X0 sub dup l stroke}for Y0 _w
Y1{X0 exch M Y1 Y0 sub dup l stroke}for}bd/bdiag{X0 _w X1{Y1 M X1 X0
sub dup neg l stroke}for Y0 _w Y1{X0 exch M Y1 Y0 sub dup neg l stroke}for}bd/AU{1
add cvi 31 or}bd/AD{1 sub cvi -32 and}bd/hf{pathbbox 1 add cvi 31 or/Y1
xd 1 add cvi 31 or/X1 xd 1 sub cvi -32 and/Y0 xd 1 sub cvi -32 and/X0
xd 2 w [] 0 setdash/_w _rz 20 div 8 div 2 mul _pq div round 8 mul def
cp _hb 0 eq{horz}if _hb 1 eq{vert}if _hb 2 eq{fdiag}if _hb 3 eq{bdiag}if
_hb 4 eq{horz vert}if _hb 5 eq{fdiag bdiag}if gr}bd
/sp{dup 0 eq{/_pb xd}{255 div/_gr xd/_pb xd 100 div/_f xd}ifelse}bd/B{/_t2
xd/_t3 xd _pb _t2 _t3 8 idiv add get 1 7 _t3 8 mod sub bitshift and
0 ne}bd/SF{2{1 add 4 mul cvi exch}repeat B{0}{1}ifelse}bd/pb{_gr setgray
_f _a/SF load setscreen}bd
/SCA{10 div dup dup/_a xd/_A xd/__na xd currentscreen/__p xd/__a xd/__f
xd __f __na/__p load setscreen}bd
/SCF{100 div dup dup/_F xd/__nf xd currentscreen/__p xd/__a xd/__f
xd __nf __a/__p load setscreen}bd
/makeOblique 0 def/bC 1 def/makeoutlinedict 7 dict def/mbf{makeoutlinedict
begin/uniqueid exch def/strokewidth exch def/newfontname exch def/basefontname
exch def/basefontdict basefontname findfont def/numentries basefontdict
maxlength 1 add def basefontdict/UniqueID known not{/numentries numentries
1 add def}if/outfontdict numentries dict def basefontdict{exch dup/FID
ne{exch outfontdict 3 1 roll put}{pop pop}ifelse}forall outfontdict/FontName
newfontname put outfontdict/PaintType 2 put outfontdict/StrokeWidth
strokewidth put outfontdict/UniqueID uniqueid put newfontname outfontdict
definefont pop end}bd/usf{findfont exch makeOblique 0 ne{/makeOblique
makeOblique sin neg def [1 0 makeOblique 1 0 0] matrix concatmatrix/makeOblique
0 def}if makefont setfont}bd/sf{/_nf xd FontDirectory _nf known{pop}{findfont
dup maxlength dict/_nfd xd{exch dup/FID ne{exch _nfd 3 1 roll put}{pop
pop}ifelse}forall _nfd dup/FontName _nf put/Encoding WE put _nf _nfd
definefont pop}ifelse _nf usf}bd/us{gs [] 0 setdash 0 setlinecap w
0 exch M _w _be add 0 l currentrgbcolor sC S gr}bd/ob{n/_h xd _cx _cy
M _w 0 l 0 _h l _w neg 0 l P tc fill}bd/st{dup/_be xd 0 ne{/_bc xd}if/_vt
xd dup stringwidth pop/_sw xd dup length dup 1 gt{1 sub}if/_sl xd _be
0 eq{_w _sw sub _sl div 0 _vt{exch}if 3 -1 roll ashow}{_be _bc div
0 _vt{exch}if 32 _w _sw sub _sl div 0 _vt{exch}if 6 -1 roll awidthshow}ifelse}bd/sb{gs
M 10 div rotate/_w xd currentpoint/_cy xd/_cx xd 0 ne{ob}if tc _cx
_cy 3 -1 roll add tr 0 0 M st dup 4 and 4 eq{4 sub/_dx _cx _px sub
def _dx neg 0 tr/_w _w _dx add def}if dup 0 eq{pop}{dup 1 eq{pop us}{2
eq{us}{dup 3 1 roll us us}ifelse}ifelse}ifelse _cx/_px xd gr}bd
/cl{gs n M dup neg 3 1 roll 0 l 0 exch l 0 l P clip n}bd/clx{gs n
_cm currentmatrix pop _m setmatrix n M dup neg 3 1 roll 0 l 0 exch
l 0 l P clip _cm setmatrix n}bd
/sbb{/_t save def _neg 0 eq{0}{1}ifelse setgray tr/_dW xd/_dH xd dup/_sW
xd 8 _bs div dup _sW add 1 sub exch div cvi/str exch string def exch
dup/_sH xd 3 -1 roll{/_im{[_sW _dW div 0 0 _sH neg _dH div 0 _sH]}bd}{/_im{[_sW
_dW div 0 0 _sH _dH div 0 0]}bd}ifelse _bs 1 eq{false}{_bs}ifelse _im{currentfile
str readhexstring pop}_bs 1 eq{imagemask}{image}ifelse _t restore}bd
/scbb{/_t save def tr/_dW xd/_dH xd dup/_sW xd 8 _bs div dup _sW add
1 sub exch div cvi 3 mul/str exch string def exch dup/_sH xd _bs 4
-1 roll{/_im{[_sW _dW div 0 0 _sH neg _dH div 0 _sH]}bd}{/_im{[_sW
_dW div 0 0 _sH _dH div 0 0]}bd}ifelse _im{currentfile str readhexstring
pop}false 3 colorimage _t restore}bd
/mr{2 copy -1 eq{_yp _ym 2 mul sub}{0}ifelse exch -1 eq{_xp _xm 2
mul sub}{0}ifelse exch tr scale}bd
/WE 256 array def StandardEncoding WE copy pop WE dup
0[/grave/acute/circumflex/tilde/macron/breve/dotaccent/dieresis/ring/cedilla
/hungarumlaut/ogonek/caron/dotlessi
]putinterval dup
39/quotesingle put dup  96/grave put dup 124/bar put dup
145[/quoteleft/quoteright/quotedblleft/quotedblright/bullet/endash/emdash]putinterval 
160[/space/exclamdown/cent/sterling/currency/yen/brokenbar/section/dieresis
/copyright/ordfeminine/guillemotleft/logicalnot/endash/registered/macron/ring
/plusminus/twosuperior/threesuperior/acute/mu/paragraph/periodcentered
/cedilla/onesuperior/ordmasculine/guillemotright/onequarter/onehalf
/threequarters/questiondown/Agrave/Aacute/Acircumflex/Atilde/Adieresis/Aring
/AE/Ccedilla/Egrave/Eacute/Ecircumflex/Edieresis/Igrave/Iacute/Icircumflex
/Idieresis/Eth/Ntilde/Ograve/Oacute/Ocircumflex/Otilde/Odieresis/circumflex
/Oslash/Ugrave/Uacute/Ucircumflex/Udieresis/Yacute/Thorn/germandbls/agrave
/aacute/acircumflex/atilde/adieresis/aring/ae/ccedilla/egrave/eacute
/ecircumflex/edieresis/igrave/iacute/icircumflex/idieresis/eth/ntilde/ograve
/oacute/ocircumflex/otilde/odieresis/tilde/oslash/ugrave/uacute/ucircumflex
/udieresis/yacute/thorn/ydieresis]putinterval
/TTEV 256 array def TTEV
0 [/G00/G01/G02/G03/G04/G05/G06/G07/G08/G09/G0a/G0b/G0c/G0d/G0e/G0f
   /G10/G11/G12/G13/G14/G15/G16/G17/G18/G19/G1a/G1b/G1c/G1d/G1e/G1f
   /G20/G21/G22/G23/G24/G25/G26/G27/G28/G29/G2a/G2b/G2c/G2d/G2e/G2f
   /G30/G31/G32/G33/G34/G35/G36/G37/G38/G39/G3a/G3b/G3c/G3d/G3e/G3f
   /G40/G41/G42/G43/G44/G45/G46/G47/G48/G49/G4a/G4b/G4c/G4d/G4e/G4f
   /G50/G51/G52/G53/G54/G55/G56/G57/G58/G59/G5a/G5b/G5c/G5d/G5e/G5f
   /G60/G61/G62/G63/G64/G65/G66/G67/G68/G69/G6a/G6b/G6c/G6d/G6e/G6f
   /G70/G71/G72/G73/G74/G75/G76/G77/G78/G79/G7a/G7b/G7c/G7d/G7e/G7f
   /G80/G81/G82/G83/G84/G85/G86/G87/G88/G89/G8a/G8b/G8c/G8d/G8e/G8f
   /G90/G91/G92/G93/G94/G95/G96/G97/G98/G99/G9a/G9b/G9c/G9d/G9e/G9f
   /Ga0/Ga1/Ga2/Ga3/Ga4/Ga5/Ga6/Ga7/Ga8/Ga9/Gaa/Gab/Gac/Gad/Gae/Gaf
   /Gb0/Gb1/Gb2/Gb3/Gb4/Gb5/Gb6/Gb7/Gb8/Gb9/Gba/Gbb/Gbc/Gbd/Gbe/Gbf
   /Gc0/Gc1/Gc2/Gc3/Gc4/Gc5/Gc6/Gc7/Gc8/Gc9/Gca/Gcb/Gcc/Gcd/Gce/Gcf
   /Gd0/Gd1/Gd2/Gd3/Gd4/Gd5/Gd6/Gd7/Gd8/Gd9/Gda/Gdb/Gdc/Gdd/Gde/Gdf
   /Ge0/Ge1/Ge2/Ge3/Ge4/Ge5/Ge6/Ge7/Ge8/Ge9/Gea/Geb/Gec/Ged/Gee/Gef
   /Gf0/Gf1/Gf2/Gf3/Gf4/Gf5/Gf6/Gf7/Gf8/Gf9/Gfa/Gfb/Gfc/Gfd/Gfe/Gff
] putinterval
end
MGXPS_2.8_EPS begin
0 0 1 0 3400 4400 400 -1 0 0 1 Cne
4280 3240 0 0 clpg
0 0 0 sc
/_o 0 def
255 255 255 bc
0 466 M 232 484 434 456 609 377 bz
786 300 927 174 1046 0 bz
S
0 655 M 232 672 434 640 609 561 bz
786 484 927 363 1046 189 bz
S
373 1270 M 605 1287 807 1255 982 1176 bz
1158 1099 1300 978 1419 804 bz
S
594 425 M 0 0 0 141 2 D
S
599 430 M 899 613 1176 751 1198 1180 bz
S
594 566 M 799 694 1070 845 1103 1112 bz
S
1103 1112 M 0 0 95 68 2 D
S
900 733 M -45 2 3 45 2 D
S
627 485 M -31 -32 -32 32 2 D
S
1143 1184 M 45 -5 -4 -44 2 D
S
942 694 M 44 -6 -6 -44 2 D
S
1012 333 M 0 0 286 -160 2 D
S
1485 899 M 0 0 147 -243 2 D
S
1169 1165 M 0 0 110 199 2 D
S
0 828 M 232 845 434 818 609 739 bz
786 661 927 536 1046 362 bz
S
[5 13 ]0 2 d
0 466 M 0 0 0 361 2 D
S
1046 1 M 0 0 0 362 2 D
S
1245 771 M 0 0 372 63 2 D
S
399 1137 M 0 0 -54 279 2 D
S
[]0 0 d
346 1416 M 578 1438 940 1408 1108 1317 bz
1284 1224 1520 1039 1618 834 bz
S
399 1137 M 626 1150 738 1112 898 1045 bz
1065 976 1121 940 1246 771 bz
S
end
EndOfTheIncludedPostscriptMagicCookie

\closepsdump

\psdump{repbdy.eps}
statusdict begin/waittimeout 0 def end{/MGXPS_2.8_EPS begin}stopped{/MGXPS_2.8_EPS
200 dict def}if MGXPS_2.8_EPS begin/bd{bind def}bind def/xd{exch def}bd/ld{load
def}bd/M/moveto ld/m/rmoveto ld/L/lineto ld/l/rlineto ld/w/setlinewidth
ld/n/newpath ld/P/closepath ld/tr/translate ld/gs/gsave ld/gr/grestore
ld/lc/setlinecap ld/lj/setlinejoin ld/ml/setmiterlimit ld/TM/concat
ld/F/fill ld/bz/curveto ld/_m matrix def/_cm matrix def/_d 0 def/_o
0 def/_pb 0 def/_hb -1 def currentscreen/_P xd dup/_A xd/__A xd dup/_F
xd/__F xd/srgbc{_neg 0 ne{3{neg 1 add 3 1 roll}repeat}if setrgbcolor}bd/ssc{gs
0 0 M 0 _yp neg _xp 0 0 _yp 3 D closepath fill gr 0 0 0 srgbc}bd/Cne{/_bs
xd/_neg where{/_neg get/_neg xd pop}{pop 0/_neg xd}ifelse dup 0 ne{setflat}{pop}ifelse
dup -1 eq{pop}{pop pop}ifelse/_rz xd/_yp xd/_xp xd/_a 0 def 0 ne{/_a
90 def -90 rotate _xp _rz div 72 mul neg 0 tr}if 0 _yp _rz div 72 mul
tr 72 _rz div dup neg scale/#copies xd 2 copy/_ym xd/_xm xd tr n 1
lj 1 lc 1 w 0 8{dup}repeat sc fc bc}bd/clpg{n M dup neg 3 1 roll 0
l 0 exch l 0 l P clip n}bd/SB{_o 0 ne _d 0 ne and{gs []0 setdash _br
_bg _bb srgbc stroke gr}if}bd/S{_r _g _b srgbc SB stroke}bd/SX{_r _g
_b srgbc _cm currentmatrix pop _m setmatrix SB stroke _cm setmatrix}bd/d{/_d
xd setdash}bd/c{3{255 div 3 1 roll}repeat}bd/sc{c/_b xd/_g xd/_r xd}bd/bc{c/_bb
xd/_bg xd/_br xd}bd/fc{c/_B xd/_G xd/_R xd}bd/tc{c srgbc}bd/f{0 eq{eofill}{fill}ifelse}bd/D{{l}repeat}bd/PT{dup
statusdict exch known{cvx exec}{pop}ifelse}bd/TT{dup statusdict exch
known{cvx exec}{pop pop}ifelse}bd/SM{_m currentmatrix pop}bd/RM{_m
setmatrix}bd/sp{dup 0 eq{/_pb xd}{255 div/_gr xd/_pb xd 100 div/_f
xd}ifelse}bd/B{/_t2 xd/_t3 xd _pb _t2 _t3 8 idiv add get 1 7 _t3 8
mod sub bitshift and 0 ne}bd/SF{2{1 add 4 mul cvi exch}repeat B{0}{1}ifelse}bd/pb{_gr
setgray _f _a/SF load setscreen}bd
/s{100 div/_t xd[_t 0 0 _t 1 _t sub _xp 2 div mul 1 _t sub _yp 2 div
mul]concat}bd
/_am matrix def/a{_am currentmatrix pop/_t1 xd/_t2 xd/_t3 exch 10
div def/_t4 exch 10 div def tr _t2 _t1 neg scale/_t xd _t 1 eq{0 0
M}if 0 0 1 _t4 _t3 arc _t 0 ne{P}if _am setmatrix}bd/an{_am currentmatrix
pop/_t1 xd/_t2 xd/_t3 exch 10 div def/_t4 exch 10 div def tr _t2 _t1
neg scale/_t xd _t 1 eq{0 0 M}if 0 0 1 _t4 _t3 arcn _t 0 ne{P}if _am
setmatrix}bd
/_em matrix def/E{_em currentmatrix pop/_t1 xd/_t2 xd/_t3 xd/_t4 xd
_t4 _t3 _t1 sub M _t4 _t3 tr _t2 _t1 neg scale 0 0 1 90 450 arc _em
setmatrix}bd
/_rm matrix def/cn{_rm currentmatrix pop tr x3 y3 scale arc _rm setmatrix}bd/rr{2
div/y3 xd 2 div/x3 xd/y2 xd/x2 xd/y1 xd/x1 xd x1 y1 y3 add M 0 0 1
180 270 x1 x3 add y1 y3 add cn x2 x3 sub y1 L 0 0 1 270 0 x2 x3 sub
y1 y3 add cn x2 y2 y3 sub L 0 0 1 0 90 x2 x3 sub y2 y3 sub cn x1 x3
add y2 L 0 0 1 90 180 x1 x3 add y2 y3 sub cn}bd
/cp{/_t xd{{gs _t 0 eq{eoclip}{clip}ifelse}stopped{gr currentflat
1 add setflat}{exit}ifelse}loop n}bd/p{P dup 0 eq{S pop}{_R _G _B srgbc
_pb 0 eq{2 eq{gs f gr S}{f}ifelse}{_pb 1 eq{_o 1 eq{gs _br _bg _bb
srgbc 1 index f gr}if _R _G _B srgbc 2 eq{gs hf gr S}{hf}ifelse n}{2
eq{true pbc S}{false pbc n}ifelse}ifelse}ifelse}ifelse}bd/px{P dup
0 eq{SX pop}{_R _G _B srgbc _pb 0 eq{2 eq{gs f gr SX}{f}ifelse}{_pb
1 eq{_o 1 eq{gs _br _bg _bb srgbc 1 index f gr}if 2 eq{gs _m setmatrix
hf gr SX}{gs _m setmatrix hf gr}ifelse n}{2 eq{true pbc SX}{false pbc
n}ifelse}ifelse}ifelse}ifelse}bd
/spc{dup 0 eq{/_pb xd}{1 eq{255 div/_gr xd/_pb xd/_pg 0 def/_pr 0
def}{/_pb xd/_pg xd/_pr xd}ifelse 1 add/_pq xd}ifelse}bd/_im{[_t 0
0 _t 0 0]}bd/bbx{_pq 75 mul _rz div/_t xd pathbbox 1 add cvi 31 or/_ury
xd 1 add cvi 31 or/_urx xd 1 sub cvi -32 and/_lly xd 1 sub cvi -32
and/_llx xd/_dH _ury _lly sub def/_dW _urx _llx sub def}bd/xp{/_pbs
xd _row _pbs mul _pbs getinterval/_s xd 0 _pbs 8 _pbs div dup _sW add
1 sub exch div cvi dup/str exch string def _pbs sub{str exch _s putinterval}for
str}bd/pbc{/__t save def{_m setmatrix}if 0 eq{eoclip}{clip}ifelse/_row
0 def bbx _llx _lly tr _t _dW mul ceiling cvi 8 add dup 8 mod sub dup/_sW
xd _t _dH mul ceiling cvi 8 add dup 8 mod sub _pg 0 ne{8 _im{_pr 8
xp}{_pg 8 xp}{_pb 8 xp/_row _row 1 add 8 mod def}true 3 colorimage}{1
_o eq{gs _br _bg _bb srgbc 2 copy true _im{_pb 1 xp/_row _row 1 add
8 mod def}imagemask gr}if _R _G _B srgbc false _im{_pb 1 xp/_row _row
1 add 8 mod def}imagemask}ifelse __t restore}bd
/hp{dup/_hb xd -1 eq{/_pb 0 def}{1 add/_pq xd/_pb 1 def}ifelse}bd/vert{X0
_w X1{dup Y0 M Y1 L stroke}for}bd/horz{Y0 _w Y1{dup X0 exch M X1 exch
L stroke}for}bd/fdiag{X0 _w X1{Y0 M X1 X0 sub dup l stroke}for Y0 _w
Y1{X0 exch M Y1 Y0 sub dup l stroke}for}bd/bdiag{X0 _w X1{Y1 M X1 X0
sub dup neg l stroke}for Y0 _w Y1{X0 exch M Y1 Y0 sub dup neg l stroke}for}bd/AU{1
add cvi 31 or}bd/AD{1 sub cvi -32 and}bd/hf{pathbbox 1 add cvi 31 or/Y1
xd 1 add cvi 31 or/X1 xd 1 sub cvi -32 and/Y0 xd 1 sub cvi -32 and/X0
xd 2 w [] 0 setdash/_w _rz 20 div 8 div 2 mul _pq div round 8 mul def
cp _hb 0 eq{horz}if _hb 1 eq{vert}if _hb 2 eq{fdiag}if _hb 3 eq{bdiag}if
_hb 4 eq{horz vert}if _hb 5 eq{fdiag bdiag}if gr}bd
/sp{dup 0 eq{/_pb xd}{255 div/_gr xd/_pb xd 100 div/_f xd}ifelse}bd/B{/_t2
xd/_t3 xd _pb _t2 _t3 8 idiv add get 1 7 _t3 8 mod sub bitshift and
0 ne}bd/SF{2{1 add 4 mul cvi exch}repeat B{0}{1}ifelse}bd/pb{_gr setgray
_f _a/SF load setscreen}bd
/SCA{10 div dup dup/_a xd/_A xd/__na xd currentscreen/__p xd/__a xd/__f
xd __f __na/__p load setscreen}bd
/SCF{100 div dup dup/_F xd/__nf xd currentscreen/__p xd/__a xd/__f
xd __nf __a/__p load setscreen}bd
/makeOblique 0 def/bC 1 def/makeoutlinedict 7 dict def/mbf{makeoutlinedict
begin/uniqueid exch def/strokewidth exch def/newfontname exch def/basefontname
exch def/basefontdict basefontname findfont def/numentries basefontdict
maxlength 1 add def basefontdict/UniqueID known not{/numentries numentries
1 add def}if/outfontdict numentries dict def basefontdict{exch dup/FID
ne{exch outfontdict 3 1 roll put}{pop pop}ifelse}forall outfontdict/FontName
newfontname put outfontdict/PaintType 2 put outfontdict/StrokeWidth
strokewidth put outfontdict/UniqueID uniqueid put newfontname outfontdict
definefont pop end}bd/usf{findfont exch makeOblique 0 ne{/makeOblique
makeOblique sin neg def [1 0 makeOblique 1 0 0] matrix concatmatrix/makeOblique
0 def}if makefont setfont}bd/sf{/_nf xd FontDirectory _nf known{pop}{findfont
dup maxlength dict/_nfd xd{exch dup/FID ne{exch _nfd 3 1 roll put}{pop
pop}ifelse}forall _nfd dup/FontName _nf put/Encoding WE put _nf _nfd
definefont pop}ifelse _nf usf}bd/us{gs [] 0 setdash 0 setlinecap w
0 exch M _w _be add 0 l currentrgbcolor sC S gr}bd/ob{n/_h xd _cx _cy
M _w 0 l 0 _h l _w neg 0 l P tc fill}bd/st{dup/_be xd 0 ne{/_bc xd}if/_vt
xd dup stringwidth pop/_sw xd dup length dup 1 gt{1 sub}if/_sl xd _be
0 eq{_w _sw sub _sl div 0 _vt{exch}if 3 -1 roll ashow}{_be _bc div
0 _vt{exch}if 32 _w _sw sub _sl div 0 _vt{exch}if 6 -1 roll awidthshow}ifelse}bd/sb{gs
M 10 div rotate/_w xd currentpoint/_cy xd/_cx xd 0 ne{ob}if tc _cx
_cy 3 -1 roll add tr 0 0 M st dup 4 and 4 eq{4 sub/_dx _cx _px sub
def _dx neg 0 tr/_w _w _dx add def}if dup 0 eq{pop}{dup 1 eq{pop us}{2
eq{us}{dup 3 1 roll us us}ifelse}ifelse}ifelse _cx/_px xd gr}bd
/cl{gs n M dup neg 3 1 roll 0 l 0 exch l 0 l P clip n}bd/clx{gs n
_cm currentmatrix pop _m setmatrix n M dup neg 3 1 roll 0 l 0 exch
l 0 l P clip _cm setmatrix n}bd
/sbb{/_t save def _neg 0 eq{0}{1}ifelse setgray tr/_dW xd/_dH xd dup/_sW
xd 8 _bs div dup _sW add 1 sub exch div cvi/str exch string def exch
dup/_sH xd 3 -1 roll{/_im{[_sW _dW div 0 0 _sH neg _dH div 0 _sH]}bd}{/_im{[_sW
_dW div 0 0 _sH _dH div 0 0]}bd}ifelse _bs 1 eq{false}{_bs}ifelse _im{currentfile
str readhexstring pop}_bs 1 eq{imagemask}{image}ifelse _t restore}bd
/scbb{/_t save def tr/_dW xd/_dH xd dup/_sW xd 8 _bs div dup _sW add
1 sub exch div cvi 3 mul/str exch string def exch dup/_sH xd _bs 4
-1 roll{/_im{[_sW _dW div 0 0 _sH neg _dH div 0 _sH]}bd}{/_im{[_sW
_dW div 0 0 _sH _dH div 0 0]}bd}ifelse _im{currentfile str readhexstring
pop}false 3 colorimage _t restore}bd
/mr{2 copy -1 eq{_yp _ym 2 mul sub}{0}ifelse exch -1 eq{_xp _xm 2
mul sub}{0}ifelse exch tr scale}bd
/WE 256 array def StandardEncoding WE copy pop WE dup
0[/grave/acute/circumflex/tilde/macron/breve/dotaccent/dieresis/ring/cedilla
/hungarumlaut/ogonek/caron/dotlessi
]putinterval dup
39/quotesingle put dup  96/grave put dup 124/bar put dup
145[/quoteleft/quoteright/quotedblleft/quotedblright/bullet/endash/emdash]putinterval 
160[/space/exclamdown/cent/sterling/currency/yen/brokenbar/section/dieresis
/copyright/ordfeminine/guillemotleft/logicalnot/endash/registered/macron/ring
/plusminus/twosuperior/threesuperior/acute/mu/paragraph/periodcentered
/cedilla/onesuperior/ordmasculine/guillemotright/onequarter/onehalf
/threequarters/questiondown/Agrave/Aacute/Acircumflex/Atilde/Adieresis/Aring
/AE/Ccedilla/Egrave/Eacute/Ecircumflex/Edieresis/Igrave/Iacute/Icircumflex
/Idieresis/Eth/Ntilde/Ograve/Oacute/Ocircumflex/Otilde/Odieresis/circumflex
/Oslash/Ugrave/Uacute/Ucircumflex/Udieresis/Yacute/Thorn/germandbls/agrave
/aacute/acircumflex/atilde/adieresis/aring/ae/ccedilla/egrave/eacute
/ecircumflex/edieresis/igrave/iacute/icircumflex/idieresis/eth/ntilde/ograve
/oacute/ocircumflex/otilde/odieresis/tilde/oslash/ugrave/uacute/ucircumflex
/udieresis/yacute/thorn/ydieresis]putinterval
/TTEV 256 array def TTEV
0 [/G00/G01/G02/G03/G04/G05/G06/G07/G08/G09/G0a/G0b/G0c/G0d/G0e/G0f
   /G10/G11/G12/G13/G14/G15/G16/G17/G18/G19/G1a/G1b/G1c/G1d/G1e/G1f
   /G20/G21/G22/G23/G24/G25/G26/G27/G28/G29/G2a/G2b/G2c/G2d/G2e/G2f
   /G30/G31/G32/G33/G34/G35/G36/G37/G38/G39/G3a/G3b/G3c/G3d/G3e/G3f
   /G40/G41/G42/G43/G44/G45/G46/G47/G48/G49/G4a/G4b/G4c/G4d/G4e/G4f
   /G50/G51/G52/G53/G54/G55/G56/G57/G58/G59/G5a/G5b/G5c/G5d/G5e/G5f
   /G60/G61/G62/G63/G64/G65/G66/G67/G68/G69/G6a/G6b/G6c/G6d/G6e/G6f
   /G70/G71/G72/G73/G74/G75/G76/G77/G78/G79/G7a/G7b/G7c/G7d/G7e/G7f
   /G80/G81/G82/G83/G84/G85/G86/G87/G88/G89/G8a/G8b/G8c/G8d/G8e/G8f
   /G90/G91/G92/G93/G94/G95/G96/G97/G98/G99/G9a/G9b/G9c/G9d/G9e/G9f
   /Ga0/Ga1/Ga2/Ga3/Ga4/Ga5/Ga6/Ga7/Ga8/Ga9/Gaa/Gab/Gac/Gad/Gae/Gaf
   /Gb0/Gb1/Gb2/Gb3/Gb4/Gb5/Gb6/Gb7/Gb8/Gb9/Gba/Gbb/Gbc/Gbd/Gbe/Gbf
   /Gc0/Gc1/Gc2/Gc3/Gc4/Gc5/Gc6/Gc7/Gc8/Gc9/Gca/Gcb/Gcc/Gcd/Gce/Gcf
   /Gd0/Gd1/Gd2/Gd3/Gd4/Gd5/Gd6/Gd7/Gd8/Gd9/Gda/Gdb/Gdc/Gdd/Gde/Gdf
   /Ge0/Ge1/Ge2/Ge3/Ge4/Ge5/Ge6/Ge7/Ge8/Ge9/Gea/Geb/Gec/Ged/Gee/Gef
   /Gf0/Gf1/Gf2/Gf3/Gf4/Gf5/Gf6/Gf7/Gf8/Gf9/Gfa/Gfb/Gfc/Gfd/Gfe/Gff
] putinterval
end
MGXPS_2.8_EPS begin
0 0 1 0 3400 4400 400 -1 0 0 1 Cne
4280 3240 0 0 clpg
0 0 0 sc
/_o 0 def
255 255 255 bc
277 776 M 208 874 98 922 53 930 bz
7 935 -9 894 4 810 bz
S
1216 557 M 1268 473 1309 407 1335 260 bz
S
82 117 M 200 0 0 176 2 D
S
66 133 M 16 16 16 -16 2 D
S
617 490 7 7 E
0 0 0 fc
0 2 p
267 278 M -15 16 15 15 2 D
S
11 w
10433 1000 div ml
278 776 M 400 589 614 463 737 464 bz
849 466 919 627 1011 656 bz
1101 682 1166 646 1215 558 bz
S
1 w
726 270 M 8 32 33 -8 2 D
S
823 584 M -32 7 7 33 2 D
S
SM
[0.7863 0.0137 -0.0137 0.7863 611 487]TM
9 9 9 9 E
0 2 px
RM
[5 13 ]0 2 d
1092 775 M 0 0 122 -202 2 D
S
[]0 0 d
620 492 M 0 0 138 -230 2 D
S
[5 13 ]0 2 d
492 418 M 0 0 -216 357 2 D
S
6 w
[]0 0 d
492 419 M 0 0 600 355 2 D
S
1 w
278 777 M 149 687 424 543 620 352 bz
816 161 898 -26 1068 3 bz
1189 24 1217 46 1245 125 bz
1287 245 1216 429 1217 556 bz
1217 696 1261 674 1234 779 bz
1202 904 1179 870 946 896 bz
709 921 423 879 278 777 bz
0 p
end
EndOfTheIncludedPostscriptMagicCookie

\closepsdump

\psdump{tau_of_n.eps}
statusdict begin/waittimeout 0 def end{/MGXPS_2.8_EPS begin}stopped{/MGXPS_2.8_EPS
200 dict def}if MGXPS_2.8_EPS begin/bd{bind def}bind def/xd{exch def}bd/ld{load
def}bd/M/moveto ld/m/rmoveto ld/L/lineto ld/l/rlineto ld/w/setlinewidth
ld/n/newpath ld/P/closepath ld/tr/translate ld/gs/gsave ld/gr/grestore
ld/lc/setlinecap ld/lj/setlinejoin ld/ml/setmiterlimit ld/TM/concat
ld/F/fill ld/bz/curveto ld/_m matrix def/_cm matrix def/_d 0 def/_o
0 def/_pb 0 def/_hb -1 def currentscreen/_P xd dup/_A xd/__A xd dup/_F
xd/__F xd/srgbc{_neg 0 ne{3{neg 1 add 3 1 roll}repeat}if setrgbcolor}bd/ssc{gs
0 0 M 0 _yp neg _xp 0 0 _yp 3 D closepath fill gr 0 0 0 srgbc}bd/Cne{/_bs
xd/_neg where{/_neg get/_neg xd pop}{pop 0/_neg xd}ifelse dup 0 ne{setflat}{pop}ifelse
dup -1 eq{pop}{pop pop}ifelse/_rz xd/_yp xd/_xp xd/_a 0 def 0 ne{/_a
90 def -90 rotate _xp _rz div 72 mul neg 0 tr}if 0 _yp _rz div 72 mul
tr 72 _rz div dup neg scale/#copies xd 2 copy/_ym xd/_xm xd tr n 1
lj 1 lc 1 w 0 8{dup}repeat sc fc bc}bd/clpg{n M dup neg 3 1 roll 0
l 0 exch l 0 l P clip n}bd/SB{_o 0 ne _d 0 ne and{gs []0 setdash _br
_bg _bb srgbc stroke gr}if}bd/S{_r _g _b srgbc SB stroke}bd/SX{_r _g
_b srgbc _cm currentmatrix pop _m setmatrix SB stroke _cm setmatrix}bd/d{/_d
xd setdash}bd/c{3{255 div 3 1 roll}repeat}bd/sc{c/_b xd/_g xd/_r xd}bd/bc{c/_bb
xd/_bg xd/_br xd}bd/fc{c/_B xd/_G xd/_R xd}bd/tc{c srgbc}bd/f{0 eq{eofill}{fill}ifelse}bd/D{{l}repeat}bd/PT{dup
statusdict exch known{cvx exec}{pop}ifelse}bd/TT{dup statusdict exch
known{cvx exec}{pop pop}ifelse}bd/SM{_m currentmatrix pop}bd/RM{_m
setmatrix}bd/sp{dup 0 eq{/_pb xd}{255 div/_gr xd/_pb xd 100 div/_f
xd}ifelse}bd/B{/_t2 xd/_t3 xd _pb _t2 _t3 8 idiv add get 1 7 _t3 8
mod sub bitshift and 0 ne}bd/SF{2{1 add 4 mul cvi exch}repeat B{0}{1}ifelse}bd/pb{_gr
setgray _f _a/SF load setscreen}bd
/s{100 div/_t xd[_t 0 0 _t 1 _t sub _xp 2 div mul 1 _t sub _yp 2 div
mul]concat}bd
/_am matrix def/a{_am currentmatrix pop/_t1 xd/_t2 xd/_t3 exch 10
div def/_t4 exch 10 div def tr _t2 _t1 neg scale/_t xd _t 1 eq{0 0
M}if 0 0 1 _t4 _t3 arc _t 0 ne{P}if _am setmatrix}bd/an{_am currentmatrix
pop/_t1 xd/_t2 xd/_t3 exch 10 div def/_t4 exch 10 div def tr _t2 _t1
neg scale/_t xd _t 1 eq{0 0 M}if 0 0 1 _t4 _t3 arcn _t 0 ne{P}if _am
setmatrix}bd
/_em matrix def/E{_em currentmatrix pop/_t1 xd/_t2 xd/_t3 xd/_t4 xd
_t4 _t3 _t1 sub M _t4 _t3 tr _t2 _t1 neg scale 0 0 1 90 450 arc _em
setmatrix}bd
/_rm matrix def/cn{_rm currentmatrix pop tr x3 y3 scale arc _rm setmatrix}bd/rr{2
div/y3 xd 2 div/x3 xd/y2 xd/x2 xd/y1 xd/x1 xd x1 y1 y3 add M 0 0 1
180 270 x1 x3 add y1 y3 add cn x2 x3 sub y1 L 0 0 1 270 0 x2 x3 sub
y1 y3 add cn x2 y2 y3 sub L 0 0 1 0 90 x2 x3 sub y2 y3 sub cn x1 x3
add y2 L 0 0 1 90 180 x1 x3 add y2 y3 sub cn}bd
/cp{/_t xd{{gs _t 0 eq{eoclip}{clip}ifelse}stopped{gr currentflat
1 add setflat}{exit}ifelse}loop n}bd/p{P dup 0 eq{S pop}{_R _G _B srgbc
_pb 0 eq{2 eq{gs f gr S}{f}ifelse}{_pb 1 eq{_o 1 eq{gs _br _bg _bb
srgbc 1 index f gr}if _R _G _B srgbc 2 eq{gs hf gr S}{hf}ifelse n}{2
eq{true pbc S}{false pbc n}ifelse}ifelse}ifelse}ifelse}bd/px{P dup
0 eq{SX pop}{_R _G _B srgbc _pb 0 eq{2 eq{gs f gr SX}{f}ifelse}{_pb
1 eq{_o 1 eq{gs _br _bg _bb srgbc 1 index f gr}if 2 eq{gs _m setmatrix
hf gr SX}{gs _m setmatrix hf gr}ifelse n}{2 eq{true pbc SX}{false pbc
n}ifelse}ifelse}ifelse}ifelse}bd
/spc{dup 0 eq{/_pb xd}{1 eq{255 div/_gr xd/_pb xd/_pg 0 def/_pr 0
def}{/_pb xd/_pg xd/_pr xd}ifelse 1 add/_pq xd}ifelse}bd/_im{[_t 0
0 _t 0 0]}bd/bbx{_pq 75 mul _rz div/_t xd pathbbox 1 add cvi 31 or/_ury
xd 1 add cvi 31 or/_urx xd 1 sub cvi -32 and/_lly xd 1 sub cvi -32
and/_llx xd/_dH _ury _lly sub def/_dW _urx _llx sub def}bd/xp{/_pbs
xd _row _pbs mul _pbs getinterval/_s xd 0 _pbs 8 _pbs div dup _sW add
1 sub exch div cvi dup/str exch string def _pbs sub{str exch _s putinterval}for
str}bd/pbc{/__t save def{_m setmatrix}if 0 eq{eoclip}{clip}ifelse/_row
0 def bbx _llx _lly tr _t _dW mul ceiling cvi 8 add dup 8 mod sub dup/_sW
xd _t _dH mul ceiling cvi 8 add dup 8 mod sub _pg 0 ne{8 _im{_pr 8
xp}{_pg 8 xp}{_pb 8 xp/_row _row 1 add 8 mod def}true 3 colorimage}{1
_o eq{gs _br _bg _bb srgbc 2 copy true _im{_pb 1 xp/_row _row 1 add
8 mod def}imagemask gr}if _R _G _B srgbc false _im{_pb 1 xp/_row _row
1 add 8 mod def}imagemask}ifelse __t restore}bd
/hp{dup/_hb xd -1 eq{/_pb 0 def}{1 add/_pq xd/_pb 1 def}ifelse}bd/vert{X0
_w X1{dup Y0 M Y1 L stroke}for}bd/horz{Y0 _w Y1{dup X0 exch M X1 exch
L stroke}for}bd/fdiag{X0 _w X1{Y0 M X1 X0 sub dup l stroke}for Y0 _w
Y1{X0 exch M Y1 Y0 sub dup l stroke}for}bd/bdiag{X0 _w X1{Y1 M X1 X0
sub dup neg l stroke}for Y0 _w Y1{X0 exch M Y1 Y0 sub dup neg l stroke}for}bd/AU{1
add cvi 31 or}bd/AD{1 sub cvi -32 and}bd/hf{pathbbox 1 add cvi 31 or/Y1
xd 1 add cvi 31 or/X1 xd 1 sub cvi -32 and/Y0 xd 1 sub cvi -32 and/X0
xd 2 w [] 0 setdash/_w _rz 20 div 8 div 2 mul _pq div round 8 mul def
cp _hb 0 eq{horz}if _hb 1 eq{vert}if _hb 2 eq{fdiag}if _hb 3 eq{bdiag}if
_hb 4 eq{horz vert}if _hb 5 eq{fdiag bdiag}if gr}bd
/sp{dup 0 eq{/_pb xd}{255 div/_gr xd/_pb xd 100 div/_f xd}ifelse}bd/B{/_t2
xd/_t3 xd _pb _t2 _t3 8 idiv add get 1 7 _t3 8 mod sub bitshift and
0 ne}bd/SF{2{1 add 4 mul cvi exch}repeat B{0}{1}ifelse}bd/pb{_gr setgray
_f _a/SF load setscreen}bd
/SCA{10 div dup dup/_a xd/_A xd/__na xd currentscreen/__p xd/__a xd/__f
xd __f __na/__p load setscreen}bd
/SCF{100 div dup dup/_F xd/__nf xd currentscreen/__p xd/__a xd/__f
xd __nf __a/__p load setscreen}bd
/makeOblique 0 def/bC 1 def/makeoutlinedict 7 dict def/mbf{makeoutlinedict
begin/uniqueid exch def/strokewidth exch def/newfontname exch def/basefontname
exch def/basefontdict basefontname findfont def/numentries basefontdict
maxlength 1 add def basefontdict/UniqueID known not{/numentries numentries
1 add def}if/outfontdict numentries dict def basefontdict{exch dup/FID
ne{exch outfontdict 3 1 roll put}{pop pop}ifelse}forall outfontdict/FontName
newfontname put outfontdict/PaintType 2 put outfontdict/StrokeWidth
strokewidth put outfontdict/UniqueID uniqueid put newfontname outfontdict
definefont pop end}bd/usf{findfont exch makeOblique 0 ne{/makeOblique
makeOblique sin neg def [1 0 makeOblique 1 0 0] matrix concatmatrix/makeOblique
0 def}if makefont setfont}bd/sf{/_nf xd FontDirectory _nf known{pop}{findfont
dup maxlength dict/_nfd xd{exch dup/FID ne{exch _nfd 3 1 roll put}{pop
pop}ifelse}forall _nfd dup/FontName _nf put/Encoding WE put _nf _nfd
definefont pop}ifelse _nf usf}bd/us{gs [] 0 setdash 0 setlinecap w
0 exch M _w _be add 0 l currentrgbcolor sC S gr}bd/ob{n/_h xd _cx _cy
M _w 0 l 0 _h l _w neg 0 l P tc fill}bd/st{dup/_be xd 0 ne{/_bc xd}if/_vt
xd dup stringwidth pop/_sw xd dup length dup 1 gt{1 sub}if/_sl xd _be
0 eq{_w _sw sub _sl div 0 _vt{exch}if 3 -1 roll ashow}{_be _bc div
0 _vt{exch}if 32 _w _sw sub _sl div 0 _vt{exch}if 6 -1 roll awidthshow}ifelse}bd/sb{gs
M 10 div rotate/_w xd currentpoint/_cy xd/_cx xd 0 ne{ob}if tc _cx
_cy 3 -1 roll add tr 0 0 M st dup 4 and 4 eq{4 sub/_dx _cx _px sub
def _dx neg 0 tr/_w _w _dx add def}if dup 0 eq{pop}{dup 1 eq{pop us}{2
eq{us}{dup 3 1 roll us us}ifelse}ifelse}ifelse _cx/_px xd gr}bd
/cl{gs n M dup neg 3 1 roll 0 l 0 exch l 0 l P clip n}bd/clx{gs n
_cm currentmatrix pop _m setmatrix n M dup neg 3 1 roll 0 l 0 exch
l 0 l P clip _cm setmatrix n}bd
/sbb{/_t save def _neg 0 eq{0}{1}ifelse setgray tr/_dW xd/_dH xd dup/_sW
xd 8 _bs div dup _sW add 1 sub exch div cvi/str exch string def exch
dup/_sH xd 3 -1 roll{/_im{[_sW _dW div 0 0 _sH neg _dH div 0 _sH]}bd}{/_im{[_sW
_dW div 0 0 _sH _dH div 0 0]}bd}ifelse _bs 1 eq{false}{_bs}ifelse _im{currentfile
str readhexstring pop}_bs 1 eq{imagemask}{image}ifelse _t restore}bd
/scbb{/_t save def tr/_dW xd/_dH xd dup/_sW xd 8 _bs div dup _sW add
1 sub exch div cvi 3 mul/str exch string def exch dup/_sH xd _bs 4
-1 roll{/_im{[_sW _dW div 0 0 _sH neg _dH div 0 _sH]}bd}{/_im{[_sW
_dW div 0 0 _sH _dH div 0 0]}bd}ifelse _im{currentfile str readhexstring
pop}false 3 colorimage _t restore}bd
/mr{2 copy -1 eq{_yp _ym 2 mul sub}{0}ifelse exch -1 eq{_xp _xm 2
mul sub}{0}ifelse exch tr scale}bd
/WE 256 array def StandardEncoding WE copy pop WE dup
0[/grave/acute/circumflex/tilde/macron/breve/dotaccent/dieresis/ring/cedilla
/hungarumlaut/ogonek/caron/dotlessi
]putinterval dup
39/quotesingle put dup  96/grave put dup 124/bar put dup
145[/quoteleft/quoteright/quotedblleft/quotedblright/bullet/endash/emdash]putinterval 
160[/space/exclamdown/cent/sterling/currency/yen/brokenbar/section/dieresis
/copyright/ordfeminine/guillemotleft/logicalnot/endash/registered/macron/ring
/plusminus/twosuperior/threesuperior/acute/mu/paragraph/periodcentered
/cedilla/onesuperior/ordmasculine/guillemotright/onequarter/onehalf
/threequarters/questiondown/Agrave/Aacute/Acircumflex/Atilde/Adieresis/Aring
/AE/Ccedilla/Egrave/Eacute/Ecircumflex/Edieresis/Igrave/Iacute/Icircumflex
/Idieresis/Eth/Ntilde/Ograve/Oacute/Ocircumflex/Otilde/Odieresis/circumflex
/Oslash/Ugrave/Uacute/Ucircumflex/Udieresis/Yacute/Thorn/germandbls/agrave
/aacute/acircumflex/atilde/adieresis/aring/ae/ccedilla/egrave/eacute
/ecircumflex/edieresis/igrave/iacute/icircumflex/idieresis/eth/ntilde/ograve
/oacute/ocircumflex/otilde/odieresis/tilde/oslash/ugrave/uacute/ucircumflex
/udieresis/yacute/thorn/ydieresis]putinterval
/TTEV 256 array def TTEV
0 [/G00/G01/G02/G03/G04/G05/G06/G07/G08/G09/G0a/G0b/G0c/G0d/G0e/G0f
   /G10/G11/G12/G13/G14/G15/G16/G17/G18/G19/G1a/G1b/G1c/G1d/G1e/G1f
   /G20/G21/G22/G23/G24/G25/G26/G27/G28/G29/G2a/G2b/G2c/G2d/G2e/G2f
   /G30/G31/G32/G33/G34/G35/G36/G37/G38/G39/G3a/G3b/G3c/G3d/G3e/G3f
   /G40/G41/G42/G43/G44/G45/G46/G47/G48/G49/G4a/G4b/G4c/G4d/G4e/G4f
   /G50/G51/G52/G53/G54/G55/G56/G57/G58/G59/G5a/G5b/G5c/G5d/G5e/G5f
   /G60/G61/G62/G63/G64/G65/G66/G67/G68/G69/G6a/G6b/G6c/G6d/G6e/G6f
   /G70/G71/G72/G73/G74/G75/G76/G77/G78/G79/G7a/G7b/G7c/G7d/G7e/G7f
   /G80/G81/G82/G83/G84/G85/G86/G87/G88/G89/G8a/G8b/G8c/G8d/G8e/G8f
   /G90/G91/G92/G93/G94/G95/G96/G97/G98/G99/G9a/G9b/G9c/G9d/G9e/G9f
   /Ga0/Ga1/Ga2/Ga3/Ga4/Ga5/Ga6/Ga7/Ga8/Ga9/Gaa/Gab/Gac/Gad/Gae/Gaf
   /Gb0/Gb1/Gb2/Gb3/Gb4/Gb5/Gb6/Gb7/Gb8/Gb9/Gba/Gbb/Gbc/Gbd/Gbe/Gbf
   /Gc0/Gc1/Gc2/Gc3/Gc4/Gc5/Gc6/Gc7/Gc8/Gc9/Gca/Gcb/Gcc/Gcd/Gce/Gcf
   /Gd0/Gd1/Gd2/Gd3/Gd4/Gd5/Gd6/Gd7/Gd8/Gd9/Gda/Gdb/Gdc/Gdd/Gde/Gdf
   /Ge0/Ge1/Ge2/Ge3/Ge4/Ge5/Ge6/Ge7/Ge8/Ge9/Gea/Geb/Gec/Ged/Gee/Gef
   /Gf0/Gf1/Gf2/Gf3/Gf4/Gf5/Gf6/Gf7/Gf8/Gf9/Gfa/Gfb/Gfc/Gfd/Gfe/Gff
] putinterval
end
MGXPS_2.8_EPS begin
0 0 1 0 3400 4400 400 -1 0 0 1 Cne
4280 3240 0 0 clpg
0 0 0 sc
/_o 0 def
255 255 255 bc
[5 13 ]0 2 d
0 529 M 219 227 483 54 661 10 bz
838 -33 887 64 1070 278 bz
S
[]0 0 d
1777 860 M 1101 435 L
1101 435 1100 448 1083 448 bz
1045 448 989 353 1070 278 bz
S
[5 13 ]0 2 d
391 1025 M 618 1264 847 1365 1117 1363 bz
1362 1359 1561 1183 1777 860 bz
S
[]0 0 d
304 744 M 234 842 125 889 80 897 bz
34 903 18 862 31 778 bz
S
1243 525 M 1295 441 1335 375 1362 228 bz
S
2 529 M 141 628 385 807 510 888 bz
595 782 774 587 749 569 bz
711 536 518 753 509 889 bz
447 964 392 1025 392 1025 bz
S
6 w
10433 1000 div ml
305 744 M 428 557 642 431 764 432 bz
876 434 947 595 1038 624 bz
1128 650 1194 614 1243 525 bz
S
1 w
644 458 7 7 E
0 0 0 fc
0 2 p
end
EndOfTheIncludedPostscriptMagicCookie

\closepsdump

\usepackage{ifthen,amstex,rotating,calc,epic}
\usepackage[dvips]{epsfig}
\usepackage{epic}

\allowdisplaybreaks

\makeatletter
\renewcommand{\thesection}{\arabic{section}}

\renewcommand{\@seccntformat}[1]{\csname the#1\endcsname.\hspace{1em}}
\renewcommand\section{\@startsection{section}{1}{\z@}
                                   {-3.5ex \@plus -1ex \@minus -.2ex}
                                   {2.3ex \@plus.2ex}
                                   {\normalfont\normalsize}}
\renewcommand\subsection{\@startsection{subsection}{2}{\z@}
                                     {-3.25ex\@plus -1ex \@minus -.2ex}
                                     {1.5ex \@plus .2ex}
                                     {\normalfont\normalsize}}
\renewcommand{\l@section}[2]{
  \ifnum \c@tocdepth >\z@
    \addpenalty{\@secpenalty}
    \addvspace{\p@}
    \setlength\@tempdima{1.5em}
    \begingroup
      \parindent \z@ \rightskip \@pnumwidth
      \parfillskip -\@pnumwidth
      \leavevmode \scshape
      \advance\leftskip\@tempdima
      \hskip -\leftskip
      #1\nobreak\hfil \nobreak\hbox to\@pnumwidth{\hss #2}\par
    \endgroup
  \fi}
\def\@maketitle{
  \newpage
  \null
  \vskip 2em
  \begin{center}
    {\sffamily \large \@title \par}
    \vskip 1.5em
      {\normalsize
      \lineskip .5em
        {\scshape \normalsize \@author}}
  \end{center}
  \par
  \vskip 1.5em}
\makeatother

\DeclareFontFamily{OT1}{rsfs}{\hyphenchar\font45 }
\DeclareFontShape{OT1}{rsfs}{m}{n}{ <-7> rsfs5 <7-10> rsfs7 <10-> rsfs10}{}
\DeclareMathAlphabet{\mathcurl}{OT1}{rsfs}{m}{n}

\newcommand\real{\mathbb{R}}%
\newcommand\dup{\mathrm{d}}%
\newcommand\poly{(-\Delta)^m\dir}%
\newcommand\power[1][\Omega]{H_{#1,m}}%
\newcommand\dir{|_{\operatorname{DIR}}}
\newcommand\finbox{~\hfill$\Box$}%
\newcommand\lowearlybox{\tag*{$\frac[0pt]{\displaystyle{\phantom{\big|}}}{%
  \displaystyle{\Box}}$}}%
\newcommand\slref[1]{\textup{\ref{#1}}}%
\newcommand\sleqref[1]{\textsl{\eqref{#1}}}%
\newcommand\earlybox{\tag*{$\Box$}}%
\newcommand\re{\operatorname{Re}}%
\newcommand\bd[1]{\mathbf{#1}}%
\newcommand\der[1][h]{\frac{\dup #1}{\dup\bd{y}}}%
\newcommand\numder[1][h]{1+\der[#1]\tr\der[#1]}%
\newcommand\bnumder[1][h]{\left(\numder[#1]\right)}%
\newcommand\matder[1][h]{I_{N-1}+\der[#1]\der[#1]\tr}%
\newcommand\bmatder[1][h]{\left(\matder[#1]\right)}%
\newcommand\sder[1][h]{\frac{\dup^2#1}{\dup\bd{y}^2}}%
\newcommand\dev[3]{\frac{\partial #1_{#2}}{\partial y_{#3}^{\phantom2}}}%
\newcommand\sdevd[2]{\frac{\partial^2#1}{\partial y^2_{#2}}}%
\newcommand\sdev[3]{\frac{\partial^2#1}{\partial y_{#2}^{\phantom2}\partial%
  y_{#3}^{\phantom2}}}%
\newcommand\tr{^{\mathrm{T}}}%
\newcommand\ip[2]{\langle #1,#2\rangle}%

\newcommand\afigwidth0%
\newcommand\afigheight0%
\newcommand\bfigwidth0%
\newcommand\bfigheight0%
\newlength\aeqspace%
\newlength\beqspace%
\newlength\border%
\newlength\tempeqspace%
\newcounter{tempactr}%
\newcounter{tempbctr}%
\newcounter{tempcctr}%
1%

\newcommand{\fig}[6]{%
  \setcounter{tempactr}{#6}%
  \setcounter{tempbctr}{200+\border/\unitlength}%
  \setcounter{tempcctr}{#6+\border/\unitlength}%
  \setlength\tempeqspace{#5}%
  \vskip\tempeqspace%
  \vskip\border%
  \centering%
  \begin{picture}(200,\value{tempactr})%
    \put(0,0){\epsfig{figure=#1.eps,%
                      width=200\unitlength}}%
    #3%
  \end{picture}%
  \vskip\border%
  \vskip\tempeqspace%
  \vskip\abovecaptionskip%
  \refstepcounter{figure}#2%
  \makebox[0pt][c]{\parbox[t]{1.5\textwidth}{\centering%
    \ifthenelse{\equal{#4}{}}{Figure \thefigure}{Figure \thefigure: #4}}}%
  \vskip\belowcaptionskip}%

\newcommand{\tab}[6]{%
  \setcounter{tempactr}{#6}%
  \setcounter{tempbctr}{200+\border/\unitlength}%
  \setcounter{tempcctr}{#6+\border/\unitlength}%
  \setlength\tempeqspace{#5}%
  \vskip\tempeqspace%
  \vskip\border%
  \centering%
  {\small%
  \begin{picture}(200,\value{tempactr})%
    \put(0,0){\begin{tabular}[b]{#2}%
      #3%
    \end{tabular}}%
  \end{picture}}%
  \vskip\border%
  \vskip\tempeqspace%
  \vskip\abovecaptionskip%
  \refstepcounter{table}#1%
  \makebox[0pt][c]{\parbox[t]{1.5\textwidth}{\centering%
    \ifthenelse{\equal{#4}{}}{Table \thetable}{Table \thetable: #4}}}%
  \vskip\belowcaptionskip}%

\newcommand{\dbledgram}[3]{{\setlength\border{0pt}
  \renewcommand\baselinestretch1\normalsize%
  \setlength\aeqspace{0pt}%
  \setlength\beqspace{0pt}%
  #1%
  \ifthenelse{\bfigheight > \afigheight}{%
    \setlength\aeqspace{\textwidth/200*(\bfigheight-\afigheight)}}{%
    \setlength\beqspace{\textwidth/200*(\afigheight-\bfigheight)}}%
  \begin{figure}[t]%
    \par%
    \noindent%
    \hfill%
    \setlength\unitlength{\textwidth/100*\afigwidth/200}%
    \begin{minipage}[t]{\textwidth/100*\afigwidth+2\border}%
      #2\aeqspace{200*\afigheight/\afigwidth}%
    \end{minipage}%
    \hfill%
    \setlength\unitlength{\textwidth/100*\bfigwidth/200}%
    \begin{minipage}[t]{\textwidth/100*\bfigwidth+2\border}%
      #3\beqspace{200*\bfigheight/\bfigwidth}%
    \end{minipage}%
    \hfill%
    \null%
  \end{figure}}}%

\newcommand{\sngldiagm}[2]{{%
  \setlength\border{0pt}%
  \renewcommand\baselinestretch1\normalsize%
  #1%
  \begin{figure}[t]%
    \par%
    \noindent%
    \hfill%
    \setlength\unitlength{\textwidth/100*\afigwidth/200}%
    \begin{minipage}[t]{\textwidth/100*\afigwidth+2\border}%
      #2{0pt}{200*\afigheight/\afigwidth}%
    \end{minipage}%
    \hfill%
    \null%
  \end{figure}}}%

\newcounter{theorem}[section]%
\renewcommand\thetheorem{\thesection.\arabic{theorem}}
\newenvironment{thm}[1]{%
    \refstepcounter{theorem}%
    \textsc{#1 \thetheorem.}{}%
    \begin{itshape}}%
  {\end{itshape}}%

\newenvironment{proof}{\textit{Proof.}{}}{}%

\newenvironment{note}{%
  \refstepcounter{theorem}%
  \textsc{Note \thetheorem.}{}}%
  {}%

\newenvironment{definition}[1][]{%
  \refstepcounter{theorem}%
  \ifthenelse{\equal{#1}{}}{\textsc{Definition \thetheorem.}{} }%
    {\textbf{Definition \thetheorem: (#1)}{} }%
  \begin{itshape}}%
  {\end{itshape}}%

\begin{document}

\begin{center}
\large
A RIEMANNIAN OFF-DIAGONAL HEAT KERNEL BOUND FOR\\
UNIFORMLY ELLIPTIC OPERATORS\\
\vspace{1cm}
\normalsize
M.P. OWEN\\
\vspace{2cm} 

\begin{minipage}{0.9\linewidth}
  \textsc{Abstract}.
  We find a Gaussian off-diagonal heat kernel estimate for uniformly
  elliptic operators with measurable coefficients acting on regions $\Omega 
  \subseteq\real^N$, where the order $2m$ of the operator satisfies $N<2m$. 
  The estimate is expressed using certain Riemannian-type metrics, and a 
  geometrical result is established allowing conversion of the estimate into 
  terms of the usual Riemannian metric on $\Omega$. Work of 
  Barbatis~\cite{Barb:U} is applied to find the best constant in this 
  expression.\par
\end{minipage}
\end{center}

\section{INTRODUCTION}

Let $H$ be a differential operator with quadratic form $\phantom{\mathbf{}}$%
%
%
\begin{equation}\label{eqn:ellipt2}
  Q(f)=\int_\Omega\sum
  \begin{Sb}
    |i|\le m \\
    |j|\le m
  \end{Sb}
  a_{i,j}(x)D^if(x)\overline{D^jf(x)}\dup^Nx,
\end{equation}
where $a_{i,j}(x)=\overline{a_{j,i}(x)}$ are complex-valued bounded 
measurable functions on a region $\Omega\subseteq\real^N$. Dirichlet boundary 
conditions are imposed upon $H$ by restricting the domain of the quadratic 
form to be the Sobolev space $W_0^{m,2}(\Omega)$, which is the closure of
$C_c^\infty(\Omega)$ in the Hilbert space $W^{m,2}(\Omega)$. Here $C_c^\infty(
\Omega)$ denotes the space of smooth, compactly supported functions on $\Omega$
and $W^{m,2}(\Omega)$ is the space of all functions $f\in L^2(\Omega)$ whose
weak derivatives $D^\alpha f$ lie in $L^2(\Omega)$ for all multi-indices
$\alpha$ such that $|\alpha|\le m$. It is equipped with the inner product
\begin{equation}
  \ip fg_{m,2}:=\sum_{|\alpha|\le m}\ip{D^\alpha f}{D^\alpha g}_2.
\end{equation}
In the special case where the coefficients $a_{i,j}$ are constant and chosen 
so that 
\[ \sum_{|i|=|j|=m}a_{i,j}\xi^{i+j}=|\xi|^{2m} \]
and $a_{i,j}=0$ whenever $|i|+|j|\le2m-1$,
the associated operator $H=:H_{\Omega,m}$ is the polyharmonic operator $\poly$
of the region $\Omega$. We denote this particular quadratic form by $Q_m$.

Throughout this paper we shall assume that the coefficients $a_{i,j}(x)$ are 
chosen in such a way that $Q$ satisfies the G\r{a}rding 
inequality
\begin{equation}\label{eqn:garding}
  \lambda Q_m(f)-c\|f\|_2^2\le Q(f)\le\mu Q_m(f)+d\|f\|_2^2,
\end{equation}
where $0<\lambda\le\mu$ and $c$, $d$ are non-negative constants. We then say
that $H$ is uniformly elliptic. Quadratic forms satisfying the G\r{a}rding 
inequality are closed on the domain $W^{m,2}(\Omega)$. For a more detailed 
account of uniformly elliptic operators see~\cite{Davi3:1995}.

\begin{note}\label{def:homog}
  If the sum in equation~\eqref{eqn:ellipt2} is only taken over non-negative 
  multi-indices $i$, $j$ with $|i|=|j|=m$ then the operator is said to be 
  homogeneous of order $2m$. If this is the case then we may set $c=d=0$ in the
  G\r{a}rding inequality~\eqref{eqn:garding}.
\end{note}

\begin{thm}{Lemma}\label{thm:unifkerbd}
  For $N<2m$ the operator $H$ has a heat kernel $K(t,x,y)$ which satisfies
  \begin{equation}
    |K(t,x,y)|\le ct^{-N/2m}e^t
  \end{equation}
  for all $x,y\in\Omega$ and all $t>0$.
\end{thm}

\begin{proof}
  See~\cite[Corollary 15 and Lemma 17]{Davi3:1995}. If the operator $H$ is 
  homogeneous then the bound is valid even without the term $e^t$. \finbox
\end{proof}


In~\cite[Lemma 19]{Davi3:1995}, Davies obtains the pointwise heat kernel bound
\begin{equation}\label{eqn:davipt}
  |K(t,x,y)|\leq c_1t^{-N/2m}\exp[-c_2|y-x|^{2m/(2m-1)}t^{-1/(2m-1)}+kt],
\end{equation}
for all $t>0$ and $x,y\in\real^N$, where $c_1,c_2,k$ are positive constants. If
the operator is homogeneous of order $2m$ then the constant $k$ may be set to 
zero. In~\cite{Barb+Davi:1996}, Barbatis and Davies find the sharp constant $c_2$ 
for this expression in the following sense: If $H$ is uniformly elliptic, 
homogeneous of order $2m$, and satisfies the G\r{a}rding inequality
\begin{equation}
  Q_m(f)\le Q(f)\le\mu Q_m(f)
\end{equation}
then
\begin{equation}
  |K(t,x,y)|\leq c_\epsilon t^{-N/2m}\exp[-(\sigma_m-O(\mu-1)-\epsilon)|y-x|^{
  2m/(2m-1)}t^{-1/(2m-1)}],
\end{equation}
where $\epsilon>0$ and
\begin{equation}
  \sigma_m=(2m-1)(2m)^{-2m/(2m-1)}\sin[\pi/(4m-2)].
\end{equation}
This is shown to be sharp by considering the case $H=\power[\real^N]$.

The Euclidean metric $d_0(x,y):=|y-x|$ is a relatively weak and unnatural way 
of expressing the heat kernel bound for non-convex regions. As an example, for
a horse-shoe shaped region whose extremities are touching, the Euclidean 
distance $d_0(x,y)$ can be made arbitrarily small for internally distant points
$x,y\in\Omega$, rendering the heat kernel bound useless at these points. 

The aim of this paper is to find heat kernel bounds which are given 
in terms of the Riemannian metric $d_g$ (see definition~\eqref{eqn:riemann}),
instead of the Euclidean metric $d_0$. This would improve 
the original bound; for a horse-shoe shaped region which touches itself or 
even overlaps with itself one may choose $x$ and $y$ to make $d_0(x,y)$ 
arbitrarily small whilst $d_g(x,y)$ remains large.
In order to find bounds involving the Riemannian metric, we first find bounds 
involving Riemannian-type metrics.


\begin{definition}\label{def:rietyp}
  For $\beta>0$ define the Riemannian-type metrics $d_{m,\beta}:\Omega^2
  \rightarrow\real_+$ by
  \begin{equation}
    d_{m,\beta}(x,y)=\sup\{\phi(y)-\phi(x):\phi\in\mathcurl E_{m,\beta}\},
  \end{equation}
  where $\mathcurl E_{m,\beta}$ denotes the set of all bounded real valued 
  smooth functions $\phi$ on a region $\Omega$ such that
  \begin{equation}
    \|\nabla\phi\|_\infty\leq1\quad\text{and}\quad\|D^i\phi\|_\infty\le\beta^{|
    i|-1}
  \end{equation}
  for all non-negative multi-indices $i$ such that $2\le|i|\le m$.
\end{definition}
 
We show, for an arbitrary uniformly elliptic operator whose 
coefficients need only be measurable, that
\begin{equation}\label{eqn:mine}
  |K(t,x,y)|\leq c_1t^{-N/2m}\exp[-c_2d_{m,\beta}(x,y)^{2m/(2m-1)}t^{-1/(2m-1)}
  +k(1+\beta^{2m})t]
\end{equation}
where $N<2m$. Note that since $d_{m,0}$ is the Euclidean metric $d_0$, setting 
$\beta=0$ in this equation retrieves equation~\eqref{eqn:davipt}.


If $\Omega$ is a region with $C^2$ boundary and radii of curvature uniformly 
bounded below by $r$ then there exists a constant $K$, dependent only upon
$m$ and $N$, such that for $\beta\geq4K/r$,
  \begin{equation*}
    \left(1-\sqrt{\frac{K}{\beta r}}\right)d_g\leq d_{m,\beta}\leq d_g,
  \end{equation*}
where $d_g$ is the standard Riemannian distance
\begin{equation}\label{eqn:riemann}
  d_g(x,y):=\inf\{l(\gamma):\gamma(0)=x,\gamma(1)=y,\gamma\subseteq\overline
  \Omega,\gamma\text{ is cts and piecewise }C^1\},
\end{equation}
and where $l(\gamma)$ is the length of the path $\gamma$.
See Theorem~\ref{thm:metricmain}. This feature allows us to convert the 
off-diagonal bound~\eqref{eqn:mine} into the bound
\begin{equation}\label{eqn:reimbd}
  |K(t,x,y)|\leq c_1t^{-N/2m}\exp[-c_2d_g(x,y)^{2m/(2m-1)}t^{-1/(2m-1)}+kt]
\end{equation}
involving the Riemannian metric.

\begin{note}
  As a special case of~\cite[Theorem 3.2.7 and Corollary 3.2.8]{Davi:1989} 
  Davies has the bound
  \begin{equation}
    0\leq K(t,x,y)\leq c_\delta t^{-N/2}\exp\{-d_g(x,y)^2/4(1+\delta)t\},
  \end{equation}
  valid for all $N$, on the heat kernel $K(t,x,y)$ of the Dirichlet Laplacian 
  on a region $\Omega\subseteq\real^N$. The regularity of the boundary of 
  $\Omega$, required for higher order operators is a genuine feature of the 
  order.\finbox
\end{note}


A more natural class of metrics for the determination of heat kernel bounds is,
in a certain sense, the class of Finsler-type metrics $d_{a,M}$ induced by the
operator itself. Here $a$ denotes the principal symbol
\begin{equation}\label{eqn:symbol}
  a(x,\xi)=\sum
  \begin{Sb}
    |i|=m \\
    |j|=m
  \end{Sb}
  a_{i,j}(x)\xi^{i+j}
\end{equation}
of the operator, and $M$ is a positive constant.

\begin{definition}\label{def:finsler}
  The Finsler-type metrics are defined by 
  \begin{equation}
    d_{a,M}(x,y)=\sup\{\phi(y)-\phi(x):\phi\in\mathcurl{F}_{a,M}\},
  \end{equation}
  where $\mathcurl{F}_{a,M}$ denotes the set of all bounded real-valued 
  smooth functions $\phi$ on $\Omega$ such that
  \begin{equation}
    a(x,\nabla\phi(x))\leq1\quad\text{and}\quad\|D^i\phi\|_\infty\le M
  \end{equation}
  for all non-negative multi-indices $i$ such that $2\le|i|\le m$. 
\end{definition}

Note that the Riemannian-type metrics $d_{m,\beta}$ are 
similar to the Finsler-type metrics induced by the polyharmonic operator 
$\poly$, whose symbol is
\begin{equation*}
  a(x,\xi)=|\xi|^{2m}.
\end{equation*}

Under the following assumptions, Barbatis~\cite{Barb:U} uses Finsler-type 
metrics to express a sharp heat kernel bound.

\begin{thm}{Assumptions}\label{thm:barbassump}
  \begin{enumerate}
    \item $H$ is uniformly elliptic and homogeneous of order $2m>N$\textup{\/;}
    \item The coefficients $a_{i,j}$ lie in the Sobolev space $W^{m,\infty}(
          \Omega)$\textup{\/;}
    \item The symbol of $H$ is strongly convex (see 
          Definition~\slref{def:strconv} below). 
  \end{enumerate}
\end{thm}

\begin{definition}\label{def:strconv}
  For $|k|=2m$ define
  \begin{equation*}
    \alpha_k(x)=\frac{k!}{(2m)!}\sum
    \begin{Sb}
      |i|=|j|=m \\
      i+j=k
    \end{Sb}
    a_{i,j}(x).
  \end{equation*}
  We may rewrite the principal symbol of $H$ as
  \begin{equation*}
    a(x,\xi)=\sum_{|k|=2m}\frac{(2m)!}{k!}a_k(x)\xi^k
  \end{equation*}
  We say that the symbol $a(x,\xi)$ is strongly convex if the quadratic form
  \begin{equation*}
    \Gamma(x,\zeta)=\sum
    \begin{Sb}
      |p|=m \\
      |q|=m
    \end{Sb}
    a_{p+q}(x)\zeta_p\zeta_q
  \end{equation*}
  is non-negative for each $x\in\Omega$, where $\zeta=(\zeta_p)_{|p|=m}\in
  \real^{\nu}$ and
  \[ \nu=\left(\frac[0pt]{n+m-1}{n-1}\right) \]
  is the number of distinct multi-indices $p$ with $|p|=m$.
\end{definition}

Barbatis proves that, under these assumptions, we have the heat kernel bound
\begin{equation}\label{eqn:barbs}
  |K(t,x,y)|\le c_\delta t^{-N/2m}\exp[-(\sigma_m-\delta)d_{a,M}(x,y)^{2m/(
  2m-1)}t^{-1/(2m-1)}]
\end{equation}
for $M$ large and for $t/d_{a,M}\le T_{\delta,M}$. We shall apply this result 
under the same assumptions, to find the sharp constant $c_2$ in the 
bound~\eqref{eqn:reimbd}.

\section{THE HEAT KERNEL BOUNDS}\label{sec:kerbd}

We shall use the technique found in~\cite{Davi3:1995} of twisting the operator 
$H$ to define $H_{\alpha\phi}=e^{\alpha\phi}He^{-\alpha\phi}$, which has 
quadratic form 
\begin{equation}
    Q_{\alpha\phi}(f)=\int_\Omega\sum
    \begin{Sb}
      |i|\le m \\
      |j|\le m
    \end{Sb}   
    a_{i,j}(x)\{e^{\alpha\phi}D^ie^{-\alpha\phi}f(x)\}\{\overline{e^{-\alpha
    \phi}D^je^{\alpha\phi}f(x)}\}\dup^Nx.
\end{equation}
Here $\alpha>0$ and $\phi\in\mathcurl{E}_{m,\beta}$ for some $\beta>0$. The 
following proposition is a generalisation of~\cite[Lemmas 1 and 2]{Davi3:1995},
treating derivatives of the function $\phi$ more delicately.

\begin{thm}{Proposition}
  The twisted quadratic form satisfies the inequality
  \begin{equation}\label{eqn:twist}
    |Q_{\alpha\phi}(f)-Q(f)|\le\epsilon Q(f)+c_\epsilon(1+\alpha^{2m}+\beta^{2m
    })\|f\|_2^2,
  \end{equation}
  where $\epsilon$ may be taken arbitrarily small.
\end{thm}

\begin{proof}
  Each term in $Q_{\alpha\phi}(f)$ may be expanded using formulae of the type 
  \begin{equation}
    e^{\alpha\phi}D^ie^{-\alpha\phi}f=D^if+\sum c_k\big(\prod_{r=1}^p\alpha D^{
    k_r}\phi\big)D^{k_0}f,
  \end{equation}
  where the sum is taken over all integers $p$ and non-negative multi-indices 
  $k_0,\dots,k_p$ such that
  \begin{equation}
    \sum_{r=0}^pk_r=i,\quad k_0\neq i,\quad\text{and}\quad k_1,\dots k_p\neq0.
  \end{equation}
  Combining these terms,
  \begin{equation}
    Q_{\alpha\phi}(f)=Q(f)+\int_\Omega\sum\nolimits'c_{k,l}(x)\alpha^p\big(
    \prod_{r=1}^pD^{k_r}\phi\big)D^{k_0}f.\alpha^q\big(\prod_{s=1}^qD^{l_s}\phi
    \big)\overline{D^{l_0}f}\dup^Nx,
  \end{equation}
  where the sum $\sum\nolimits'$ is taken over all integers $p,q$ and 
  non-negative multi-indices $k_0,\dots,k_p,l_0,\dots,l_q$ such that 
  \begin{equation}
    \sum_{r=0}^p|k_r|\le m,\quad\sum_{s=0}^q|l_s|\le m,\quad\text{and}\quad|k_
    0|+|l_0|\le2m-1.
  \end{equation}
  Hence 
  \begin{equation}
    |Q_{\alpha\phi}(f)-Q(f)|\le\sum\nolimits'\|c_{k,l}\|_\infty\|\alpha^p\big(
    \prod_{r=1}^p D^{k_r}\phi\big)D^{k_0}f\|_2\|\alpha^q\big(\prod_{s=1}^qD^{l_
    s}\phi\big)D^{l_0}f\|_2.
  \end{equation}
  Now
  \begin{align*}
    \|\alpha^p\big(\prod_{r=1}^pD^{k_r}\phi\big)D^{k_0}f\|_2&\le\big(\prod_{r=1
    }^p\|D^{k_r}\phi\|_\infty\big)\|\alpha^pD^{k_0}f\|_2 \\
    &\le\big(\prod_{r=1}^p\beta^{|k_r|-1}\big)\|\alpha^pD^{k_0}f\|_2 \\
    &=\|\alpha^p\beta^{s-p-|k_0|}D^{k_0}f\|_2,
  \end{align*}
  where $s:=\sum_{r=0}^pk_r\le m$, and so
  \begin{align}
    \|\alpha^p\big(\prod_{r=1}^pD^{k_r}\phi\big)D^{k_0}f\|_2^2&\le\int_{\real^N
    }\alpha^{2p}\beta^{2s-2p-2|k_0|}(i\xi)^{2k_0}|\hat{f}(\xi)|^2\dup^N\xi
    \notag\\
    &\le\int_{\real^N}\alpha^{2p}\beta^{2s-2p-2|k_0|}|\xi|^{2|k_0|}|\hat{f}(\xi
    )|^2\dup^N\xi\notag \\
    &\le\int_{\real^N}[\epsilon|\xi|^{2s}+c_\epsilon(\alpha^{2s}+\beta^{2s})]|
    \hat{f}(\xi)|^2\dup^N\xi \label{eqn:earlystop}\\
    &\le\int_{\real^N}[\epsilon|\xi|^{2m}+c_\epsilon(\alpha^{2m}+\beta^{2m}+1)]
    |\hat{f}(\xi)|^2\dup^N\xi \label{eqn:normalstop}
  \end{align}
  where $\epsilon$ may be arbitrarily small in~\eqref{eqn:earlystop} provided
  $|k_0|\le s-1$ and arbitrarily small in~\eqref{eqn:normalstop} provided $|k_
  0|\le m-1$. Each term 
  in $\sum\nolimits'$ has at least one of $|k_0|\le m-1$ or $|l_0|\le m-1$ and 
  so is dominated by
  \begin{equation*}
    \epsilon Q_m(f)+c_\epsilon(1+\alpha^{2m}+\beta^{2m})\|f\|_2^2.
  \end{equation*}
  Using the G\r{a}rding inequality~\eqref{eqn:garding} this in turn is 
  dominated by
  \begin{equation*}
    \epsilon Q(f)+c_\epsilon(1+\alpha^{2m}+\beta^{2m})\|f\|_2^2. \earlybox
  \end{equation*}
\end{proof}

\begin{note}\label{not:homog}
  If $H$ is homogeneous (see Definition~\ref{def:homog}), then $2s=2m$ 
  and the above proof finishes with equation~\eqref{eqn:earlystop} instead 
  of~\eqref{eqn:normalstop}. This yields the twisted form inequality
  \begin{equation}
    |Q_{\alpha\phi}(f)-Q(f)|\le\epsilon Q(f)+c_\epsilon(\alpha^{2m}+\beta^{2m
    })\|f\|_2^2,
  \end{equation}
  instead of equation~\eqref{eqn:twist}. This will induce a corresponding 
  change in inequalities~\eqref{eqn:unifsemi1}, \eqref{eqn:unifsemi2} and 
  \eqref{eqn:unifker}.\finbox
\end{note}

\begin{thm}{Lemma}
  There exist positive constants $c,k>0$ such that
  \begin{equation}\label{eqn:unifsemi1}
    \|e^{-H_{\alpha\phi}t}\|\leq\exp[k(1+\alpha^4+\beta^4)t]
  \end{equation}
  and
  \begin{equation}\label{eqn:unifsemi2}
    \|H_{\alpha\phi}e^{-H_{\alpha\phi}t}\|\leq ct^{-1}\exp[k(1+\alpha^4+\beta^4
    )t].
  \end{equation}
  for all $t>0$, $\alpha,\beta>0$, and $\phi\in\mathcurl{E}_{m,\beta}$.
\end{thm}

\begin{proof}
  See~\cite[Lemmas 6 and 7]{Davi3:1995}.\finbox
\end{proof}

\begin{thm}{Lemma}
  Let $d_{m,\beta}$ be the Riemannian-type metrics of 
  Definition~\slref{def:rietyp}. There exist positive constants $c_1,c_2,k$ 
  such that
  \begin{equation}\label{eqn:unifker}
    |K(t,x,y)|\leq c_1t^{-N/2m}\exp[-c_2d_{m,\beta}(x,y)^{2m/(2m-1)}t_{\phantom
    {g}}^{-1/(2m-1)}+k(1+\beta^{2m})t]
  \end{equation}
  for all $\beta,t>0$ and all $x,y\in\Omega$.
\end{thm}

\begin{proof}
  Put $f_t=\exp[-H_{\alpha\phi}t]f$. Then
  \begin{align*}
    \|f_t\|_\infty&\leq c\|(-\Delta)^{m/2}f_t\|_2^{N/2m}\|f_t\|_2^{1-N/2m} \\
    &\le cQ(f_t)^{N/4m}\|f_t\|_2^{1-N/2m} \\
    &\le c\{\re Q_{\alpha\phi}(f_t)+(1+\alpha^{2m}+\beta^{2m})\|f_t\|_2^2\}^{
    N/4m}\|f_t\|_2^{1-N/2m} \\
    &\le c\{\|H_{\alpha\phi}f_t\|_2\|f_t\|_2+(1+\alpha^{2m}+\beta^{2m})\|f_t\|
    _2^2\}^{N/4m}\|f_t\|_2^{1-N/2m} \\
    &\leq c\{t^{-1}+(1+\alpha^{2m}+\beta^{2m})\}^{N/4m}\exp[k(1+\alpha^{2m}+
    \beta^{2m})t]\|f\|_2 \\
    &\leq ct^{-N/4m}\exp[k(1+\alpha^{2m}+\beta^{2m})t]\|f\|_2.
  \end{align*}
  Therefore
  \begin{equation}
    \|\exp[-H_{\alpha\phi}t]\|_{\infty,2}\leq ct^{-N/4m}\exp[k(1+\alpha^{2m}+
    \beta^{2m})t].
  \end{equation}
  By duality,
  \begin{equation}
    \|\exp[-H_{\alpha\phi}t]\|_{\infty,1}\leq ct^{-N/2m}\exp[k(1+\alpha^{2m}+
    \beta^{2m})t].
  \end{equation}
  But $\exp[-H_{\alpha\phi}t]$ has kernel
  \begin{equation}
    K_{\alpha\phi(x)}(t,x,y)=e^{\alpha\phi(x)}K(t,x,y)e^{-\alpha\phi(x)}
  \end{equation}
  for all $t>0$ and $x,y\in\Omega$. Equivalently,
  \begin{equation}
    |K(t,x,y)|\leq ct^{-N/2m}\exp[\alpha(\phi(y)-\phi(x))+k(1+\alpha^{2m}+\beta
    ^{2m})t].
  \end{equation}
  Taking the infimum over all $\phi\in\mathcurl{E}_{m,\beta}$ in this bound, 
  we see that
  \begin{equation*}
    |K(t,x,y)|\leq ct^{-N/2m}\exp[-\alpha d_{m,\beta}(x,y)+k(1+\alpha^{2m}+
    \beta^{2m})t].
  \end{equation*}
  Optimising with respect to $\alpha$ gives
  \begin{equation*}
    |K(t,x,y)|\leq ct^{-N/2m}\exp[-k'd_{m,\beta}(x,y)^{2m/(2m-1)}t_{\phantom{g}}
    ^{-1/(2m-1)}+k(1+\beta^{2m})t]. \earlybox
  \end{equation*}
\end{proof}

\begin{thm}{Theorem}
  Let $\Omega\subseteq\real^N$ be a region whose boundary 
  is $C^2$, with radii of curvature uniformly bounded below, and suppose that
  $N<2m$. There exist positive constants $c_1$, $c_2$, $k$ such that
  \begin{equation}\label{eqn:geodbd}
    |K(t,x,y)|\leq c_1t^{N/2m}\exp[-c_2d_g(x,y)^{2m/(2m-1)}t_{\phantom{g}}^{
    -1/(2m-1)}+kt]
  \end{equation}
\end{thm}

\begin{proof}
  This follows by setting $\beta=4K/r$ and applying the main result, 
  Theorem~\ref{thm:metricmain}, of Section~\ref{sec:metrics}.\finbox
\end{proof}

\begin{note}
  The statement~\eqref{eqn:geodbd} is equivalent to the existence of positive
  constants $c_1,c_2,T$ such that for $t\leq Td_g(x,y)$
  \begin{equation}
    |K(t,x,y)|\leq c_1t^{-N/2m}\exp[-c_2d_g(x,y)^{2m/(2m-1)}t_{\phantom{g}}^{
    -1/(2m-1)}]
  \end{equation}
  See Lemma~\ref{thm:conequiv}.\finbox
\end{note}

\section{SHARP CONSTANTS FOR THE HEAT KERNEL BOUND}

\begin{thm}{Lemma}\label{thm:conequiv}
  Let $H$ be a uniformly elliptic operator acting in $L^2(\Omega)$ where 
  $\Omega\subseteq\real^N$ and $N<2m$. Let $K(t,x,y)$ be the 
  heat kernel of $H$, and let $c_2$ be fixed. The following conditions are 
  equivalent:
  \begin{enumerate}
    \item For all $\epsilon>0$ there exist positive constants $c_1,T$ such that 
          \begin{equation}
            |K(t,x,y)|\leq c_1t^{-N/2m}\exp[-(c_2-\epsilon)d_g(x,y)^{2m/(2m-1)}
            t_{\phantom{g}}^{-1/(2m-1)}]
          \end{equation}
          for $t/d_g(x,y)<T$.
    \item For all $\epsilon>0$ there exist positive constants $c_1,k$ such that
          \begin{equation}
            |K(t,x,y)|\leq c_1t^{-N/2m}\exp[-(c_2-\epsilon)d_g(x,y)^{2m/(2m-1)}
            t_{\phantom{g}}^{-1/(2m-1)}+kt]
          \end{equation}
          for $t>0$ and $x,y\in\Omega$.
  \end{enumerate}
\end{thm}

\begin{proof}
  Suppose that (i) holds. Let $k=(c_2-\epsilon)T^{-2m/(2m-1)}+1$. Then for $t/d
  _g\geq T$
  \begin{equation*}
    -(c_2-\epsilon)d_g^{2m/(2m-1)}t_{\phantom{g}}^{-1/(2m-1)}\geq-kt+t
  \end{equation*}
  and hence using Lemma~\ref{thm:unifkerbd},
  \begin{align*}
    |K(t,x,y)|&\leq ct^{-N/2m}\exp[t] \\
    &\leq ct^{-N/2m}\exp[-(c_2-\epsilon)d_g^{2m/(2m-1)}t_{\phantom{g}}^{-1/(2m-
    1)}+kt].
  \end{align*}
  Conversely, suppose that (ii) holds. Let $T=(\epsilon/k)^{(2m-1)/2m}$. For $t
  /d_g(x,y)\leq T$ we have $kt\leq\epsilon d_g^{2m/(2m-1)}t_{\phantom{g}}^{-1/(
  2m-1)}$ and so
  \begin{align*}
  |K(t,x,y)|&\leq c_1t^{-N/2m}\exp[-c_2d_g^{2m/(2m-1)}t_{\phantom{g}}^{-1/(2m-
  1)}+\epsilon d_g^{2m/(2m-1)}t_{\phantom{g}}^{-1/(2m-1)}] \\
  &= c_1t^{-N/2m}\exp[-(c_2-\epsilon)d_g^{2m/(2m-1)}t_{\phantom{g}}^{-1/(2m-1)
  }]. \earlybox
  \end{align*}
\end{proof}

\begin{thm}{Theorem}
  Suppose that $\Omega$ is a region with $C^2$ boundary and whose radii of 
  curvature are bounded below. Suppose also that $H$ satisfies the 
  assumptions~\slref{thm:barbassump}. For $\epsilon>0$ there exist positive 
  constants $c_\epsilon$, $k_\epsilon$ such that
  \begin{equation}\label{eqn:bestcon}
    |K(t,x,y)|\le c_\epsilon t^{-N/2m}\exp[-(\sigma_m-\epsilon)\mu^{-1/(2m-1)}d
    _g(x,y)^{2m/(2m-1)}t^{-1/(2m-1)}+k_\epsilon t].
  \end{equation}
  Equivalently, for $\epsilon>0$ there exist positive constants $c_\epsilon$, 
  $T_\epsilon$ such that
  \begin{equation}
    |K(t,x,y)|\le c_\epsilon t^{-N/2m}\exp[-(\sigma_m-\epsilon)\mu^{-1/(2m-1)}d
    _g(x,y)^{2m/(2m-1)}t^{-1/(2m-1)}]
  \end{equation}
  for $t/d_g(x,y)<T_\epsilon$.
\end{thm}

\begin{proof}
  For $\epsilon>0$ fixed let $\beta\ge4K/r$ be large enough such that the 
  result~\eqref{eqn:barbs} of Barbatis~\cite{Barb:U} is valid for $M:=\mu^{-
  1/2m}\beta^{m-1}$, and such that
  \begin{equation*}
    (\sigma_m-\epsilon)\le(\sigma_m-\delta)\left(1-\sqrt{\frac{K}{\beta r}}
    \right).
  \end{equation*}
  The G\r{a}rding inequality implies that the symbol $a(x,\xi)$ defined in 
  equation~\eqref{eqn:symbol} satisfies
  \begin{equation*}
    \lambda|\xi|^{2m}\le a(x,\xi)\le\mu|\xi|^{2m}
  \end{equation*}
  for all $\xi\in\real^N$. Let $\phi\in\mathcurl{E}_{m,\beta}$ and define $\psi
  =\mu^{-1/2m}\phi$. Then
  \begin{equation*}
    a(x,\nabla\psi(x))=\mu^{-1}a(x,\nabla\phi(x))\le|\nabla\phi(x))|^{2m}\le1
  \end{equation*}
  and
  \begin{align*}
    \|D^i\psi\|_\infty&=\mu^{-1/2m}\|D^i\phi\|_\infty\le\mu^{-1/2m}\beta^{|i|-1
    }\le M,
  \end{align*}
  so $\psi\in\mathcurl{F}_{a,M}$. Hence
  \begin{align}\label{eqn:relbarb}
    \mu^{-1/2m}d_{m,\beta}(x,y)&=\sup\{\mu^{-1/2m}\phi(y)-\mu^{-1/2m}\phi(x):
    \phi\in\mathcurl{E}_{m,\beta}\} \notag\\
    &\le\sup\{\psi(y)-\psi(x):\psi\in\mathcurl{F}_{a,M}\} \notag\\
    &=d_{a,M}(x,y).
  \end{align}
  For $t/d_g(x,y)\le(1-\sqrt{K/\beta r})\mu^{-1/2m}T_{\delta,M}=:T_\epsilon$ 
  we see, using equation~\eqref{eqn:relbarb} and Theorem~\ref{thm:metricmain}, 
  that
  \begin{align*}
    t/d_{a,M}(x,y)&\le\mu^{1/2m}t/d_{m,\beta}(x,y) \\
    &\le\mu^{1/2m}(1-\sqrt{K/\beta r})^{-1}t/d_g(x,y) \\
    &\le T_{\delta,M},
  \end{align*}
  and so using~\eqref{eqn:barbs},
  \begin{align*}
    \lefteqn{|K(t,x,y)|}& \\
    &\le c_\delta t^{N/2m}\exp[-(\sigma_m-\delta)d_{a,M}(x,y)^{2m/(2m
    -1)}t^{-1/(2m-1)}] \\
    &\le c_\delta t^{N/2m}\exp[-(\sigma_m-\delta)\mu^{-1/(2m-1)}d_{m,\beta}(x,y
    )^{2m/(2m-1)}t^{-1/(2m-1)}] \\
    &\le c_\delta t^{N/2m}\exp[-(\sigma_m-\delta)\mu^{-1/(2m-1)}\left(1-\sqrt{
    \frac{K}{\beta r}}\right)d_g(x,y)^{2m/(2m-1)}t^{-1/(2m-1)}] \\
    &\le c_\epsilon t^{N/2m}\exp[-(\sigma_m-\epsilon)\mu^{-1/(2m-1)}d_g(x,y)^{2
    m/(2m-1)}t^{-1/(2m-1)}].
  \end{align*}
  Using Lemma~\ref{thm:conequiv} this is equivalent to the 
  bound~\eqref{eqn:bestcon}.\finbox
\end{proof}

\section{A COMPARISON OF RIEMANNIAN-TYPE METRICS WITH THE STANDARD RIEMANNIAN 
  METRIC FOR A HIGHLY NON-CONVEX REGION}
\label{sec:metrics}

The purpose of this section is to prove a geometrical result concerning 
metrics, which can be applied in Section~\ref{sec:kerbd}. The Riemannian-type 
metrics of Definition~\ref{def:rietyp} are used in Section~\ref{sec:kerbd} to 
express certain Gaussian heat kernel estimates. It is beneficial to compare 
these metrics for different values of $\beta$ in order to convert the bound 
into terms of the standard Riemannian metric on a region. The comparison is 
particularly interesting when one notes that $d_{m,0}=d_0$ is the standard 
Euclidean metric
\begin{equation}
  d_0(x,y)=|y-x|,
\end{equation}
and $d_{m,\infty}$ is the standard Riemannian metric $d_g$, defined in 
equation~\eqref{eqn:riemann}.
If $\Omega$ is convex then all the above metrics are identical. For non-convex 
regions however this is not the case, and a useful comparison is non-trivial. 
Clearly $d_{m,\beta}$ is an increasing function of $\beta$.

Let $\Omega$ be a region in $\real^N$, whose boundary $\partial\Omega$ is 
$C^2$, with radii of curvature uniformly bounded below. We shall prove that for
$\beta\geq4K/r$,
\begin{equation}\label{eqn:metricbd}
  \left(1-\sqrt{\frac{K}{\beta r}}\right)d_g\leq d_{m,\beta}\leq d_g
\end{equation}
where $r>0$ is the greatest lower bound of the radii of curvature of the
boundary of $\Omega$, and $K$ depends only on $m$ and $N$. This result is also 
valid for locally Euclidean Riemannian manifolds. See Note~\ref{not:loceucl}.


We now develop tools for local representation of the surface 
$\partial\Omega$ of $\Omega$.

\begin{definition}
  We say that a region $\Omega\subseteq\real^N$ has $C^2$ boundary $\partial
  \Omega$ if $\partial\Omega=\partial\overline\Omega$ and if for each point 
  $p\in\partial\Omega$, there exist a set $U=U(p)$ open in $\real^N$ and 
  containing $p$, a local coordinate system $\bd{y}=(y_1,
  \dots,y_{N-1})$ and $y_N$, with $(\bd{y},y_N)=(\bd{0},0)$ at $p$, and a 
  function $h=h(.,p)$ such that $\partial\Omega\cap U$ has a representation
  \begin{equation}
    y_N=h(\bd{y})\qquad\bd{y}\in G,\quad h\in C^2(\overline G),
  \end{equation}
  where $G=G(p)$ is open in $\real^{N-1}$ and convex.
\end{definition}

The surface $\partial\Omega$ within $U$ may equivalently be represented by the non-degenerate,
bijective $C^2$ map $\sigma:G\rightarrow\partial\Omega\cap U$ defined in local coordinates by
\begin{equation}
  \underline\sigma(\bd{y})=(\bd{y},h(\bd{y})).
\end{equation}%
\sngldiagm{%
  \renewcommand\afigwidth{45}%
  \renewcommand\afigheight{32}}{%
  \fig{repbdy}{}{%
    \put(88,56){$p$}%
    \put(154,116){$U$}%
    \put(70,96){$\Omega$}%
    \put(38,86){$\bd{x}$}%
    \put(-4,118){$x_N$}%
    \put(200,100){$\partial\Omega$}%
    \put(142,22){\begin{rotate}{-29}$G$\end{rotate}}%
    \put(112,38){\begin{rotate}{-29}$\bd{y}$\end{rotate}}%
    \put(116,90){\begin{rotate}{-29}$y_N=h(\bd{y})$\end{rotate}}}%
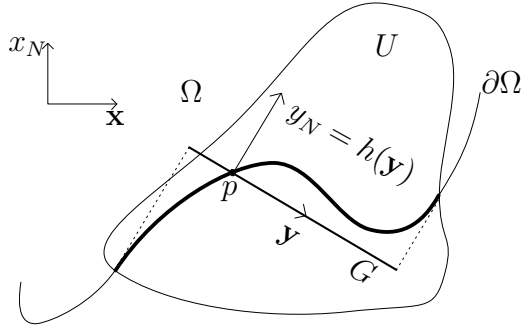
    {Representation of the boundary}}

\begin{thm}{Notation}
  \begin{enumerate}
    \item
      The matrix derivative of a function $f:\real^{N-1}\rightarrow\real^M$ is 
      defined by
      \begin{equation}
        \der[f]:=
        \begin{pmatrix}
          \dev{f}{1}{1}   & \dev{f}{2}{1} & \dots  & \dev{f}{M}{1} \\
          \vdots          &               & \ddots & \vdots \\
          \dev{f}{1}{N-1} &       \hdotsfor{2}     & \dev{f}{M}{N-1} \\
        \end{pmatrix}
      \end{equation}
      where $f_i:=\pi_i\circ f\ (i=1,\dots,M)$ are the coordinate functions of 
      $f$;
    \item
      The second derivative of a function $h:\real^{N-1}\rightarrow\real$ is 
      the matrix defined by
      \begin{equation}
        \sder:=\der[\phantom{\bd{y}}]\left(\der\tr\right)=
        \begin{pmatrix}
          \sdevd{h}{1}     & \dots  & \sdev{h}{1}{N-1} \\
          \vdots           & \ddots & \vdots           \\
          \sdev{h}{N-1}{1} & \dots  & \sdevd{h}{N-1}   \\
        \end{pmatrix}.
      \end{equation}
  \end{enumerate}
\end{thm}

Let the $y_N$-axis point into $\Omega$. The unit normal $n=n(p)$ to the
surface at $p=\sigma(\bd{y})$ is defined in local coordinates by
\begin{equation}
   \underline{n}(\bd{y})=\left(\der\tr,-1\right)\bnumder^{-\frac{1}{2}}.
\end{equation}
The $C^1$ map $\tau:G\times\real\rightarrow\real^N$ defined by
\begin{equation}
   \tau(\bd{y},u)=\sigma(\bd{y})+un(\bd{y})
\end{equation}
is non-degenerate at $(\bd{0},0)$. 

Define $\pi:\partial\Omega\times\real\rightarrow\partial\Omega\times\real$ by
\begin{alignat}{2}
  \pi(p,u)  &= (p,0)\qquad   && p\in\partial\Omega,\quad u\in\real, \\
  \intertext{and $\rho:\partial\Omega\times\real\rightarrow\real^N$ by}
  \rho(p,u) &= p+un(p)\qquad && p\in\partial\Omega,\quad u\in\real.
\end{alignat}

\begin{definition}
  For $\delta>0$ we define the $\delta$-neighbourhoods of $\partial\Omega$ and 
  $\Omega$ by
  \begin{equation}
    (\partial\Omega)_\delta:=\{z\in\real^N:d(z,\partial\Omega)<\delta\}
  \end{equation}
  and
  \begin{equation}
    \Omega_\delta:=\{z\in\real^N:d(z,\Omega)<\delta\}.
  \end{equation}
\end{definition}

\begin{thm}{Proposition}\label{thm:nbdy}
  The $\delta$-neighbourhood $(\partial\Omega)_\delta$ of the boundary is the 
  image of $\partial\Omega\times(-\delta,\delta)$ under the map $\rho$. 
  Similarly $\Omega_\delta=\Omega\cup\rho(\partial\Omega\times[0,\delta))$.
\end{thm}

\begin{proof}
  Suppose that $z=\rho(p,u)=p+un(p)$ for some $p\in
  \partial\Omega$ and some $-\delta<u<\delta$. Then $|z-p|=|u|<\delta$ so 
  $d(z,\partial\Omega)<\delta$.

  Conversely, suppose that $0\leq d:=d(z,\partial\Omega)<\delta$. Let $p\in
  \partial\Omega$ be such that $|z-p|=d$. Using the representation of the 
  surface $\underline\sigma(\bd{y})=(\bd{y},h(\bd{y}))$ in the local 
  coordinate 
  system based at $p$ we see that the function $\bd{y}\mapsto|\underline{z}-(
  \bd{y},h(\bd{y}))|^2$ is minimized at $\bd{y}=\bd{0}$. Thus for $i=1,\dots,N
  -1$,
  \begin{align*}
    0&=\left.\frac{\partial\phantom{y_i}}{\partial y_i}|\underline{z}-(\bd{y},
    h(
    \bd{y}))|^2\right|_{\bd{y}=0} \\
    &=\left.2\left\langle\underline{z}-(\bd{y},h(\bd{y})),-\left(0,\dots,0,1,0,
    \dots,0,\frac{\partial h}{\partial y_i}\right)\right\rangle\right|_{\bd{y}=0
    } \\
    &=-2\left\langle\underline{z}-(\bd{0},0),\left(0,\dots,0,1,0,\dots,0,\left.
    \frac{\partial h}{\partial y_i}\right|_{\bd{y}=0}\right)\right\rangle
  \end{align*}
  so $z-p$ is normal to the surface. The vector $(z-p)/d$ has unit 
  modulus so $(z-p)/d=\pm n(p)$, the sign being dependent on whether 
  $(z-p)/d$ is inward or outward pointing. Thus
  $z=p\pm dn(p)\in\rho(\partial\Omega\times(-\delta,\delta))$.

  The proof that $\Omega_\delta=\Omega\cup\rho(\partial\Omega\times[0,\delta))$ 
  is similar.\finbox
\end{proof}

\begin{thm}{Condition}\label{thm:hyp}
  Let $\Omega$ be a region in $\real^N$ with $C^2$ boundary such that there 
  exists an $r>0$ whereby 
  \begin{equation}\label{eqn:hyp}
    B(p\pm rn(p);r)\cap\partial\Omega=\emptyset
  \end{equation}
  for all $p\in\partial\Omega$.
\end{thm}

This condition is slightly stronger than requiring that the radii of curvature
at points of the boundary are bounded below by $r$. This is done to exclude
certain regions for which the results of this section still hold, but
which require a more technical treatment. See Note~\ref{not:loceucl}.

\begin{thm}{Lemma}\label{thm:equiv}
  Equation~\sleqref{eqn:hyp} in Condition~\slref{thm:hyp} is equivalent to 
  bijectivity of the restriction $\rho:\partial\Omega\times(-r,r)\rightarrow(
  \partial\Omega)_r$.
\end{thm}

\begin{proof}
  Suppose that for some $p\in\partial\Omega$, $B(p+rn(p);r)\cap\partial\Omega
  \not=\emptyset$. By writing $B(p+rn(p);r)=\bigcup_{0<d<r}B_d$ where $B_d:=B
  (p+dn(p);d)$ we see that $B_d\cap\partial\Omega\neq\emptyset$ for some $d<r$. 
  Let $x\in B_d\cap\partial\Omega$. Then $d(x,p+dn(p))<d$ so, by 
  Proposition~\ref{thm:nbdy}, $p+dn(p)\in(\partial\Omega)_d$. Hence
  \[ \rho(p,d)\in\rho(\partial\Omega\times(-d,d)), \]
  and $\rho$ is not injective. The same conclusion is drawn if $B(p-rn(p);r)
  \cap\partial\Omega\not=\emptyset$ for some $p\in\partial\Omega$.
  \par
  Conversely suppose that the restriction of $\rho$ is not injective. Then
  \begin{equation}\label{eqn:inj}
    p+un(p)=q+vn(q)
  \end{equation}
  for some
  \begin{equation}\label{eqn:cond}
    (p,u)\not=(q,v).
  \end{equation}
  Without loss of generality we may assume that 
  $|v|\leq|u|$. Assume also that $0\leq u<r$. Then
  \begin{align}
    |p+rn(p)-q|&=|vn(q)+(r-u)n(p)| \notag\\
    &\leq|v|+r-u \notag\\
    &\leq r.
  \end{align}
  Moreover, using inequality~\eqref{eqn:inj} we see that~\eqref{eqn:cond} 
  implies that $un(p)\not=vn(q)$ and hence the inequality is strict. Thus
  \begin{equation}
    B(p+rn(p);r)\cap\partial\Omega\not=\emptyset.
  \end{equation}
  If we assume that $-r<u\leq0$ then similarly we obtain
  \begin{equation}
    B(p-rn(p);r)\cap\partial\Omega\not=\emptyset. \earlybox
  \end{equation}
\end{proof}

\begin{thm}{Proposition}
  If $\Omega$ is bounded then it satisfies Condition~\slref{thm:hyp}.
\end{thm}

\begin{proof}
  Since $\tau$ is non-degenerate at $(\bd0,0)$, an application of the inverse 
  function theorem shows that $\tau$ it is injective and non-degenerate in an 
  open neighbourhood $N$ of $(\bd{0},0)$. The set $\tau(N)$ is open because 
  $\tau$ is non-degenerate, and since $\tau$ is continuous, $\tau^{-1}(\tau(N)
  \cap U)$ is open. See Figure~\ref{fig:tau_of_n}. If $G$ and $r'>0$ are small 
  enough,
  \begin{equation*}
    G\times(-r',r')\subseteq\tau^{-1}(\tau(N)\cap U)\subseteq N.
  \end{equation*}
  Again, since $\tau$ is non-degenerate in $N$, the coordinate neighbourhood 
  $V:=\tau(G\times(-r',r'))$ is open and so there exists $r''_p>0$ such that 
  $B(p;\frac13r''_p)\subseteq V$. See Figure~\ref{fig:coord}. The collection 
  $\{B(p;\frac13r''_p):p\in\partial\Omega\}$ forms an open covering of 
  $\partial\Omega$ and has a finite subcovering $\{B(p_i;\frac 13r''_i)\}_{i=1
  }^m$ by compactness. Let $r=\frac13\min_{i=1,\dots,m}r''_i$.
\dbledgram{%
  \renewcommand\afigwidth{45}%
  \renewcommand\afigheight{35}%
  \renewcommand\bfigwidth{42}%
  \renewcommand\bfigheight{27}}{%
  \fig{tau_of_n}{\label{fig:tau_of_n}}{%
    \put(154,126){$\partial\Omega$}%
    \put(64,104){$p$}%
    \put(108,32){$\tau(N)$}}%
    {$\tau(N)$ for a typical region $\Omega$.}}{%
  \fig{coord}{\label{fig:coord}}{%
    \put(86,72){$p$}%
    \put(202,100){$\partial\Omega$}%
    \put(112,14){$\tau(N)\cap U$}%
    \put(-4,98){$\tau(G\times(-r',r'))$}}%
    {Construction of\par$G\times(-r',r')$.}}%

  This construction has been chosen so that if $p\in\partial\Omega$ then 
  $p\in B(p_i;\frac 13r''_i)$ for some $i$ and then $B(p;2r)\subseteq B(p_i;r''
  _i)\subseteq V_i$.

  We shall now prove that $\rho$ is injective. For suppose otherwise, then 
  $\rho(p,u)=\rho(q,v)$ where $p,q\in\partial\Omega$, $-r\leq u,v<r$. By the
  construction above, $B(p;2r)\subseteq V_i$ for some $i$. Since $p+un(p)=q+vn
  (q)$ we see that $|p-q|\leq|u|+|v|<2r$ so $q\in V_i$. Thus $p,q\in V_i\cap
  \partial\Omega$ and so
  \begin{gather*}
    p=\tau_i(\bd{x},0)=\sigma_i(\bd{x}) \\
    q=\tau_i(\bd{y},0)=\sigma_i(\bd{y})
  \end{gather*}
  for some $\bd{x},\bd{y}\in G_i$. Now $\tau_i(\bd{x},u)=\rho(p,u)=\rho(q,v)=
  \tau_i(\bd{y},v)$, and since $\tau_i$ is injective we see that $\bd{x}=\bd{y}$ 
  and $u=v$. Moreover, $p=\sigma_i(\bd{x})=\sigma_i(\bd{y})=q$.\finbox
\end{proof}


From this point onwards we shall assume that all regions satisfy
Condition~\ref{thm:hyp}.

\begin{definition}
  We say that a real symmetric $(N-1)\times(N-1)$ matrix $A$ is non-negative, 
  and write $A\geq0$ if $\bd{a}A\bd{a}\tr\geq0$ for all $\bd{a}\in\real^{N-1}$.
\end{definition}

\begin{thm}{Proposition}\label{thm:spheretouch}
  Let $p\in\partial\Omega$ and let $h:G\rightarrow\real$ be the representation 
  of the surface $\partial\Omega$ in the local coordinate system $(\bd{y},y_N)$ 
  based about $p$. Then
  \begin{align}
    \textup{(i)}  & \quad\matder\geq rn_N\sder; \\
    \textup{(ii)} & \quad\matder\geq -rn_N\sder,
  \end{align}
  where
  \begin{equation}
    n_N:=-\bnumder^{-1/2}
  \end{equation}
  is the $N$-th entry of the unit normal $\underline{n}(\bd{y})$ in local coordinates.
\end{thm}

\begin{proof}
  Let $\underline{\tilde\sigma}(\bd{y})=(\bd{y},g(\bd{y}))$ represent, in local 
  coordinates, the surface of the sphere, radius $r$ touching the surface $(\bd
  {y},h(\bd{y}))$ at the point $\underline{q}=(\bd{y_0},h(\bd{y_0}))$ where the 
  unit normal is $\underline{n}=(\bd{n},n_N)$. Then
  \begin{equation}
    (\underline{\tilde\sigma}(\bd{y})-\underline{q})(\underline{\tilde\sigma}(
    \bd{y})-2r\underline{n}-\underline{q})\tr=0.
  \end{equation}
  Differentiating with respect to $\bd{y}=(y_1,\dots,y_{N-1})$ we see that
  \begin{align}
    \bd{0}\tr&=\left(I_{N-1}\left|\der[g]\right.\right)(\underline{\tilde\sigma
    }(\bd{y})-2r\underline{n}-\underline{q})\tr+\left(I_{N-1}\left|\der[g]
    \right.\right)(\underline{\tilde\sigma}(\bd{y})-\underline{q})\tr \notag\\
    &=2\left( I_{N-1}\left|\der[g]\right.\right)(\underline{\tilde\sigma}(\bd
    y)-r\underline{n}-\underline{q})\tr\notag\\
    &=2(\bd{y}\tr-r\bd{n}\tr-\bd{y_0}\tr)+2(g-rn_N-h(\bd{y_0}))\der[g].
  \end{align}
  Differentiating the transpose,
  \begin{equation}
    \bd{0}=\matder[g]+(g-rn_N-h(\bd{y_0}))\sder[g].
  \end{equation}
  Since the sphere touches the surface at $\bd{y}=\bd{y_0}$, we see that 
  \begin{equation}
    g(\bd{y_0})=h(\bd{y_0}),\quad\der[g](\bd{y_0})=\der(\bd{y_0}),\quad
    \text{and}\quad\sder[g](\bd{y_0})\leq\sder(\bd{y_0}).
  \end{equation}
  Thus, at $\bd{y}=\bd{y_0}$,
  \begin{equation}
    \matder=\matder[g]=rn_N\sder[g]\geq rn_N\sder.
  \end{equation}
  This result holds for all $\bd{y_0}\in G$. The proof of part (ii) uses the 
  fact that a ball of radius $r$ fits inside the region.\finbox
\end{proof}

For $\bd{a}\in\real^{N-1}$ define
\begin{equation}
  \alpha=\frac{\displaystyle{\bnumder^{-1/2}\bd{a}\sder\bd{a}\tr}}{
  \displaystyle{\bd{a}\bmatder\bd{a}\tr}},\qquad
  \beta=\frac{\displaystyle{\bd{a}\der[\underline{n}]\der[\underline{n}]\tr\bd{
  a}\tr}}{\displaystyle{\bd{a}\bmatder\bd{a}\tr}},
\end{equation}
and let $A$ be the real symmetric matrix
\begin{equation}
  A=\matder.
\end{equation}

\begin{thm}{Lemma}\label{thm:matrices}
  \begin{align*}
    \textup{(i)}&\quad\der[\underline{n}]\der[\underline{n}]\tr=\bnumder^{-1}
    \sder A^{-1}\sder\mathrm;\\
    \textup{(ii)}&\quad\alpha^2,\beta\in\left[0,1/r^2\right]\mathrm; \\
    \textup{(iii)}&\quad1+2\alpha r+\beta r^2\geq0\text{ and }1-2\alpha r+\beta
    r^2\geq0\mathrm; \\
    \textup{(iv)}&\quad\left(1-\frac{|u|}r\right)^2\leq1+2\alpha u+\beta u^2
    \leq\left(1+\frac{|u|}r\right)^2\text{ for all }u\in[-r,r].
  \end{align*}
\end{thm}

\begin{proof}
  (i) By differentiating the expression for the normal, we see that
  \begin{align*}
    \der[\underline{n}]&=\bnumder^{-3/2}\sder\left[\bnumder\left(I_{N-1}|\bd{0}
    \tr\right)-\der\left(\der\tr,-1\right)\right] \\
    &=\bnumder^{-1/2}\sder\left[\left(I_{N-1}|\bd{0}\tr\right)-\bnumder^
    {-1}\der\left(\der\tr,-1\right)\right].
  \end{align*}
  Hence
  \begin{align*}
  \der[\underline{n}]\der[\underline{n}]\tr&=\bnumder^{-1}\sder\left[I_{N-1}-2
  \bnumder^{-1}\der\der\tr\right. \\
  &\quad\left.+\bnumder^{-2}\der\bnumder\der\tr\right]\sder \\
  &=\bnumder^{-1}\sder\left[I_{N-1}-\der\der\tr\bnumder^{-1}\right]\sder \\
  &=\bnumder^{-1}\sder A^{-1}\sder.
  \end{align*}
  (ii) Proposition~\ref{thm:spheretouch} immediately implies that 
  \begin{equation*}
    -\frac1r\leq\alpha\leq\frac1r.
  \end{equation*}
  Since $A$ is a positive definite matrix, we may pre- and post-multiply the 
  results of Proposition~\ref{thm:spheretouch} by $A^{-1/2}$ to get
  \begin{gather}
    I_{N-1}+rn_NA^{-1/2}\sder A^{-1/2}\geq0; \label{eqn:curvI}\\
    I_{N-1}-rn_NA^{-1/2}\sder A^{-1/2}\geq0. \label{eqn:curvII}
  \end{gather}
  The matrices in inequalities~\eqref{eqn:curvI} and~\eqref{eqn:curvII} 
  commute, so their product is also positive. Hence
  \begin{equation}
    I_{N-1}-{r}^2n_N^2A^{-1/2}\sder A^{-1}\sder A^{-1/2}\geq0.
  \end{equation}
  By pre- and post-multiplying by $A^{1/2}$, we see that
  \begin{equation}
    \bmatder-{r}^2n_N^2\sder A^{-1}\sder\geq0,
  \end{equation}
  so using part (i),
  \begin{align*}
    r^2\bd{a}\der[\underline{n}]\der[\underline{n}]\tr\bd{a}\tr&=r^2n_N^2\bd{a}
    \sder A^{-1}\sder\bd{a}\tr\\
    &\le\bd{a}\bmatder\bd{a}\tr.
  \end{align*}
  Thus $\beta\leq1/r^2$.

  (iii) Squaring the left hand side matrix of inequality~\eqref{eqn:curvI} will 
  yield a positive matrix
  \begin{equation*}
    I_{N-1}+2rn_NA^{-1/2}\sder A^{-1/2}+{r}^2n_N^2A^{-1/2}\sder 
    A^{-1}\sder A^{-1/2}.
  \end{equation*}
  Pre- and post-multiplying by $A^{1/2}$, we see that 
  \begin{equation*}
    \bd{a}\bmatder\bd{a}\tr+2rn_N\bd{a}\sder\bd{a}\tr+{r}^2n_N^2\bd{a}\sder A^{
    -1}\sder\bd{a}\tr\geq0
  \end{equation*}
  for all $\bd{a}\in\real^{N-1}$. Thus $1+2\alpha r+\beta r^2\geq0$. Squaring 
  the left hand side of inequality~\eqref{eqn:curvII}, we see that $1-2\alpha 
  r+\beta r^2\geq0$.

  (iv) By part (iii), the polynomial $1+2\alpha u+\beta u^2$ dominates 
  \[ 1-2\frac u{r}+\frac{u^2}{{r}^2} \]
  at $u=0,r$. Moreover, since $\beta\leq1/r^2$, this holds true for 
  all $u\in[0,r]$. Thus
  \begin{equation*}
    1+2\alpha u+\beta u^2\geq\left(1-\frac u{r}\right)^2
  \end{equation*}
  for all $u\in[0,r]$. Similarly,
  \begin{equation*}
    1+2\alpha u+\beta u^2\geq\left(1+\frac u{r}\right)^2
  \end{equation*}
  for all $u\in[-r,0]$. Also, by part (ii),
  \begin{equation*}
    1+2\alpha u+\beta u^2\leq 1+2\frac{|u|}{r}+\frac{|u|^2}{{r}^2}=\left(1+
    \frac{|u|}{r}\right)^2. \lowearlybox
  \end{equation*}
\end{proof}

Let $P:G\times(-r,r)\rightarrow G\times(-r,r)$ be the projection defined by 
$P(\bd{y},u)=(\bd{y},0)$. Let
\begin{equation}
  \tau'(.,.):G\times(-r,r)\rightarrow M_N(\real)
\end{equation}
denote the Jacobian matrix 
\[ \frac{\dup\tau}{\dup(\bd{y},u)} \]
of $\tau$.

\begin{thm}{Proposition} \label{thm:cvatur}
  Let $(\bd{y},u)\in G\times(-r,r)$ and let $v=(\bd{a},a_N)\in\real^N$. Then
  \begin{equation}
    \left(1-\frac{|u|}{r}\right)|P(v)(\tau'\circ P)(\bd{y},u)|
    \leq\left|v\tau'(\bd{y},u)\right|.
  \end{equation}
  An alternative formulation of the above in local coordinates is
  \begin{equation}
    \bd{a}\bmatder\bd{a}\tr\left(1-\frac{|u|}{r}\right)^2
    \leq|v\underline{\tau}'(\bd{y},u)|^2.
  \end{equation}
\end{thm}

\begin{proof}
  \begin{align*}
    |v\,\underline{\tau}'(\bd{y},u)|^2 &= \left|v\left[
    \left(
    \begin{array}[c]{c|c}
      I_{N-1}&\der \\[2pt] \hline
      \multicolumn{2}{c}{\underline{n}(\bd{y})}
    \end{array}
    \right)+u\left(
    \begin{array}[c]{c}
      \der[\underline{n}] \\[2pt] \hline
      \bd0 
    \end{array}
    \right)\right]\right|^2 \\
    &=v\left(
    \begin{array}[c]{c|c}
      \matder&\bd0\tr \\[2pt] \hline
      \bd0&1\phantom{\tr}
    \end{array}
    \right)v\tr \\
    &\qquad\qquad+2uv\bnumder^{-1/2}\left(
    \begin{array}[c]{c|c}
      \sder&\bd0\tr \\[2pt] \hline
      \bd0&0\phantom{\tr}
    \end{array}
    \right)v\tr \\
    &\qquad\qquad+u^2v\left(
    \begin{array}[c]{c|c}
      \der[\underline{n}]\der[\underline{n}]\tr&\bd{0}\tr \\[2pt] \hline
      \bd0&0\phantom{\tr}
    \end{array}
    \right)v\tr \\
    &=\bd{a}\bmatder\bd{a}\tr+a_N^2 \\
    &\qquad+2u\bnumder^{-1/2}\bd{a}\sder\bd{a}\tr+u^2\bd{a}\der[\underline{n}]\der[\underline{n}]
    \tr\bd{a}\tr \\
    &=\bd{a}\bmatder\bd{a}\tr[1+2\alpha u+\beta u^2]+a_N^2
  \end{align*}
  The proposition follows by applying Lemma~\ref{thm:matrices}.\finbox
\end{proof}


\sngldiagm{%
  \renewcommand\afigwidth{50}%
  \renewcommand\afigheight{44}}{%
  \fig{nurhopi}{}{%
  \put(154,0){$\nu_1$}%
  \put(98,72){$\rho$}%
  \put(122,96){$\rho^{-1}$}%
  \put(60,110){$\pi$}%
  \put(200,92){$(\partial\Omega)_r$}%
  \put(162,150){$\partial\Omega\times(-r,r)$}}%
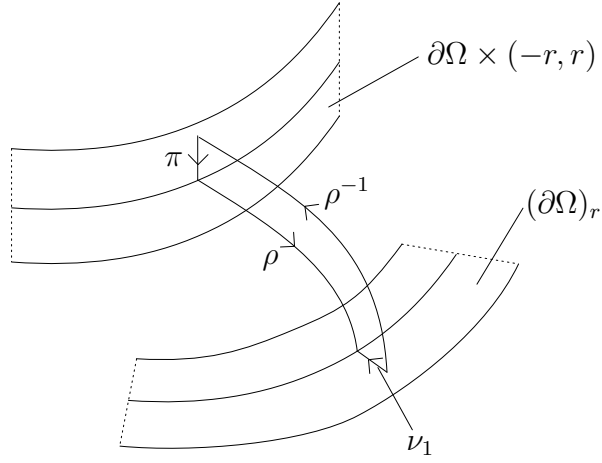
  {Construction Of $\nu_1$}}%

We shall now remove references to local parametrisations of the 
boundary. If a point $z\in\real^N$ is within a distance $r$ of $\partial\Omega$
then there is a unique nearest point of
$\partial\Omega$ to $z$. We shall define maps $\nu_1$, $\nu_2$ which formalize 
the notion of a nearest point to $\partial\Omega$ and $\overline\Omega$ 
respectively. Define $\nu_1:(\partial\Omega)_r\rightarrow\partial\Omega$ by
\begin{equation}
  \nu_1:=\rho\pi\rho^{-1}.
\end{equation}
This is well defined due to Lemma~\ref{thm:equiv}. Define $\nu_2:\Omega_r\rightarrow\overline\Omega$ by
\begin{equation}
  \nu_2(\omega)=\begin{cases}
     \omega,                       & \omega\in\Omega \\
     \nu_1(\omega), & \omega\in\rho(\partial\Omega\times[0,r)).
   \end{cases}
\end{equation}
Note that $\rho(\partial\Omega\times[0,r))\cap\Omega=\emptyset$ because of the construction of $\rho$.

Let $x\in\Omega$ and define $d_x:\Omega_{r}\rightarrow\real_+$ by
\begin{equation}
  d_x(z):=d_g(x,\nu_2 z)
\end{equation}
where we recall that the standard Riemannian metric $d_g:\overline\Omega^2
\rightarrow\real_+$ is defined by
\begin{equation}
  d_g(x,y):=\inf\{l(\gamma):\gamma(0)=x,\gamma(1)=y,\gamma\subseteq\overline
  \Omega,\gamma\text{ is cts and piecewise }C^1\}.
\end{equation}%
\sngldiagm{%
  \renewcommand\afigwidth{50}%
  \renewcommand\afigheight{39}}{%
  \fig{distfn}{}{%
  \put(30,100){$x$}%
  \put(200,60){$\Omega_r$}%
  \put(186,118){$\partial\Omega$}}%
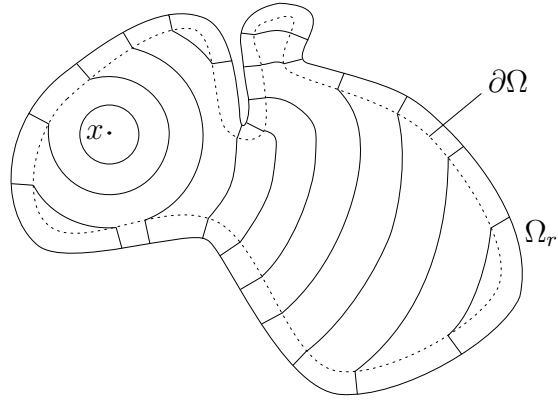
  {Contour sketch of $d_x$ for a typical region $\Omega$}}

\begin{thm}{Construction}\label{thm:constr}
  Let $0\leq\delta\leq r$ be fixed. Let $\gamma_1$ satisfy either of the 
  following conditions, and define $\gamma_2$ accordingly.
  \begin{enumerate}
    \item Let $\gamma_1:[0,1]\rightarrow(\partial\Omega)_\delta$ be a $C^1$ 
          curve such that $\gamma_1\not\subseteq\Omega$. Let 
          \begin{equation}
            T_0=\inf\{t\in[0,1]:\gamma_1(t)\not\in\Omega\}
          \end{equation}
          and let 
          \begin{equation}
            T_1=\sup\{t\in[0,1]:\gamma_1(t)\not\in\Omega\}.
          \end{equation}
          Define $\gamma_2:[0,1]\rightarrow\overline\Omega$ by
          \begin{equation}
            \gamma_2(t)=\begin{cases}
              \nu_1\gamma_1(t),\quad & t\in(T_0,T_1) \\
              \gamma_1(t),                          & \text{otherwise\textup;}
              \end{cases}
          \end{equation}
    \item Let $\gamma_1:[0,1]\rightarrow\Omega$ be a $C^1$ curve. Then let 
          $\gamma_2=\gamma_1$.
  \end{enumerate}
\end{thm}

\begin{thm}{Lemma}\label{thm:newcurv}
  In both cases of Construction~\slref{thm:constr}, $\gamma_2$ is a 
  piecewise $C^1$ curve such that $\gamma_2\subseteq\overline\Omega$, $\gamma
  _2(0)=\nu_2[\gamma_1(0)]$, $\gamma_2(1)=\nu_2[\gamma_1(1)]$ and
  \begin{equation}
    \left(1-\frac\delta r\right)l(\gamma_2)\leq l(\gamma_1).
  \end{equation}
\end{thm}

\begin{proof}
  We shall only prove the lemma for case (i), as case (ii) is trivial. For 
  $t\in(T_0,T_1)$ let $p\in\partial\Omega$ be the nearest point of
  $\partial\Omega$ to $\gamma_1(t)$ and let $(\bd{y},h(\bd{y}))$ represent the 
  surface $\partial\Omega$ in the local coordinate system based at $p$. Since 
  $\gamma_1$ is continuous, there exists an open interval $I\subseteq(T_0,T_1)
  $, containing $t$ such that $\gamma_1(I)\subseteq\tau(G\times(-r,r))$. Let 
  $\gamma_3,\gamma_4$ be paths defined, for $s\in I$, by%
  \sngldiagm{%
    \renewcommand\afigwidth{90}%
    \renewcommand\afigheight{20}}{%
    \fig{constr}{}{%
    \put(63,36){$\gamma_1$}%
    \put(52,23){$\gamma_2$}%
    \put(182,34){$\gamma_3$}%
    \put(171,18){$\gamma_4$}%
    \put(97,31){$\tau$}%
    \put(95,10){$\tau^{-1}$}%
    \put(20,5){$\Omega$}}%
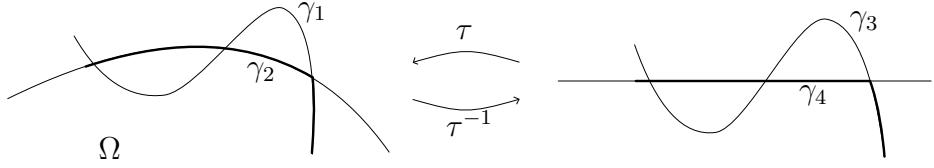
    {Construction of the paths $\gamma_1$, $\gamma_2$, $\gamma_3$, and%
    $\gamma_4$}}

  \begin{equation}
    \gamma_3:=\tau^{-1}\gamma_1\qquad\gamma_4:=P\gamma_3.
  \end{equation}
  Note that $\gamma_2|_I=\tau\gamma_4|_I$ and that $\gamma_2,\gamma_3,\gamma_4
  $ are $C^1$ on $I$ because $\tau$ is $C^1$. By writing $\gamma_3(s)=(\bd{a}(
  s),a_N(s))$ we see that $|a_N(t)|=d(\gamma_1(t),\partial\Omega)<\delta$. 
  Thus by Proposition~\ref{thm:cvatur},
  \begin{align*}
    \left(1-\frac\delta r\right)|\gamma_2'(t)|
    &<    \left(1-\frac{|a_N(t)|}{r}\right)|\gamma_4'(t)\tau'(\gamma_4(t))| \\
    &=    \left(1-\frac{|a_N(t)|}{r}\right)|P\gamma_3'(t)(\tau'\circ P)(\gamma_3
    (t))| \\
    &\leq |\gamma_3'(t)\tau'(\gamma_3(t))| \\
    &=    |\gamma_1'(t)|.
  \end{align*}
  Integration yields the result.\finbox
\end{proof}

\begin{thm}{Corollary}\label{thm:dist}
  \begin{enumerate}
    \item
      Let $z_1\in\Omega_{\delta}$, where $0<\delta\leq r$. Then
      \begin{equation}
        \limsup_{z_2\rightarrow z_1}\frac{|d_x(z_2)-d_x(z_1)|}{|z_2-z_1|}\left
         (1-\frac{\delta}r\right)\leq1\mathrm;
      \end{equation}
    \item
      Let $x\in\Omega$ and let $z\in B(x;\delta)$, where $0<\delta\leq r$. 
      Then
      \begin{equation}
        d_g(x,\nu_2 z)\le\left(\frac1\delta-\frac1r\right)^{-1}\mathrm;
      \end{equation}
    \item
      Let $x,y\in\Omega$ be such that $d_g(x,y)<2r$. Then
      \begin{equation}
        |y-x|\geq\frac{2rd_g(x,y)}{2r+d_g(x,y)}.
      \end{equation}
  \end{enumerate}
\end{thm}

\begin{proof}
  (i) Let $z_2\in B(z_1;\delta-d(z_1,\overline\Omega))$ and let $\gamma_1$ be 
      the straight line joining $z_1$ and $z_2$. Then $\gamma_1$ satisfies
      one of the two conditions in Construction~\ref{thm:constr}. Accordingly,
      \begin{align*}
        |d_x(z_2)-d_x(z_1)| &= |d_g(x,\nu_2 z_2)-d_g(x,\nu_2 z_1)| \\
        &\leq d_g(\nu_2 z_1,\nu_2 z_2) \\
        &\leq l(\gamma_2) \\
        &\leq \left(1-\frac\delta r\right)^{-1}l(\gamma_1) \\
        &=    \left(1-\frac\delta r\right)^{-1}|z_2-z_1|.
      \end{align*}
 (ii) Let $\gamma_1$ be the straight line joining $x$ and $z$. Then $\gamma_1$ 
      satisfies one of the conditions in Construction~\ref{thm:constr}. Hence, 
      by Lemma~\ref{thm:newcurv}
      \begin{equation*}
        d_g(x,\nu_2 z)\left(1-\frac\delta r\right)\leq\left(1-\frac\delta r
        \right)l(\gamma_2)\leq l(\gamma_1)=\delta.
      \end{equation*}
(iii) Let $\gamma_1$ be the straight line joining $x$ and $y$. Then setting 
      $\delta=\frac12|y-x|<r$, we see that $\gamma_1$ satisfies one of the
      conditions in Construction~\ref{thm:constr}. Hence, by 
      Lemma~\ref{thm:newcurv}
      \[ \left(1-\frac{\frac12|y-x|}r\right)d_g(x,y)\leq|y-x|. \]
      Rearranging this gives the desired inequality.\finbox
\end{proof}


Let $k_\delta$ be approximate identities, defined as follows:

\begin{definition}
  Let $B(0;1)$ denote the unit ball in $\real^N$, and let $k:B(0;1)\rightarrow
  \real$ be smooth, non-negative and have unit integral. For $\delta>0$ define 
  $k_\delta:B(0;\delta)\rightarrow\real$ by 
  \begin{equation}
    k_\delta(z)=\delta^{-N}k\left(\frac z\delta\right).
  \end{equation}
\end{definition}

Let $x\in\Omega$ and let $\beta>K/r$ where 
\begin{equation}
  K=K_{m,N,k}:=\sup_{1\le|j|\le m-1}\bigg(\int_{B(0;1)}|D^jk|\bigg)^{1/|j|}.
\end{equation}

Define $f_{m,\beta,x}:\Omega\rightarrow\real$ by
\begin{align}
  f_{m,\beta,x}(y):=& \left(1-\frac{K}{\beta r}\right)\int_{B(y;K/\beta)}d_x
  (z)k_{K/\beta}(z-y)\dup^Nz \notag\\
  =& \left(1-\frac{K}{\beta r}\right)\int_{B(0;K/\beta)}d_x(z+y)k_
  {K/\beta}(z)\dup^Nz.
\end{align}

\begin{thm}{Lemma}\label{thm:fisin}
  The functions $f_{m,\beta,x}$ belong to the class $\mathcurl{E}_{m,\beta}$ 
  (see Definition~\slref{def:rietyp}).
\end{thm}

\begin{proof}
  Suppose that $y_1,y_2\in\Omega$ and that $|j|$ is a non-negative multi-index 
  such that $1\le|j|\le m-1$. Then by Corollary~\ref{thm:dist} (i),
  \begin{align*}
    \lefteqn{\limsup_{y_2\rightarrow y_1}\frac{|D^jf_{m,\beta,x}(y_2)-D^jf_
    {x,\beta}(y_1)|}{|y_2-y_1|}} \quad & \\
    &=\limsup_{y_2\rightarrow y_1}\frac{|\left(1-\frac{K}{\beta r}\right)
    \int_{B(0;K/\beta)}\{d_x(z+y_2)-d_x(z+y_1)\}D^jk_{K/\beta}
    (z)\dup^Nz|}{|y_2-y_1|} \\
    &\leq\limsup_{y_2\rightarrow y_1}\int_{B(0;K/\beta)}\frac{|d_x(z+y_2
    )-d_x(z+y_1)|}{|(z+y_2)-(z+y_1)|}\left(1-\frac{K}{\beta r}\right)|D^jk_
    {K/\beta}(z)|\dup^Nz \\
    &=\int_{B(0;K/\beta)}\left[\limsup_{y_2\rightarrow y_1}\frac{|d_x(z+y
    _2)-d_x(z+y_1)|}{|(z+y_2)-(z+y_1)|}\left(1-\frac{K}{\beta r}\right)\right
    ]|D^jk_{K/\beta}(z)|\dup^Nz \\
    &\leq\int_{B(0;K/\beta)}|D^jk_{K/\beta}(z)|\dup^Nz=K^{-|j|}\beta^{|j|
    }\int_{B(0;1)}|D^jk(z)|\dup^Nz\le\beta^{|j|}.
  \end{align*}
  Since $f_{m,\beta,x}$ is $C^\infty$ and satisfies the above  
  inequality we see that $|\nabla f_{m,\beta,x}(y).\bd{u}\tr|\leq1$ for all 
  $\bd{u}\in\real^N$ with $|\bd{u}|=1$, by taking $j=0$, and hence that $|
  \nabla f_{m,\beta,x}(y)|\leq1$. Also, letting $j\le i$ be any non-negative 
  multi-index such that $|j|=|i|-1$,
  \begin{equation*}
    |D^if_{m,\beta,x}(y)|\le\beta^{|i|-1}. \earlybox
  \end{equation*}
\end{proof}

Recall the definition of the Riemannian-type metrics $d_{m,\beta}:\Omega^2
\rightarrow\real_+$
\begin{equation}
   d_{m,\beta}(x,y):=\sup\{\phi(y)-\phi(x):\phi\in\mathcurl{E}_{m,\beta}\}
\end{equation}

\begin{thm}{Lemma}\label{thm:penult}
  Let $\Omega$ satisfy Condition~\slref{thm:hyp}. Let $x,y\in\Omega$ and let 
  $\beta>K/r$. Then
  \begin{equation}
    d_{m,\beta}(x,y)\ge d_g(x,y)\left(1-\frac{K}{\beta r}\right)-\frac{2K}
    \beta.
  \end{equation}
\end{thm}

\begin{proof}
  If $z\in B(x;K/\beta)$ then by Corollary~\ref{thm:dist} (ii),
  \begin{equation*}
    d_x(z)=d_g(x,\nu_{\overline\Omega}z)\le\left(\frac{\beta}{K}-\frac1r
    \right)^{-1},
  \end{equation*}
  so
  \begin{equation*}
    f_{m,\beta,x}(x)=\int_{B(x;K/\beta)}d_x(z)\left(1-\frac{K}{\beta r}
    \right)k_{K/\beta}(z-x)\dup^Nz\leq\frac{K}{\beta}.
  \end{equation*}
  Similarly, if $z\in B(y;K/\beta)$ then
  \begin{equation*}
    d_x(z)=d_g(x,\nu_{\overline\Omega}z)\ge d_g(x,y)-d_g(y,\nu_{\overline\Omega
    }z)\ge d_g(x,y)-\left(\frac{\beta}{K}-\frac1r\right)^{-1},
  \end{equation*}
  so
  \begin{align*}
    f_{m,\beta,x}(y)&=\int_{B(y;K/\beta)}d_x(z)\left(1-\frac{K}{\beta{r}}
    \right)k_{K/\beta}(z-y)\dup^Nz \\
    &\geq d_g(x,y)\left(1-\frac{K}{\beta r}\right)-\frac{K}{\beta}.
  \end{align*}
  Thus using Lemma~\ref{thm:fisin},
  \begin{equation*}
    d_{m,\beta}(x,y)\geq f_{m,\beta,x}(y)-f_{m,\beta,x}(x)\geq d_g(x,y)\left(1-
    \frac{K}{\beta r}\right)-\frac{2K}{\beta}.\lowearlybox
  \end{equation*}
\end{proof}

The dependence of $K$ upon $k$ can be removed by taking the infimum of all 
values of $K_{m,N,k}$, where $k$ is smooth, non-negative, supported in the unit 
ball and has unit integral. For example, $K_{2,N}=N^2$.


\begin{thm}{Theorem}\label{thm:metricmain}
  Let $\Omega$ satisfy Condition~\slref{thm:hyp}. For $\beta\geq4K/r$ we 
  have
  \begin{equation*}
    \left(1-\sqrt{\frac{K}{\beta r}}\right)d_g\leq d_{m,\beta}\leq d_g.
  \end{equation*}
\end{thm}

\begin{proof}
  Let $\gamma\subseteq\overline\Omega$, $\gamma(0)=x$, $\gamma(1)=y$, $\gamma$ 
  is piecewise $C^1$ and let $\phi$ be such that $|\nabla\phi|\leq1$. Then
  \begin{align*}
    |\phi(y)-\phi(x)|&=|\phi(\gamma(1))-\phi(\gamma(0))| \\
    &=\left|\int_0^1\frac{\dup\phantom t}{\dup t}\phi(\gamma(t))\dup t\right| \\
    &=\left|\int_0^1\nabla\phi.\gamma'(t)\dup t\right| \\
    &\leq \int_0^1|\nabla\phi||\gamma'(t)|\dup t \\
    &\leq l(\gamma).
  \end{align*}
  Taking the supremum over all $\phi$ in $\mathcurl{E}_{m,\beta}$, we see that
  \[ d_{m,\beta}(x,y)\leq d_g(x,y). \]

  Let $\beta>4K/r$ and let $\epsilon:=\sqrt{K/\beta r}<1/2$. For large 
  distances $d_g(x,y)$, Lemma \ref{thm:penult} is a useful result. We
  therefore consider two cases:

  (i) Suppose $d_g(x,y)\le2\epsilon r/(1-\epsilon)$. Then we may use 
  Corollary~\ref{thm:dist} (iii) to obtain
  \begin{align*}
    d_{m,\beta}(x,y)&\geq|y-x| \\
    &\ge\frac{2rd_g(x,y)}{2r+d_g(x,y)} \\
    &\ge\frac{2rd_g(x,y)}{2r/(1-\epsilon)} \\
    &=\left(1-\sqrt{\frac{K}{\beta r}}\right)d_g(x,y).
  \end{align*}

  (ii) For $d_g(x,y)\ge2\epsilon r/(1-\epsilon)$ we see that
  \begin{align*}
    d_g(x,y)\left(\epsilon-\frac{K}{r\beta}\right)&\ge\frac{2\epsilon r}{1-
    \epsilon}(\epsilon-\epsilon^2) \\
    &=2\epsilon^2r \\
    &=\frac{2K}\beta.
  \end{align*}
  Hence
  \begin{equation*}
    d_g(x,y)\left(1-\frac{K}{r\beta}\right)-\frac{2K}\beta\ge(1-\epsilon)d_
    g(x,y)
  \end{equation*}
  Using Lemma~\ref{thm:penult},
  \begin{equation*}
    d_{m,\beta}(x,y)\geq\left(1-\sqrt{\frac{K}{\beta r}}\right)d_g(x,y).
  \end{equation*}
  In both cases, we see that 
  \begin{equation*}
    \left(1-\sqrt{\frac{K}{\beta r}}\right)d_g(x,y)\leq d_{m,\beta}(x,y)\leq 
    d_g(x,y). \lowearlybox
  \end{equation*}
\end{proof}

\begin{note}\label{not:loceucl}
  Suppose that $\Omega$ is an $N$-dimensional locally Euclidean Riemannian 
  manifold which posesses a locally injective isometric mapping into $\real^N$ 
  (i.e. $\Omega$ is a covering space of some non-simply connected open subset 
  of $\real^N$). The notion of a Euclidean metric in such manifolds 
  degenerates, and so heat kernel bounds involving such metrics are not useful. 
  Heat kernel bounds involving the Riemannian metric can be found by adapting 
  the result of this section and using the methods of Section~\ref{sec:kerbd}. 
  In fact, by using a covering space, it can be seen that 
  Theorem~\ref{thm:metricmain} is valid where $r$ is the greatest lower bound 
  of the radii of curvature at all points of the boundary.\finbox
\end{note}

\textit{Acknowledgments.} I wish to thank E B Davies for suggesting 
this problem and for his invaluable guidance and support during my research. I 
thank also Gerassimos Barbatis for a number of useful comments. This research 
was funded by an EPSRC studentship.

\bibliographystyle{plain}

\renewcommand\baselinestretch1
\begin{flushright}
  \begin{minipage}{7cm}\centering
    MARK P. OWEN\\
    \em
    Department of Mathematics\\
    King's College London\\
    Strand\\
    London WC2R  2LS\\
    ENGLAND
  \end{minipage}
\end{flushright} 

\end{document}